\documentclass[a4paper,12pt]{book}
\usepackage{geometry}
\usepackage[british]{babel}
\usepackage[utf8]{inputenc}
\usepackage{amsmath}
\usepackage{amsfonts}
\usepackage{multirow}
\usepackage{rotating}
\usepackage{dsfont}
\usepackage{amssymb}
\usepackage{graphicx}
\usepackage{setspace}
\usepackage[arrow, matrix, curve]{xy}
\usepackage{oldgerm}
\usepackage{makeidx}
\usepackage{needspace}
\usepackage{mathrsfs}
\usepackage{tocbibind}
\usepackage{caption}
\usepackage{microtype}

\numberwithin{figure}{section}
\usepackage[hyperref,amsmath,thmmarks,thref]{ntheorem}

\theorembodyfont{}

\usepackage{index}
\newindex{default}{idx}{ind}{Index}
\newcommand{\begriff}[1]{{\index{#1}}\emp{#1}}
\newcommand{\emp}[1]{\textbf{\emph{#1}}}
\newcommand{\ra}{\rightarrow}
\newcommand{\z}{\overline{z}}

\newcommand{\y}{\overline{y}}
\newcommand{\f}{\overline{f}}

\renewcommand{\c}{\overline{c}}

\renewcommand{\ra}{\rightarrow}
\newtheorem{satz}{Theorem}[chapter]
\newtheorem{ko}{Corollary}[chapter]
\newtheorem{bem}{Remark}[chapter]
\newtheorem{defbem}{Definition and remark}[chapter]

\newtheorem{lemma}{Lemma}[chapter]
\newtheorem{definition}{Definition}[chapter]

\theoremstyle{nonumberplain}
\newtheorem{beme}{Remark}

\theoremstyle{nonumberplain}
\theoremsymbol{\ensuremath{\Box}}
\newtheorem{Bew}{Proof}

\DeclareMathOperator{\Real}{Re}
\DeclareMathOperator{\Imag}{Im}
\DeclareMathOperator{\kgV}{lcm}
\DeclareMathOperator{\Aut}{Aut}
\DeclareMathOperator{\Isom}{Isom}
\DeclareMathOperator{\Out}{Out}
\DeclareMathOperator{\Inn}{Inn}

\DeclareMathOperator{\Kern}{Kern}
\DeclareMathOperator{\Diffeo}{Diffeo}
\DeclareMathOperator{\Homo}{Homeo}
\DeclareMathOperator{\Bild}{image}
\DeclareMathOperator{\id}{id}

\DeclareMathOperator{\SL}{SL}
\DeclareMathOperator{\PSL}{PSL}
\DeclareMathOperator{\Stab}{Stab}

\DeclareMathOperator{\GL}{GL}
\DeclareMathOperator{\pr}{pr}
\DeclareMathOperator{\proj}{proj}
\DeclareMathOperator{\SO}{SO}
\DeclareMathOperator{\Aff}{Aff}

\DeclareMathOperator{\re}{ref}
\DeclareMathOperator{\der}{der}

\DeclareMathOperator{\Mod}{Mod}
\DeclareMathOperator{\Deck}{Deck}
\DeclareMathOperator{\ggT}{gcd}
\DeclareMathOperator{\HH}{\mathds H}
\DeclareMathOperator{\R}{\mathds R}
\DeclareMathOperator{\C}{\mathds C}

\DeclareMathOperator{\D}{\mathds D}
\DeclareMathOperator{\Z}{\mathds Z}
\DeclareMathOperator{\N}{\mathds N}
\DeclareMathOperator{\K}{\mathds K}

\DeclareMathOperator{\Harm}{Harm}
\DeclareMathOperator{\pw}{pw}
\DeclareMathOperator{\jac}{jac}

\newcommand{\FakRaum}[2]{
\raisebox{0.7ex}{\ensuremath{#1}}
\ensuremath{\mkern-3mu}\big/\ensuremath{\mkern-3mu}
\raisebox{-0.6ex}{\ensuremath{#2}}}

\title{Teichmüller discs in Schottkyspace}
\date{}
\author{Diego De Filippi\\ Karlsruher Institut für Technologie (KIT)\\
Institut für Algebra und Geometrie}
\begin{document}
\maketitle

\tableofcontents
\begin{spacing}{1.2}

\chapter{Introduction}

The Teichmüller space $\mathcal T_g$ classifies marked Riemann surfaces
of genus $g$. There is more than one approach to Teichmüller space,
e.g.\ by
using quasiconformal self maps or hyperbolic structures on a surface
$X$. It can also be defined by choosing a symplectic set of generators
of the fundamental group~$\pi_1(X)$ or an orientation preserving
diffeomorphism $f\colon X_{\re}\ra X$ as marking of the corresponding Riemann
surface $X$, where $X_{\re}$ is a fixed surface of the same genus. Among
other things, Teichmüller space is interesting because it helps
to understand the moduli space $\mathcal M_g$ that classifies the
isomorphism classes of (non-marked) Riemann surfaces of genus $g$ since we
get $\mathcal M_g$ as quotient space of $\mathcal T_g$ by the action of
the mapping class group $\Mod_g$.\\ One theme in the study of $\mathcal
T_g$ is the investigation of so-called Teichmüller discs. 
These are holomorphic, isometric
embeddings of the complex upper half-plane $\HH$ (respectively the unit disc
$\D$) into $\mathcal T_g$ (with respect to the Teichmüller metric on
$\mathcal T_g$ and the hyperbolic metric in $\HH$). There is one of them
for every point in $\mathcal T_g$ and for every complex direction. In
chapter \ref{grundlagen} we explain how to construct such embeddings.
Strongly related to it is the notion of a flat surface, i.e.\ a surface
whose transition maps are locally of
the form $z\mapsto \pm z+c$ for some constant $c\in\C$. We get such a
flat structure on every compact Riemann surface by means of a quadratic
holomorphic
differential $q$. Varying the flat structure by affine maps
on the charts changes the point in Teichmüller space, and hence we get
a set of points in $\mathcal T_g$ that turn out to be the points of a
Teichmüller disc $\Delta_q$. Moreover, the set of self maps of
the surface that are locally affine with respect to the flat structure
induced by $q$ form a group $\Aff^+(X,q)$ that is isomorphic to the
stabilizer of $\Delta_q$ in $\Mod_g$. The projection of $\Aff^+(X,q)$ to
$\PSL_2(\R)$ is known as projective Veech group
$\overline{\Gamma}(X,q)$. Now, the image of $\Delta_q$ in $\mathcal M_g$
is an algebraic curve if and only if $\overline{\Gamma}(X,q)$ is a
lattice. The Teichmüller discs that satisfy this requirement are subject
of great interest. We call their images in~$\mathcal M_g$ Teichmüller
curves.\\ The present work is concerned with the images of such Teichmüller discs
in the Schottky space $\mathcal S_g$. The latter can be seen as an
intermediate 
covering between $\mathcal T_g$ and $\mathcal M_g$. It classifies pairs
that consist of a Riemann surface and a Schottky covering of it, that is 
Schottky coverings are minimal planar coverings, whose deck
transformation group is free of rank $g$. The properties of Schottky
groups and of the Schottky space are explained and proved at the
beginning of chapter \ref{schottky}. These results are, indeed, already known (a
reference for Schottky groups and the more general class of Kleinian
groups is e.g.\ \cite{klein}), but they are exposed and proved here in a
way that is suitable for our purposes as well as for reasons of
self-containedness. After having dealt with Schottky groups, we explain how $\mathcal S_g$ is a (non
normal) covering between $\mathcal T_g$ and $\mathcal M_g$. We get a map
$s_\alpha\colon \mathcal T_g\ra\mathcal S_g$ which is not canonical since it
depends on a homomorphism $\alpha\colon \pi_g\ra F_g$, where $\pi_g=\langle
a_1,b_1,\dots,a_g,b_g|\Pi_{i=1}^g[a_i,b_i]\rangle$ is the abstract
fundamental group, and $F_g$ is the free group of rank $g$. The map
$s_\alpha$ is the quotient map after the subgroup $\Mod_g(\alpha) \subseteq
\Mod_g$ given by
\[\Mod_g(\alpha)=\{\varphi\in\Mod_g\mid\alpha\circ\tilde{\varphi}\equiv\alpha
  \operatorname{mod} \Inn(F_g), \text{ with }
  \tilde{\varphi}\in\Aut^+(\pi_g):
  [\tilde{\varphi}]=\varphi\in\Out^+(\pi_g)\},\] 
where we identify $\Mod_g \cong \Out^+(\pi_g)$ according to the theorem
of Dehn-Nielsen. If the
intersection of $\Aff^+(X,q)$ and $\Mod_g(\alpha)$ is trivial, then the
corresponding Teichmüller disc $\Delta_q$ is isomorphic to its image
under $s_\alpha$ in $\mathcal S_g$, hence an analytic disc.
But if for a flat surface $X$ the group $\overline{\Gamma}(X,q)$ is a
lattice (and, therefore, the corresponding Teichmüller disc leads to a
Teichmüller curve), then there exists an $\alpha$ such that
$\Aff^+(X,q)\cap\Mod_g(\alpha)\neq\{\id\}$, according to a result of
Herrlich and Schmithüsen\,\cite[5.21]{frankgabi}. This result is proven
in \cite{frankgabi} using compactifications of $\mathcal M_g$,
$\mathcal T_g$ and $\mathcal S_g$. Here we generalize this result to the
case where $\overline{\Gamma}(X,q)$ contains a parabolic element and we
give a constructive proof of it. Moreover, we show that there exists for
each parabolic element $\tau\in\Aff^+(X,q)$ an $\alpha$ such that
$\langle\tau^n\rangle\leq\Aff^+(X,q)\cap\Mod_g(\alpha)$ for some $n\in\N$.
Provided that $(X,q)$ is a translation surface (that means, the
transition maps are locally of the form $z\mapsto z+c$) we'll also show
that \hbox{$\Aff^+(X,q)\cap\Mod_g(\alpha)$} is either trivial or cylic
parabolic (Cor.\ \ref{isomzuzallg}). For $\Aff^+(X,q)\cap\Mod_g(\alpha)\neq\{\id\}$, this
leads to the fact that the image of $\Delta_q$ under $s_\alpha$ is
birational equivalent to $\HH/\langle\tau\rangle$ for a parabolic
element $\tau$, i.e.\ a punctured disc.\\ These statements will be proven
first for origamis (chapter \ref{kapori}) since this case can be handled
quite concretely and since 
at least the statement that
$\Aff^+(X,q)\cap\Mod_g(\alpha)$ is cyclic parabolic or trivial can be
shown by elementary methods.
A main ingredient to find an $\alpha$ with
$\Aff^+(X,q)\cap\Mod_g(\alpha)\neq\{\id\}$ is Thm.~\ref{hsss}, which
states that in every origami $O$ of genus $g$ there is a non-separating system
of $g$ simply closed curves that are all horizontal, which means that
for any $\varepsilon>0$ each of them is freely homotopic to a curve that lies
in an $\varepsilon$-neighbourhood of the horizontal sides of the squares
of $O$. It then follows that there is a symplectic set of generators
$(a_1,b_1,\dots,a_g,b_g)$ of the fundamental group of the origami such that
$a_1,\dots,a_g$ are freely homotopic to the horizontal curves. If we
choose $\alpha$ such that $\Kern(\alpha)$ is the normal subgroup
generated by $a_1,\dots,a_g$, then the element in $\Aff^+(O)$ that is given by a
multiple Dehn twist on all horizontal cylinders also lies in
$\Mod_g(\alpha)$. This result can be extended to parabolic elements with
any eigendirection and any flat surface. Since
the proof of the existence of a system of horizontal curves as described
above is constructive, it can (partially) be used to give an algorithm
that yields such a system of curves (chapter~\ref{kapalgorithmus}).\\
The two crucial ingredients to our proof that $\Aff^+(X,q)\cap\Mod_g(\alpha)$
is cyclic parabolic in the case that~$(X,q)$ is a translation surface is the fact that the
matrix $M_\varphi$ that describes the action of an element
$\varphi\in\Mod_g(\alpha)$ on the homology $H_1(X,\Z)$ is of the form
$M_\varphi=\begin{pmatrix}
I_g&A\\
0_g&I_g
\end{pmatrix}$
(Lemma~\ref{EW1}) and that any eigenvalue of the matrix of an element
$f\in\Aff^+(X,q)\hookrightarrow\Mod_g$ is also an eigenvalue of $M_f$
(Thm.~\ref{nichthypallgemein}). As already mentioned above, we give an
elementary proof of this latter fact for origamis, but rely on some more
sophisticated techniques in the case of general translation
surfaces~(chapter~\ref{kaphalbtrans}). Moreover, the very last example
in chapter~\ref{kapbsp} demonstrates that the above relation between the
respective eigenvalues does not hold for arbitrary translation surfaces.
So for flat surfaces, it is still an open question if
$\Aff^+(X,q)\cap\Mod_g(\alpha)$ can also be a non-cyclic group. A
promising strategy to look for such groups is probably to find (at
least) two parabolic elements with different eigendirections such that
for an $\alpha$ they both lie in $\Aff^+(X,q)\cap\Mod_g(\alpha)$.
Parabolic elements $f$ that act trivially on the homology (such elements
don't exist for translation surfaces) may be good candidates because in
this case the matrix $M_f$ is always of the form given above -- a
necessary condition for $f$ to lie in $\Mod_g(\alpha)$. As in the last
example of chapter \ref{kapbsp} it seems to be more probable that this
element is in $\Mod_g(\alpha)$ for several $\alpha$, and it should thus
be easier to find another element that does what it should.\\
I would like to thank Tobias Columbus for proofreading this english version of my PHD.

\chapter{Teichmüller discs}
\label{grundlagen}
In this chapter we recall some concepts that are important for the
subject of this thesis. More specifically, we introduce Teichmüller discs, Teichmüller
curves and projective Veech groups along with the necessary tools such as
holomorphic quadratic differentials and quasiconformal maps. Our
presentation largely follows \hbox{\cite[§2]{frankgabi}.}
\section{Riemann surfaces, moduli space and Teichmüller space}
\begin{definition}
\label{defrf}
\begin{enumerate}
\item 
  Let $X$ be a surface. A covering $\textfrak{U}$ by charts $U_i$ with
  transition maps $z_i\colon U_i\ra\C$ is called \begriff{complex atlas} if
  for all $(U_i,z_i), (U_j,z_j)$ the map
  \[ z_j\circ z_i^{-1}\colon z_i(U_i\cap U_j)\rightarrow z_j(U_i\cap U_j) \]
  is holomorphic.
\item 
  We call two atlases biholomorphic equivalent if their union is also a
  complex atlas. An equivalence class of biholomorphic equivalent
  atlases is called \begriff{complex structure}. In a complex structure
  there is a unique maximal atlas.
\item 
  A surface $X$ together with a complex structure $\Sigma$ is called
  \begriff{Riemann surface}. Most of the time, we'll omit the complex
  structure in the notation. The \begriff{genus} of a Riemann surface
  $(X,\Sigma)$ is the genus of $X$ known from the classification of
  surfaces (see \cite[§ 4.5]{kinsey}).
\end{enumerate}
\end{definition}

\begin{definition}
\label{flach}
Definition \ref{defrf} can be strengthened by further requiring that the
transition maps are of one of the following forms for every connected
component of an intersection $U_i\cap U_j$ (for $c\in\C$):\\\\
\begin{tabular}[c]{rcrlcrcrl}
(a)&$z_j\circ z_i^{-1}:$&$z_i(U_i\cap U_j)$&$\ra z_j(U_i\cap U_j)$&& (b)& $z_j\circ z_i^{-1}:$&$z_i(U_i\cap U_j)$&$\ra z_j(U_i\cap U_j)$\\
&&$z$&$\mapsto \pm z+c$&&&&$z$&$\mapsto z+c$
\end{tabular}\\\\
A complex structure satisfying condition (a) is called a \begriff{flat
  structure}, and a surface $X$ endowed with such a flat structure is
called a \begriff{flat surface}. Since the transition maps of a flat
structure are biholomorphic, a flat structure always induces a complex
structure. In case (b) we refer to the respective structure and surface
as \begriff{translation structure} and \begriff{translation surface}.
Again, we have that a translation
structure always induces a flat structure.
\end{definition}
From now on we'll focus on compact Riemann surfaces and compact
Riemann surfaces with finitely many marked points.
\begin{definition}
\begin{enumerate}
\item 
  A map $f\colon X\rightarrow Y$ between Riemann surfaces is called
  \begriff{holomorphic} if for all $p\in X$ the map $f$ is holomorphic on
  the charts, that is,  $w \circ f \circ z^{-1}$ is holomorphic for any charts $(U,z)$ of $p$ and $(V,w)$
  of $f(p)$ of the maximal atlas of $X$ respectively $Y$ with $f(U)\subseteq
  V$.\\
  We define \begriff{biholomorphic} maps in a completely analogous
  manner.
\item 
  Two Riemann surfaces $X$ and $Y$ are called \begriff{isomorphic} if
  there is a biholomorphic map $h\colon X\rightarrow Y$.
\item 
  The isomorphism classes of compact Riemann surfaces of genus $g$ form
  the so called \begriff{moduli space} of Riemann surfaces of genus $g$
  that we will denote by $\mathcal{M}_g$. We write $[X]$ for the isomorphism class of
  $X$.
\end{enumerate}
\end{definition}
\begin{definition}
\label{teich}
Take a compact Riemann surface $X_{\re}$ of genus $g$.
\begin{enumerate}
\item 
  A pair $(X,f)$ consisting of a Riemann surface $X$ of genus $g$ and an
  orientation preserving homeomorphism $f\colon X_{\re}\rightarrow X$
  that is not smooth in at most finitely many points is called
  \begriff{marked Riemann surface}.
\item 
  We define the following equivalence relation on the set of marked
  surfaces: $$(X,f) \sim (Y,h) :\Leftrightarrow h\circ f^{-1}\text{ is
    homotopic to a biholomorphic map }k:X\rightarrow Y.$$
\[
\begin{xy}
\xymatrix{
& X_{\re} \ar[ld]_f \ar[rd]^h & \\
X \ar[rr]^{k\sim h\circ f^{-1}} & & Y
}
\end{xy}
\]
\item We define \begriff{Teichmüller space} as the set of equivalence
  classes $[X,f]:= (X,f)/\sim$ of marked Riemann surfaces.
\end{enumerate}
\end{definition}
Teichmüller deformations are another possibility to define Teichmüller
space. For this purpose, we first introduce holomorphic quadratic
differentials (see also \cite[§ 4]{strebel}) and quasiconformal maps.

\section{Holomorphic quadratic differentials}
\begin{definition}
Take a Riemann surface $X$ with complex structure $\Sigma$, and let
$\textfrak A$ be the maximal atlas in $\Sigma$. A \begriff{holomorphic quadratic differential} $q$ on $X$ is given by a family
\begin{align*}
\{(\varphi_i(z_i) dz^2_i)_{i\in I}\mid&(U_i, h_i)\in\textfrak A, \varphi_i\colon h_i(U_i)\ra\C\text{ holomorphic},\\
&\forall i,j\in
I:\varphi_j(z_j)=\varphi_i(z_i)\left(\frac{dz_i}{dz_j}\right)^2 \text{
  for }h_i^{-1}(z_i)= h_j^{-1}(z_j)\}.
\end{align*}
We denote the $\C$-vector space of holomorphic quadratic differentials
on $X$ by $Q(X)$.
\end{definition}

\begin{bem}
The zeros of a quadratic holomorphic differential $q$, i.e.\ the points
$P\in X$ with $\varphi_i(h_i(P))=0$, are well-defined because they do
not depend on the chart. This follows from the fact that the transition
maps are biholomorphic since $\frac{dz_i}{dz_j}\neq0$.\\ But in in
general we have $\varphi_i(h_i(P))\neq\varphi_j(h_j(P))$ for $P\in
U_i\cap U_j\subseteq X$.
\end{bem}
Now let $X$ be a Riemann surface of genus $g\geq2$, $q\neq0$ a holomorphic quadratic differential on $X$, and let
\[X^*:=X-\{P\in X \vert P \text{ is zero of } q\}.\]
\begin{bem}
\label{hqm}
A flat structure $\mu$ on $X^*$ is given as follows: Let $U\subseteq
X^*$ be a simply connected chart of the maximal atlas of $X$ with
transition map $h$, and let $z_0\in h(U)$. Then we can define a new
chart contained in $U$ that lies in $\textfrak A$ by: 
$$P\mapsto\int_{z_0}^z\sqrt{\varphi(\xi)}d\xi;\quad z:= h(P).$$
Since the neighbourhood $U$ is simply connected, the integral does not
depend on the path between $z_0$ and $z$ in $U$. The charts that we get
in this way form a flat structure on $X^*$. The function
$\sqrt{\varphi(z)}$ exists in the whole chart $U$, but is only
well-defined up to multiplication by $(-1)$. For a chart we, therefore, have to choose a branch of it.
\end{bem}
\begin{Bew}
Let $(U_{\nu}, h_{\nu})$ and $(U_{\eta},h_{\eta})$ be two simply
connected charts of $P$ with the corresponding transition maps. Define
$z_{\nu}:= h_{\nu}(P)$ and $z_{\eta}:= h_{\eta}(P)$, fix a point $P_0\in
U_\nu\cap U_\eta$, and let $z_{\nu_0}:=h_\nu(P_0)$,
$z_{\eta_0}=h_\eta(P_0)$ and $h_{\nu\eta}:= h_{\eta}\circ h_{\nu}^{-1}$.
Then we have:
$$\int_{z_{\nu_0}}^{z_{\nu}}\sqrt{\varphi_{\nu}(\xi)}d\xi = \int_{z_{\nu_0}}^{z_{\nu}}\sqrt{\varphi_{\eta}(h_{\nu\eta}(\xi))}\frac{dh_{\nu\eta}(\xi)}{d\xi}d\xi =
\int_{z_{\eta_0}}^{z_{\eta}}\sqrt{\varphi_{\eta}(\xi)}d\xi$$
As both $P_0$ and the sign of $\sqrt{\varphi(z)}$ may be chosen
arbitrarily, it follows that these charts give rise to a flat structure.
\end{Bew}
The flat structure induced by $q$ induces an euclidean metric on $X^*$.
We define $\|q\|$ to be the area of $X$ with respect to this metric. It
follows from the proof of the above remark that $\|r\cdot
q\|=|r|\cdot\|q\|$ for $r\in\R$ and $q\in Q(X)$.

\section{Quasiconformal maps and the theorems of Teichmüller}
This chapter follows mainly \cite{farbmarg}. Further references that
also give a good introduction to the subject are \cite{ahlqc} and \cite{imtan}.
\begin{definition}
Take open sets $U,V\subseteq\C$ and an orientation preserving
homeomorphism $f:U\ra V$ that is real differentiable except at finitely many points.
For a point $p\in U$, where $f$ is real differentiable, define
\[K_f(p):=\frac{|f_z(p)|+|f_{\z}(p)|}{|f_z(p)|-|f_{\z}(p)|}=\frac{1+\frac{|f_{\z}(p)|}{|f_z(p)|}}{1-\frac{|f_{\z}(p)|}{|f_z(p)|}}\]
as the \begriff{complex dilatation} of $f$ at the point $p$.
\end{definition}
Roughly speaking, the complex dilatation tells how far $f$ is from being
conformal at the point $p$. For $K_f(p)=1$, $f$ is conformal at $p$. Since $f$ is orientation preserving, we have
\[\jac(f)_p=|f_z(p)|^2-|f_{\z}(p)|^2>0,\]
and, therefore, $\frac{|f_{\z}(p)|}{|f_z(p)|}<1$. It follows that $K_f(p)\geq1$.\\
The complex dilatation $K_f(p)$ can also be interpreted as follows: The map between the tangent spaces
\[df_p:T_pU\ra T_pV\]
maps the unit circle to an ellipse. $K_f(p)$ is then the ratio between the major axis to the minor axis of this ellipse.
\begin{definition}
For $U,V\subseteq\C$, let $f:U\ra V$ be as above. The \begriff{complex dilatation of $f$} is defined as
\[K(f):=\sup_{p\in U}K_f(p).\]
If it is finite, then $f$ is called \begriff{quasiconformal}.
\end{definition}
We can define the complex dilatation of a map between Riemann surfaces
analogous to the case of $\C$, simply as complex dilatation of the
corresponding map on the charts. This definition makes sense since the
transition maps on Riemann surfaces are biholomorphic, hence, the
complex dilatation of a map between Riemann surfaces does not depend on
the charts (see \cite[Prop. 11.3]{farbmarg}).\\ Teichmüller showed that
in every equivalence class of markings on a Riemann surfaces there
exists a unique map that is differentiable except at finitely many
points and has minimal complex dilatation.
\begin{definition}
Let $X$ and $Y$ be two compact Riemann surfaces of genus $g$. A \begriff{Teichmüller map} $f:X\ra Y$ is a homeomorphism such that there are holomorphic quadratic differentials $q_X$ on $X$ and $q_Y$ on $Y$ and a $K\geq1$ with the following properties:
\begin{itemize}
\item $f$ maps the zeros of $q_X$ to the zeros of $q_Y$,
\item in points that are not zeros of $q_X$ the map $f$ has locally the form
\[f(z)=f(x+iy)=\sqrt{K}\cdot x+\frac{1}{\sqrt{K}}\cdot iy\]
with respect to suitable charts in the flat structures that come from $q_X$ and $q_Y$.
\end{itemize}
Note that the constant $K$ does not depend on the charts.
\end{definition}
\begin{bem}
For a Teichmüller map $\;f\colon x+iy\mapsto \sqrt{K}\cdot x+\frac{1}{\sqrt{K}}\cdot iy\;$, we have $K_f=K$.
\end{bem}
\begin{Bew}
We have
\[f(z)=\frac{1}{2}\sqrt{K}\cdot (z+\z)+\frac{1}{2\sqrt{K}}\cdot (z-\z)\]
and, therefore, \[f_z(z)=\frac{1}{2}\left(\sqrt{K}+\frac{1}{\sqrt{K}}\right)\quad\text{ and }\quad f_{\z}(z)=\frac{1}{2}\left(\sqrt{K}-\frac{1}{\sqrt{K}}\right).\]
The complex dilatation is independent of $z$, and we have
\[K_f=\frac{\left(\sqrt{K}+\frac{1}{\sqrt{K}}\right)+\left(\sqrt{K}-\frac{1}{\sqrt{K}}\right)}{\left(\sqrt{K}+\frac{1}{\sqrt{K}}\right)-\left(\sqrt{K}-\frac{1}{\sqrt{K}}\right)}=K.\]
\end{Bew}

\begin{satz}[Teichmüller's Existence Theorem]
\label{teichex}
Let $X$ and $Y$ be two compact Riemann surfaces of genus $g\geq1$ and $h:X\ra Y$ a homeomorphism. Then there is a Teichmüller map $f\colon X\ra Y$ in the homotopy class of $h$ .
\end{satz}

\begin{Bew}
See \cite[11.8]{farbmarg}.
\end{Bew}

\begin{satz}[Teichmüller's Uniqueness Theorem]
\label{teicheind}
Let $X$ and $Y$ be two compact Riemann surfaces of genus $g\geq1$, $f:X\ra Y$ a Teichmüller map and $h$ a quasiconformal map homotopic to $f$.
Then we have
\[K_h\geq K_f.\]
Equality holds if and only if $f\circ h^{-1}$ is conformal. For $g\geq2$, this second assertion means $f=h$.
\end{satz}

\begin{Bew}
See \cite[11.9]{farbmarg}.
\end{Bew}
It follows from these two theorems of Teichmüller that every point $[X,h]\in\mathcal T_g$ corresponds to a unique marked Riemann surface $(X,f)$ with $(X,f)\sim(X,h)$, where $f:X_{\re}\ra X$ is a Teichmüller map. With this we can define the following metric on $\mathcal T_g$:
\begin{definition}
Let $x,y\in\mathcal T_g$, represented by $(X,f_1)$ and $(X,f_2)$, where
$f_1$ and $f_2$ are the Teichmüller maps. Moreover, let $f$ be the
Teichmüller map homotopic to $f_2\circ f_1^{-1}$, which exists by
Thm.~\ref{teichex} and is unique by Thm.~\ref{teicheind}. We define the
\begriff{Teichmüller metric} as
\[d_T(x,y):=\log K_f.\]
\end{definition}
This is actually a metric because for quasiconformal maps $f:X\ra Y$ and $g\colon Y\ra Z$, we have $K_{g\circ f}\leq K_g\cdot K_f$ (therefore, the triangle inequality holds), $K_f=K_{f^{-1}}$ (therefore, $d_T$ is symmetric) and  $K_{\id}=1$, see \cite[11.3]{farbmarg}. 

\section{Teichmüller deformations}
The flat structure that we constructed on a Riemann surface $X$ with quadratic holomorphic differential $q$ can be deformed by composing the transition maps with the map
\begin{equation}
\label{teichdef}
\vartheta_K\colon x+iy \mapsto \sqrt{K}\cdot x+\frac{1}{\sqrt{K}}\cdot iy,\quad K>1.
\end{equation}
By doing this, we get on the (topological) surface $X$ a new flat structure, that induces in general also a new complex structure. 
We denote this new Riemann surface, that is, the (topological) surface
$X$ with this new complex structure, by $X_K$ or $X_{K,q}$ if the quadratic holomorphic differential $q$ shall be emphasized.
\begin{definition}
We get (with $X_{\re}:=X$) a marked Riemann surface $(X_K,f_K)$, where
\[f_K\colon X\overset{\id}{\rightarrow} X_{K,q}\]
is the map that is topologically the identity (and, therefore,
corresponds to the map $\vartheta_K$ on the charts). We call $(X_K,f_K)$
a \begriff{Teichmüller deformation} of $X$ with respect to $q$ with
dilatation $K$.
\end{definition}
Note that $f_K$ is also the Teichmüller map, and we have
\[d_T([X,\id],[X_K,f_K])=\log K.\]

\begin{bem}
Take some holomorphic quadratic differential $q$ and a real multiple
$rq$ for some $r>0$, and let $K>1$. Then $[X_{K,q},f_K]=[X_{K,rq},f_K]$
in $\mathcal{T}_g$ and, therefore, it suffices to consider the set
\[\Sigma_X:=\{q\in Q(X)\mid||q||=1\}.\]
\end{bem}
\begin{satz}
Let $X$ be a Riemann surface of genus $g$. Every point in $\mathcal T_g$
can be obtained uniquely as Teichmüller deformation of $X$. Therefore, we have the following bijection:
\[\Sigma_X\times(1,\infty)\cup\lbrace 0 \rbrace\stackrel{1:1}{\longleftrightarrow}\mathcal{T}_g.\]
\end{satz}
\begin{Bew}
Let $[Y,f]\in\mathcal T_g$. By Thm. \ref{teichex}, we can assume
w.l.o.g. that $f\colon X\ra Y$ is the Teichmüller map. Let $q_X\in\Sigma_X$
and $q_Y\in\Sigma_Y$ be the corresponding (normed) holomorphic quadratic
differentials on $X$ and $Y$, and let $(X_K,f_K)$ be the Teichmüller
deformation of $X$ with respect to $q_X$ with dilatation $K:=K_f$. Then
$f_K$ and $f$ correspond on the charts to the map $\vartheta_K$, and so
$f_K\circ f^{-1}$ is the identity, hence biholomorphic. Therefore, 
$(Y,f)\sim(X_K,f_K)$, which implies $[Y,f]=[X_K,f_K]$.\\
Now, let $[X_{K,q},f_K]=[X_{K,\tilde{q}},t_K]$ be the two
representations. Then there is a map \hbox{$h\colon X_{K,\tilde{q}}\ra
  X_{K,\tilde{q}}$} homotopic to the identity such that $h\circ t_K\circ
f^{-1}$ is a biholomorphic map.
\[
\begin{xy}
\xymatrix{
& X \ar[ld]_{f_K} \ar[rd]^{t_K} & \\
X_{K,q} \ar[rr]^{t_K\circ f_K^{-1}} & & X_{K,\tilde{q}}
}
\end{xy}
\]
Since $X_{K,q}$ and $X_{K,\tilde{q}}$ are the same topologic surface,
$t_K$ is also a quasiconformal map from $X$ to $X_{K,q}$ and $h\in
Q(X_{K,q})$. Since $t_K$ and $f_K$ are topologically the identity, we
have by Thm. \ref{teicheind} that $h\circ t_K= f_K$, and, therefore, we
finally have $h=\id$ and $f_K=t_K$.
\end{Bew}

\section{The complex structure of Teichmüller space}
On the Teichmüller space we can define a complex structure. One
possibility to do this, is to lift the complex structure of the Schottky
space, introduced in chapter~\ref{schottky}, to the Teichmüller space
via the map $s_\alpha$ (see Thm.~\ref{teichquotabb}) because this map
is, by Thm.~\ref{teichquotabb}b, the universal covering of the Schottky
space. Since the forgetful map from the Schottky space to the moduli
space is analytic (Thm.~\ref{vergissanalytisch}), composition of it with
$s_\alpha$ (altogether the quotient map after the action of the mapping class group
$\Mod_g$ defined later) is also analytic. Another -- more classic -- possibility
to construct the complex structure on $\mathcal T_g$ is described
in detail in \cite{nag}. We give here a short sketch of it.\\ We can
associate to a marked Riemann surface $(X,f)$ a so called
\begriff{Beltrami differential} by $\mu_f(z):=\frac{f_{\z}}{f_z}$ (see
\cite[1.3.1]{nag}). The Beltrami differential $\mu_f$ is an element in
the Banach space $L^{\infty}_{(-1,1)}$ of the $(-1,1)$-forms on $X$ with
norm $\|\mu_f\|_\infty<1$. Now, let $B(X)$ be the open unit ball in
$L^{\infty}_{(-1,1)}$ and \[M(X):=\{(X,f)\mid X\text{ Riemann surface, }
  f:X_{\re}\ra X\text{ quasiconformal}\}/\sim\] with
$((X_1,f_1)\sim(X_2,f_2):\Leftrightarrow f_1\circ f_2^{-1}$ is
conformal$)$. Then the map \[M(X)\ra B(X),\quad(X,f)\mapsto \mu_f\] is
bijective \cite[2.1.4]{nag}. Since there is a natural projection
$M(X)\ra\mathcal T_g$, we have a projection $\Phi\colon B(X)\ra\mathcal
T_g$. It can be shown that there is a structure of a complex manifold on
$\mathcal T_g$ such that $\Phi$ becomes a holomorphic map
\cite[3.1]{nag}.

\section{Teichmüller discs}
\begin{definition}
Let $q$ be a holomorphic quadratic differential on a compact Riemann surface $X$ of genus $g$. Then we can define the following map:
\begin{align*}
\gamma \colon [0,\infty)&\rightarrow \mathcal{T}_g\\
t&\mapsto [X_K,f_K]\text{ with }K=e^t
\end{align*}
The image of $\gamma$ in $\mathcal{T}_g$ is called \begriff{geodesic ray with respect to q}.
\end{definition}
The map $\gamma$ is an isometric embedding with respect to the Teichmüller metric in $\mathcal T_g$.
\begin{definition}
Let $g\geq2$. A \begriff{Teichmüller disc}~$\Delta_{\iota}$ is the image of a holomorphic isometric embedding
\[\iota \colon \mathds{D} \hookrightarrow \mathcal{T}_g,\]
where $\mathds{D}:=\{z\in\C:\vert z\vert<1\}$ endowed with the Poincar\'{e} metric.
\end{definition}
There are several possibilities to describe such an embedding. We expose two of them:
\paragraph{First possibility: As set of geodesic rays}
\begin{align*}
\iota_1 :\mathds{D} &\rightarrow \mathcal{T}_g\\
z&\mapsto (X_{K,e^{-i\varphi}\cdot q},f_K)\text{ with }z=\frac{K-1}{K+1}\cdot e^{i\varphi}, K\in[1,\infty), \varphi\in[0,2\pi)
\end{align*}
For every $\varphi$, we get a geodesic ray.
\paragraph{Second possibility: By affine deformations}~\\
Here we compose the transition maps of the flat structure $\mu$ with an affine map $B$ (therefore, we identify $\mathds{C}$ with $\mathds{R}^2$) and get a new flat structure $\mu_B$. In other words: $\SL_2(\mathds{R})$ acts on the flat structures of $X^*$, more precisely: $B=\begin{pmatrix}a&b\\c&d\end{pmatrix}\in\SL_2(\mathds{R})$ acts on $\R^2$ via
$$\begin{pmatrix}
x\\
y\\
\end{pmatrix}
\mapsto
\begin{pmatrix}
a & b\\
c& d\\
\end{pmatrix}
\cdot
\begin{pmatrix}
x\\
y\\
\end{pmatrix}.$$
In $\C$ this means
$$(x+iy)\mapsto (ax+by)+i(cx+dy).$$
We then write $B\circ (X,\mu):= (X,\mu_B)$. Let $X_{\re}$ be the Riemann surface with the complex structure induced by the flat structure $\mu$, then
\[(X,\mu)\stackrel{\id}{\rightarrow}(X,\mu_B)\]
defines a marked Riemann surface, which yields the following point in Teichmüller space:
\[P_B:= [(X,\mu_B),\id].\] 
In the next theorem we summarize some properties of the action of $\SL_2(\mathds R)$ on the set of flat structures:

\begin{satz}
\label{sloperation}
\begin{enumerate}
\item $B_1\circ (B_2\circ (X,\mu))= (X,\mu_{B_1B_2})$ for $B_1,B_2\in\SL_2(\mathds R)$.
\item $B\in \SL_2(\mathds{R})$ leaves the point in the Teichmüller space fixed if and only if $B\in \SO_2(\mathds{R})$.
\item For $B=
\begin{pmatrix}
\sqrt{K} & 0\\
0&\frac{1}{\sqrt{K}}\\
\end{pmatrix}$, we get the Teichmüller deformation $\vartheta_K$ defined in (\ref{teichdef}).
\item The following map is injective:
\[\tilde \iota_2\colon \SO_2(\mathds{R})\backslash\SL_2(\mathds{R}) \rightarrow \mathcal{T}_g,\quad\SO_2(\mathds{R})\cdot B \mapsto P_B\]
\item The map
\[\hat\iota_2\colon \SO_2(\mathds{R})\backslash\SL_2(\mathds{R})\rightarrow\mathds{H},\quad\SO_2(\mathds{R})\cdot A \mapsto -\overline{A^{-1}(i)}\]
is bijective, where we mean by $A(z)$ the action via Möbius transformations.
\item The map $\hat\iota_2^{-1}$ inverse to the above map
  $\hat{\iota}_2$ is induced by
$$\mathds{H} \rightarrow \SL_2(\mathds{R}),~~~~t \mapsto \frac{1}{\sqrt{\operatorname{Im}(t)}}
\begin{pmatrix}
1&\operatorname{Re}(t)\\
0& \operatorname{Im}(t)
\end{pmatrix}.$$
\end{enumerate}
\end{satz}
\begin{Bew}
\begin{enumerate}
\item and (c) are obvious.
\item Let $z=x+iy\in\mathds C$ and $B\in\SO_2(\mathds R)$. Then $B$ is of the form $B=
\begin{pmatrix}
a & b\\
-b & a\\
\end{pmatrix}$ with $a^2+b^2=1$ and $a,b\in\mathds R$.\\
Then we have $Bz=(ax+by)+i(-bx+ay)=(a-ib)z$. Therefore, $B$ acts
holomorphic on $\mathds C$, and it follows that $P_B=P_I$.\\
Now, let $B=
\begin{pmatrix}
a & b\\
c & d\\
\end{pmatrix}\in\SL_2(\mathds R)$, and assume $P_B=P_I$:\\
In this case, the map
$$x+iy\mapsto(ax+by)+i(cx+dy)=a\frac{z+\overline{z}}{2}-ib\frac{z-\overline{z}}{2}+ic\frac{z+\overline{z}}{2}+d\frac{z-\overline{z}}{2}$$
is holomorphic. Therefore, its partial derivative with respect to $\overline{z}$ is the zero map, and so
$$\frac{\overline{z}}{2}(a-d)+i\frac{\overline{z}}{2}(b+c)=0.$$
It follows that $a=d$ and $b=-c$, and we thus have $B\in\SO_2(\mathds R)$.
\item[(d)] Let $A,B\in \SL_2(\mathds R)$ with $P_A=P_B$. It follows from a) and b) that $BA^{-1}\in\SO_2(\mathds R)$, and, therefore, $B\in\SO_2(\mathds R)\cdot A$.
\item[(e)] Let $A=\begin{pmatrix}a & b\\c & d\\\end{pmatrix}\in\SL_2(\mathds R)$. Then we have 
$$A\in\Stab(i)\Leftrightarrow\frac{ai+b}{ci+d}=i\Leftrightarrow(a=d\land b=-c)\Leftrightarrow A\in\SO_2(\mathds R).$$
$\SO_2(\mathds R)=\Stab(i)$ implies that the map is well-defined and injective.\\
Since $\SL_2(\mathds R)$ acts transitively on $\mathds H$, the map
$A\mapsto A(i)$ is surjective, too.\\
The inversion of $A$, the complex conjugation and the multiplication by
$-1$ will become important later on. Observe that they preserve the bijectivity of the map.
\item[(f)] Let $A=\begin{pmatrix}
a & b\\
c & d\\
\end{pmatrix}\in\SL_2(\mathds R)$. We have $-\overline{A^{-1}(i)}=-\overline{\frac{di-b}{-ci+a}}=\frac{1}{a^2+c^2}(ab+cd+i):=t$ and 
\[\frac{1}{\sqrt{\Imag(t)}}\begin{pmatrix}
1 & \Real(t)\\
0 & \Imag(t)\\
\end{pmatrix}=\sqrt{a^2+c^2}\begin{pmatrix}
1 & \frac{ab+cd}{a^2+c^2}\\
0 & \frac{1}{a^2+c^2}\\
\end{pmatrix}:=\tilde A\in\SL_2(\mathds R).\]
Therefore, we have
\[A\tilde A^{-1}=\sqrt{a^2+c^2}\begin{pmatrix}
a & b\\
c & d\\
\end{pmatrix}\cdot\begin{pmatrix}
\frac{1}{a^2+c^2} & -\frac{ab+cd}{a^2+c^2}\\
0 & 1\\
\end{pmatrix}=\sqrt{a^2+c^2}\begin{pmatrix}
\frac{a}{a^2+c^2} & -a\frac{ab+cd}{a^2+c^2}+b\\
\frac{c}{a^2+c^2} & -c\frac{ab+cd}{a^2+c^2}+d\\
\end{pmatrix}\]
\[=\sqrt{a^2+c^2}\begin{pmatrix}
\frac{a}{a^2+c^2} & \frac{-acd+bc^2}{a^2+c^2}\\
\frac{c}{a^2+c^2} & \frac{-cab+da^2}{a^2+c^2}\\
\end{pmatrix}=\frac{1}{\sqrt{a^2+c^2}}\begin{pmatrix}
a & -c\\
c & a\\
\end{pmatrix}\in\SO_2(\mathds R).\]
\end{enumerate}
\end{Bew}
Now, we define the map $\iota_2:=\tilde\iota_2\circ\hat\iota_2^{-1}$, in other words:
\begin{align*}
\iota_2 :\; \mathds{H} &\rightarrow \mathcal{T}_g\\
t&\mapsto P_{A_t},
\end{align*}
where we choose $A_t$ such that $-\overline{A_t^{-1}(i)}=t$.
\begin{satz}
$\iota_1$ and $\iota_2$ are Teichmüller embeddings and define the same Teichmüller disc:
\[\Delta_q := \Delta_{\iota_1}= \iota_1 (\mathds{D})=\Delta_{\iota_2}=\iota_2(\mathds{H})\]
\end{satz}
\begin{Bew}
We have $\iota_2 = \iota_1 \circ f$, where $f(z):=\frac{i-z}{i+z}$ maps $\mathds{H}$ to $\mathds{D}$, see \cite[Prop.~2.12]{frankgabi}.
By  \linebreak\cite[Prop.~2.11]{frankgabi}, $\iota_2$ (and, therefore, also
$\iota_1$) is an isometry (with respect to both the Teichmüller metric
and the hyperbolic metric), and, by \cite[Cor. 2.15]{frankgabi},
$\iota_1$ (and, therefore, also $\iota_2$) is holomorphic. The assertion
follows.
\end{Bew}

\section{Teichmüller curves and projective Veech groups}
\label{teichveech}
\begin{definition}
\label{Veechfl}
Let $X$ be a Riemann surface, $q$ a holomorphic quadratic differential
on $X$ and $\Delta_q$ the Teichmüller disc induced by $q$. We call $C_q$ a \begriff{Teichmüller curve} and $(X,q)$ a \begriff{Veech surface} if the image $C_q$ of $\Delta_q$ under the natural projection $\mathcal{T}_g\rightarrow \mathcal{M}_g,[X,f]\mapsto [X]$ is an algebraic curve.
\end{definition}
In order to say something about when $(X,q)$ leads to a Teichmüller curve, we first introduce the notion of a projective Veech group:\\
Take a pair $(X,q)$, where $X$ is a Riemann surface of genus $g$ and $q$ is a holomorphic quadratic differential on it, and let $\mu$ be the corresponding flat structure. Moreover, let $\Aff^+(X,\mu)$ be the group of orientation preserving maps that are affine with respect to the flat structure $\mu$, i.e.\ have the following form on the charts:
$$z\mapsto A\cdot z + t~~~~\textrm{for}~~A\in \SL_2(\mathds{R}),~~ t \in \mathds{C}$$
Since we have a flat structure, $A$ does not depend on the chart up to sign. So we get the following group homomorphism:
$$D\colon \Aff^+(X,\mu) \rightarrow \PSL_2(\mathds{R}), ~~~f\mapsto [A]$$

\begin{definition}
\label{projveech}
$\overline{\Gamma}(X,\mu):=D(\Aff^+(X,\mu))$ is called \begriff{projective Veech group} of $(X,\mu)$.
\end{definition}
If $\mu$ is a translation structure, we even have a group homomorphism:
\[\der\colon \Aff^+(X,\mu) \rightarrow \SL_2(\mathds{R}), ~~~f\mapsto A\]
Therefore, we define (for later use in chapter \ref{kapori}) also the Veech group.
\begin{definition}
\label{projveech-2}
Let $(X,\mu)$ be a translations surface. Then we call $\Gamma(X,\mu):=\der(\Aff^+(X,\mu))$ the \begriff{Veech group} of $(X,\mu)$.
\end{definition}
Now, we would like to analyze how the projective Veech group acts on a Teichmüller disc. Therefore, we first introduce the mapping class group.
\begin{definition}
Take a compact, oriented surface $S$ of genus $g$. Let $\Homo^+(S)$ be the group of orientation preserving homeomorphisms $S\ra S$, and let $\Homo^0(S)$ denote the subgroup of homeomorphisms that are isotopic to the identity. Then we define the \begriff{mapping class group}
\[\Mod_g:=\Mod(S):= \Homo^+(S)/\Homo^0(S).\]
\end{definition}
\begin{bem}
\begin{enumerate}
\item In the above definition of $\Mod_g$, we could also take the group of orientation preserving diffeomorphisms  $\Diffeo^+(S)$ and the subgroup of the diffeomorphisms $\Diffeo^0(S)$ isotopic to the identity instead of $\Homo^+(S)$ and $\Homo^0(S)$, see \cite[§2.1]{farbmarg}.
\item The mapping class group can also be defined for non-compact and/or bounded surfaces, see \cite[§2.1]{farbmarg}.
\item $\Aff^+(X,\mu)$ is a subgroup of $\Homo^+(X)$.
\end{enumerate}
\end{bem}
Let $X$ be a Riemann surface of genus $g$. The group $\Homo^+(X)$ acts on $\mathcal T_g$ (where we take $X$ as reference surface) via
\begin{align*}
\Homo^+(X)\times\mathcal T_g&\ra\mathcal T_g\\
(\varphi,[Y,f])&\mapsto [Y,f\circ\varphi^{-1}].
\end{align*}
We denote $[Y,f\circ\varphi^{-1}]$ by $\rho_\varphi([Y,f])$.\\
Therefore, also $\Aff^+(X,\mu)$ acts, as subgroup of $\Homo^+(X)$, on $\mathcal T_g$.\\
Since $\Homo^0(X)$ acts trivially on $\mathcal T_g$, the action of $\Homo^+(X)$ descends to an action of $\Mod_g$ on $\mathcal T_g$. As a subgroup of $\Homo^+(X)$, $\Aff^+(X,\mu)$, therefore, also acts on $\mathcal T_g$. The quotient map induced by this action is the natural projection $\mathcal{T}_g\rightarrow \mathcal{M}_g$.

\begin{satz}
\label{affoperiertauftg}
$\Aff^+(X,\mu)$ stabilizes $\Delta_q$ and acts on it via
\[\varphi \in \Aff^+(X,\mu), B\in \SL_2(\mathds{R}) \Longrightarrow \rho_{\varphi}(P_B)=P_{BA^{-1}},\]
where $A\in \SL_2(\mathds{R})$ is the preimage of $D(\varphi)=[A].$
\end{satz}
\begin{Bew}
\[
\begin{xy}
\xymatrix{(X,\mu )\ar[r]^{\varphi^{-1}}\ar[rrd]_{\id} & (X,\mu )\ar[r]^{\id} & (X,\mu_B)\\
& & (X,\mu_{BA^{-1}})\ar[u]}
\end{xy}
\]
The map $(X,\mu_{BA^{-1}})\rightarrow (X,\mu_B)$ in the diagram is an affine map as composition of affine maps. Moreover, we have: $D(BA^{-1}(BA^{-1})^{-1})=I$, and so the map is holomorphic, which means $[(X,\mu_{BA^{-1}}), \id] = [(X,\mu_B), \id\circ\varphi^{-1}]$. By definition, we have $\rho_{\varphi}([(X,\mu_B), \id])= [(X,\mu_B), \id\circ\varphi^{-1}]$, and so the assertion $\rho_{\varphi}([(X,\mu_B), \id])= [(X,\mu_{BA^{-1}}), \id]$ follows.
\end{Bew}

\begin{bem}
\label{affstab}
We even have $\Aff^+(X,\mu)\cong\Stab(\Delta_q):=\Stab_{\Mod_g}(\Delta_q)$.
\end{bem}
\begin{Bew}
See \cite[Thm 1, Lemma 5.2]{earlegard}.
\end{Bew}

\begin{satz}
We have
\[\overline{\Gamma}(X,\mu) \cong \FakRaum{\Stab(\Delta_q)}{\Stab_{\pw}(\Delta_q)},\]
where $\Stab_{\pw}(\Delta_q):=\{\varphi\in\Mod_g\colon \varphi\vert\Delta_q=\id_{\Delta_q}\}$.
\end{satz}

\begin{Bew}
Let \[\pi:\Stab(\Delta_q)\ra\FakRaum{\Stab(\Delta_q)}{\Stab_{\pw}(\Delta_q)}\] be the quotient map. By the proof of Thm. \ref{affoperiertauftg}, the action of $\varphi\in\Aff^+(X,\mu)$ on $\mathcal T_g$ only depends on $[A]$, and, therefore, $\pi\circ\rho$ factorizes over $\overline{\Gamma}(X,\mu)$. Let $\overline{\rho}$ the homomorphism induced by $\rho$ such that the following diagram commutes:
\[\begin{xy}
\xymatrix{\Aff^+(X,\mu)\ar[d]_D\ar[rr]^{\rho}&&\Stab(\Delta_q)\ar[d]^\pi\\
\overline{\Gamma}(X,\mu)\ar[rr]^{\overline{\rho}}&&\FakRaum{\Stab(\Delta_q)}{\Stab_{\pw}(\Delta_q)}&}
\end{xy}
\]
By remark \ref{affstab}, $\rho$ is an isomorphism. It remains to show that $\overline{\rho}$ is injective. For this purpose, let \hbox{$\varphi\in\Aff^+(X,\mu)$} with $\overline{\rho}(D(\varphi))=\id\vert_{\Delta_q}$ and let $A\in\SL_2(\R)$ be a preimage of $[A]:= D(\varphi)$. Then we have in particular $\overline{\rho}_\varphi(P_I)=P_I=P_{A^{-1}}$ and, by Thm. \ref{sloperation}b, 
\[A^{-1}= \begin{pmatrix}
a & b\\
-b & a
\end{pmatrix}\in\SO_2(\R),\quad\text{ for } a,b\in\R\text{ and }a^2+b^2=1.\]
Now, take  
$B= \begin{pmatrix}
\lambda & 0\\
0&\frac{1}{\lambda}
\end{pmatrix}$ with $\pm1\neq\lambda\in\R$. By assumption, we have $\overline{\rho}_\varphi(P_B)=P_B=P_{BA^{-1}}$. Therefore, the affine transformation $f$ that is induced by $B$, i.e. the transformation
\[f:z\mapsto \lambda\Real(z)+i\frac{1}{\lambda}\Imag(z)\]
is biholomorphic equivalent to the affine transformation $g$ induced by $BA^{-1}=\begin{pmatrix}
\lambda a & \lambda b\\
-\frac{b}{\lambda}&\frac{a}{\lambda}
\end{pmatrix},$ that is
\[g:z\mapsto \lambda a\Real(z)+\lambda b\Imag(z)+i\left(\frac{a}{\lambda}\Imag(z)-\frac{b}{\lambda}\Real(z)\right).\]
It is to show that $A=\pm I$.\\
We have
\begin{align*}
g\circ f^{-1}\colon z\mapsto\;&a\Real(z)+\lambda^2 b\Imag(z)+i\left(a\Imag(z)-\frac{b}{\lambda^2}\Real(z)\right)\\
&=\frac{1}{2}\left[a(z+\overline{z})-i\lambda^2 b(z-\overline{z})+a(z-\overline{z})-i\frac{b}{\lambda^2}(z+\overline{z})\right],
\end{align*}
and since $g\circ f^{-1}$ is holomorphic,
\[0=\left(g\circ f^{-1}\right)_{\overline{z}}=\frac{1}{2}ib\left(\lambda^2-\frac{1}{\lambda^2}\right).\]
Since $\lambda\neq\pm1$, it follows that $b=0$, and so we have $A=\pm I$.
\end{Bew}
We now compare the action of $\overline{\Gamma}(X,\mu)$ on $\Delta_q$ and the action of $\overline{\Gamma}(X,\mu)$ on $\mathds{H}$:
\begin{satz}
\label{rar}
Let $A\in\overline{\Gamma}(X,\mu) \leq \PSL_2(\mathds{R})$. Then the following diagram commutes:
\[
\begin{xy}
\xymatrix{\mathds{H}\ar[rr]^{t \mapsto -\overline{t}}\ar[d]^A && \mathds{H}\ar[rr]^{\iota}\ar[d]^{RAR^{-1}} && \Delta_q\ar[d]^{\rho_A}\\
\mathds{H}\ar[rr]^{t \mapsto -\overline{t}} && \mathds{H}\ar[rr]^{\iota} && \Delta_q}
\end{xy}
\]
where $R=R^{-1}=
\begin{pmatrix}
-1 & 0 \\
0 & 1
\end{pmatrix}$
(or $z\mapsto - \overline{z}$).
\end{satz}
\begin{Bew} Let $t\in\mathds{H}, B\in \SL_2(\mathds{R})$ with $-\overline{B^{-1}(i)}=-\overline{t}$. Then we have $\iota(-\overline{t})=P_B$.
\[
\begin{xy}
\xymatrix{
t\ar[r]^{t \mapsto -\overline{t}}\ar[d]^A & -\overline{t}\ar[r]^{\iota}\ar[d]^{RAR^{-1}} & P_B\ar[d]^{\rho_A}\\
A(t)\ar[r]^{t \mapsto -\overline{t}} & -\overline{A(t)}\ar@.[r] & P_{BA^{-1}}}
\end{xy}
\]
It is to show that the punctured arrow in the diagram corresponds to the map $\iota$. This is equivalent to $-\overline{(BA^{-1})^{-1}(i)}=-\overline{A(t)}$. However, this raedily follows from the computation
\[-\overline{(BA^{-1})^{-1}(i)}= -\overline{AB^{-1}(i)}= -A(\overline{B^{-1}(i)})=-A(\overline{t})=-\overline{A(t)}.\]
\end{Bew}
Let us now reconsider the question about when we get an algebraic curve as image of a Teichmüller disc $\Delta_q$ under the projection $\pr\colon \mathcal{T}_g \rightarrow \mathcal{M}_g$.
\begin{satz}
\label{algkurveechgitter}
$\pr(\Delta_q)$ is an algebraic curve if and only if the corresponding projective Veech group is a lattice in $\PSL_2(\mathds{R})$ (that is a Fuchsian group with finite covolume).
\end{satz}
\begin{Bew} See also \cite[Cor 3.3]{mcmullen}.\\
By Thm. \ref{rar}, the map $\pr\circ\iota \colon \mathds{H}\rightarrow \pr(\Delta_q) \subseteq \mathcal{M}_g$ factorizes through
$$\overline{\Gamma}^*(X,\mu) := R\overline{\Gamma}(X,\mu)R^{-1}$$
since $\rho(\Aff^+(X,\mu)) \cong \Stab(\Delta_q)\leq\Mod_g$ ($\Mod_g$ acts on the fibers of $\mathcal{T}_g$).\\
Since $\FakRaum{\mathds{H}}{\overline{\Gamma}^*(X,\mu)}$ is a surface of finite type, it is an algebraic curve if and only if $\overline{\Gamma}^*(X,\mu)$ -- or, equivalently, $\overline{\Gamma}(X,\mu)$ -- is a lattice in $\PSL_2(\mathds{R})$.\\
Since the map $\FakRaum{\mathds H}{\overline{\Gamma}^*(X,\mu)}\rightarrow \mathcal{M}_g$ is birational on the image \cite[§2]{mcmullen}, $\pr(\Delta_q)$ is an algebraic curve if and only if $\FakRaum{\mathds H}{\overline{\Gamma}^*(X,\mu)}$ is one.
\end{Bew}

\chapter{Schottky spaces}
\label{schottky}
\section{Schottky groups}
\begin{definition}
A group $\Gamma\leq \PSL_2(\mathds{C})$ of Möbius transformations is
called \begriff{Schottky group} if there is a $g\geq 1$ and disjoint,
simply connected domains $D_1,D'_1,\dots,D_g,D'_g$ bounded by disjoint Jordan curves $C_i:= \partial D_i, C'_i:= \partial D'_i$, and generators $\gamma_1,\dots,\gamma_g$ of $\Gamma$ such that $\gamma_i(C_i)=C'_i$ and $\gamma_i(D_i)=\mathds{P}^1(\mathds{C}) \backslash \overline{D'_i}$.\\
Such generators define a so called \begriff{Schottky base}.
\begin{figure}[h]
\begin{center}
\includegraphics[scale=0.45]{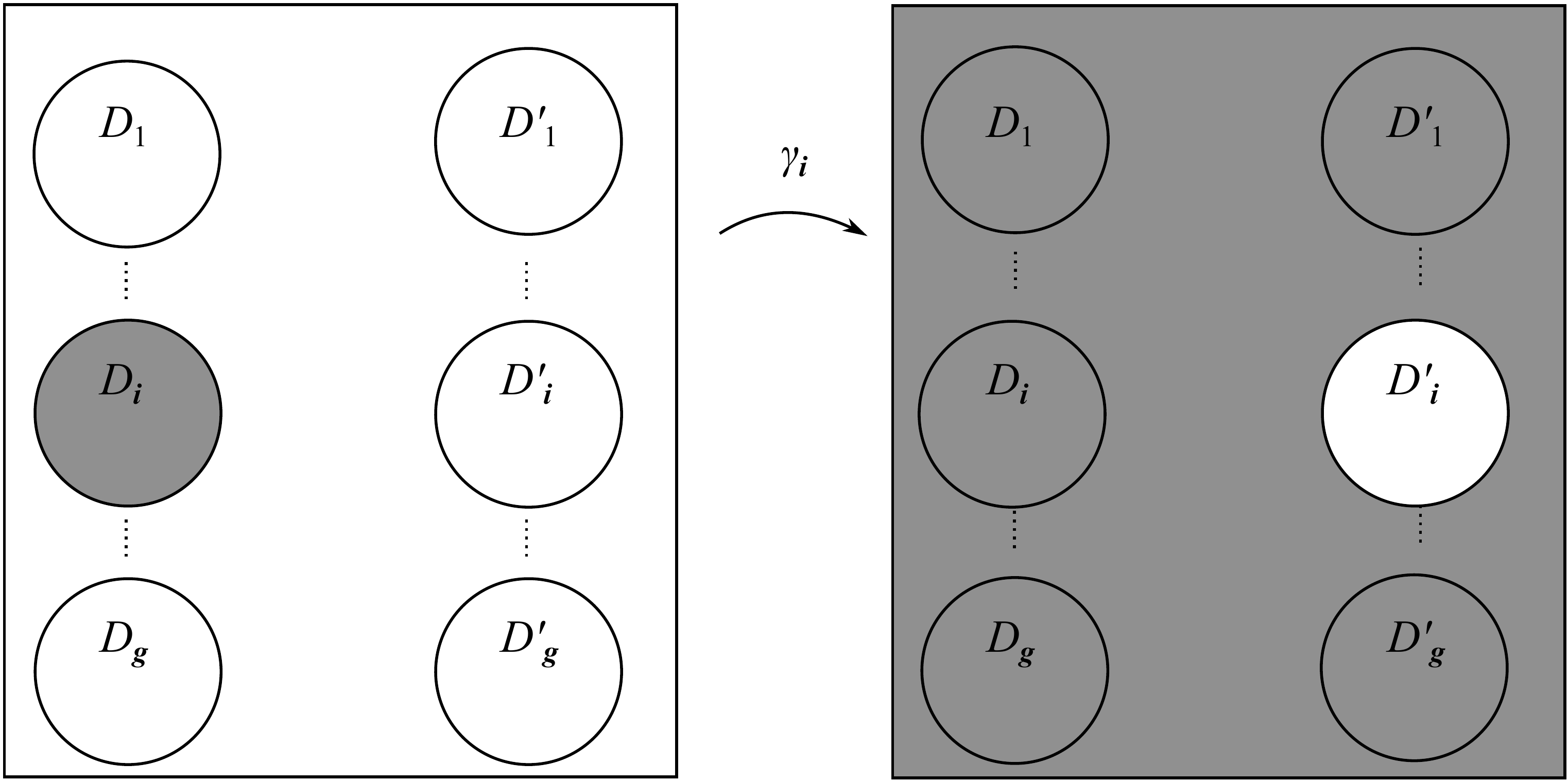}
\end{center}
\caption{A generator of a Schottky group}
\label{Schottkygruppe}
\end{figure}
\end{definition}

\begin{satz}
\label{schottkygruppe}
Let $(\gamma_1,\dots,\gamma_g)$ be a Schottky base of a Schottky group $\Gamma$ and take $D_i,D'_i,C_i,C'_i$ for $i\in\{1,\dots,g\}$ as above. Moreover, let $A:=\mathds{P}^1(\mathds{C}) \backslash \bigcup_{i=1}^g\overline{D_i\cup D'_i}$. We then have:
\begin{enumerate}
\item $\gamma_i$ is loxodromic.
\item $\gamma_i$ has its attracting fixed point in $D'_i$ and its repelling fixed point in $D_i$.
\item $\Gamma$ is free with generators $\gamma_1,\dots,\gamma_g$.
\item $\Gamma$ is Kleinian.
\item $\overline{A}$ is a fundamental domain of $\Gamma$.
\end{enumerate}
\end{satz}

\begin{Bew}~
Let $A_i:=\mathds{P}^1(\mathds{C}) \backslash (\overline{D_i\cup D'_i})$ for $i\in\{1,\dots,g\}$.
\begin{enumerate}
\item Assume that $\gamma_i$ is parabolic. We can then assume that $\gamma_i\colon z\mapsto z+1$.\\
We first show that $\gamma_i(\overline{A_i})\cap\overline{A_i}=C'_i$ and deduce $\infty\notin\overline{A_i}$ afterwards:\\
We have $\gamma_i(A_i)\cap A_i=\emptyset$ since $A_i\subset\mathds{P}^1(\mathds{C})\backslash\overline{D_i}$ and so \[\gamma_i(A_i)\subset\gamma_i(\mathds{P}^1(\mathds{C}) \backslash\overline{D_i})=D'_i\subset\mathds{P}^1(\mathds{C})\backslash A_i.\]
Furthermore, we have $\gamma_i(\partial A_i)\cap\partial A_i=C_i'$ because $C_i'\subseteq\gamma_i(\partial A_i)\cap\partial A_i$ since $\partial A_i=C_i\cup C_i'$ and $\gamma_i(C_i)=C_i'$. Moreover,
\[\gamma_i(C_i')\subset\gamma_i(\mathds{P}^1(\mathds{C})\backslash\overline{D_i})=D'_i\subset\mathds{P}^1(\mathds{C})\backslash\partial A_i\]
because $C'_i\subset\mathds{P}^1(\mathds{C}) \backslash\overline{D_i}$. It follows that $\gamma_i(\partial A_i)\cap\partial A_i\subseteq C_i'$ and, therefore, $\gamma_i(\overline{A_i})\cap\overline{A_i}=C'_i$.\\
If $\infty$  was in $\overline{A_i}$, we would have $\infty\in C'_i$ and $\infty\in C_i$ which contradicts to $C_i\cap C'_i=\emptyset$.\\
It follows that $\infty\in D_i\cup D'_i$. Let $\infty\in D'_i$. Then $D_i$ is bounded and there is, hence, some $z_0\in \overline{D_i}$ such that
$$\Real z_0\leq\Real z ~~ \forall z\in D_i.$$
That means $z_0\notin\gamma_i(D_i)$. On the other hand, we have $\overline{D_i}\subseteq\mathds{P}^1(\mathds{C}) \backslash \overline{D'_i}=\gamma_i(D_i)$, a \hbox{contradiction.}\\
If $\infty\in D_i$, then $D'_i$ is bounded, and there is a $z_0\in\overline{D'_i}$ such that
\[\Real z_0\geq\Real z ~~ \forall z\in \overline{D'_i}.\]
That means $z_0\notin\gamma_i^{-1}(\overline{D'_i})$. On the other hand, we have $\overline{D'_i}\subseteq\mathds{P}^1(\mathds{C}) \backslash D_i=\gamma^{-1}_i(\overline{D'_i})$, a contradiction.\\
Assume that $\gamma_i$ is elliptic: Then we can assume that $\gamma_i\colon z\mapsto e^{i\theta} z$. We have
$$\begin{array}{l}
\gamma_i(D_i)=\mathds{P}^1(\mathds{C}) \backslash \overline{D'_i} \Leftrightarrow\\ \gamma_i(\mathds{P}^1(\mathds{C}) \backslash D_i)=\overline{D'_i} \Rightarrow\\
\gamma_i(\overline{A})\subset \overline{D'_i} \text{ and } \gamma_i^n(\overline{A})\subset \overline{D'_i}~\forall n\in\mathds N.
\end{array}$$
On the other hand, let $z_0\in A$, and let $U(z_0)\subseteq A$ be an open neighbourhood of $z_0$. There is an $N\in\mathds N$ such that $\gamma_i^N(z_0)=e^{iN\theta}\cdot z_0\in U(z_0)\Rightarrow\gamma_i^N(z_0)\notin \overline{D'_i}$, which is a contradiction.
It follows that $\gamma_i$ is loxodromic.
\item W.l.o.g. let $\gamma_i\colon z\mapsto k\cdot z,~~\vert k \vert >1$. Then $\infty$ is the attracting and 0 the repelling fixed point. $\gamma_i(A_i)\cap A_i=\emptyset=\gamma_i(C_i)\cap C_i$, and so $0, \infty\notin \overline{A_i}$. Since $\gamma_i^n(D'_i)\subset D'_i$ for all $n\in\mathds N$ and since $\gamma_i^n(z)\rightarrow\infty$ for $z\in\mathds P^1(\mathds C)\backslash\{0\}$, it turns out that $D'_i$ cannot be bounded, and so $\infty\in D'_i$. We assume $0\in D'_i$. Since $D'_i$ is open, we have $B_{\varepsilon}(0)\subseteq D'_i$ for some $\varepsilon>0$. Then $\gamma_i^n(B_{\varepsilon}(0))=B_{\varepsilon\cdot\vert k \vert^n}(0)\subseteq D'_i~\forall n\in\mathds N$. Hence, $D'_i\supseteq\mathds C$, a contradiction. It follows that $0\in D_i$.
\item Let $F_g=\langle\overline{\gamma}_1,\dots,\overline{\gamma}_g\rangle$ be the free group and $\Psi\colon F_g\ra\Gamma$ the homomorphism defined by $\overline{\gamma}_i\mapsto\gamma_i$ for $i\in\{1,\dots,g\}$. We want to show that $\Psi$ is injective.\\
Let $\overline{\gamma}\in F_g\setminus\{1\}$ be given as reduced word, i.e.\ of the form
\hbox{$\overline{\gamma}=\overline{\gamma}_{i_1}^{j_1}\cdots\overline{\gamma}_{i_n}^{j_n}$} for $n\in\mathds N$, $i_1,\dots,i_n\in\{1,\dots,g\}$ and $j_1,\dots,j_n\in\mathds Z\backslash \{0\}$ and $i_k\neq i_{k+1}$ for $k\in\{1,\dots,n-1\}$.\\
Let $\gamma=\Psi(\overline{\gamma})$. We must show that $\gamma\neq\id$. We have 
\[\gamma_i(A)\subset\gamma_i(A_i)\subset\gamma_i(\mathds P^1(\mathds C)\backslash\overline{D_i})=D'_i\quad\text{ and }\quad\gamma_i(D'_i)\subset\gamma_i(\mathds P^1(\mathds C)\backslash\overline{D_i})= D'_i\] and analoguously $\gamma_i^{-1}(A)\subset\gamma_i^{-1}(A_i)\subset D_i$ and $\gamma_i^{-1}(D_i)\subset D_i$. It, thus, follows that
\begin{align*}
\gamma(A)&= \gamma_{i_1}^{j_1}\circ\dots\circ\gamma_{i_n}^{j_n}(A)\\
&\subset\gamma_{i_1}^{j_1}\circ\dots\circ\gamma_{i_{n-1}}^{j_{n-1}}(D_{i_n}\cup D'_{i_n})\\
&\subset\gamma_{i_1}^{j_1}\circ\dots\circ\gamma_{i_{n-1}}^{j_{n-1}}(A_{i_{n-1}})\\
&\subset\gamma_{i_1}^{j_1}\circ\dots\circ\gamma_{i_{n-2}}^{j_{n-2}}(D_{i_{n-1}}\cup D'_{i_{n-1}})\subset\dots\subset(D_{i_1}\cup D'_{i_1}).
\end{align*}
So we have $\gamma(A)\cap A=\emptyset$, and this implies $\gamma\neq\id$.
\item By the proof of c), no element in $\bigcup_{\gamma\in\Gamma}\gamma(A)$ can be a limit point and, therefore, $\Gamma$ is Kleinian. Moreover, no element in $\bigcup_{\gamma\in\Gamma}\gamma(\overline{A})$ can be a limit point: Therefore, let $z\in\partial\gamma(A)$ for some $\gamma\in\Gamma$. Then we have $z\in\gamma(C_i)$ or $z\in\gamma(C'_i)$, and so $z$ has an open neighbourhood intersecting $\gamma(A)$, $\partial\gamma(A)$ and $\gamma\circ\gamma_i^{-1}(A)$ respectively $\gamma\circ\gamma_i(A)$, lying in their union and containing other point of the orbit of $z$.
\item We already know that $\gamma(A)\cap A=\emptyset~\forall\gamma\in\Gamma\backslash\{\id\}$. It remains to show that $\bigcup_{\gamma\in\Gamma}\gamma(\overline{A})$ is the whole region of discontinuity. In other words:
$$\mathds P^1(\mathds C)\backslash\bigcup_{\gamma\in\Gamma}\gamma(\overline{A})=\Lambda:=\{\text{limit points of }\Gamma\}.$$
By the proof of d), we only need to show the inclusion "$\subseteq$":\\
W.l.o.g. let $\infty\in A$ and $z\in \mathds P^1(\mathds C)\backslash\bigcup_{\gamma\in\Gamma}\gamma(\overline{A})$. Then we have the following (instead of some $D'_i$, we could also have $D_i$, in this case, replace $\gamma_i$ by $\gamma_i^{-1}$):
\[\begin{array}{l}
z\notin\overline{A}\Rightarrow z\in D'_{i_1}\text{ for an } i_1\in\{1,\dots,g\}\\
z\notin\gamma_{i_1}(\overline{A})\Rightarrow z\in\gamma_{i_1}(D'_{i_2})\subset D'_{i_1}\text{ for an } i_2\in\{1,\dots,g\}\\
z\notin\gamma_{i_1}\circ\gamma_{i_2}(\overline{A})\Rightarrow z\in\gamma_{i_1}\circ\gamma_{i_2}(D'_{i_3})\subset\gamma_{i_1}(D'_{i_2})\subset D'_{i_1}\text{ for an } i_3\in\{1,\dots,g\} \text{ etc.}
\end{array}\]
For some sequence $(\gamma_{i_j}^{e_j})_{j\in\N}$ with $i_j\in\{1,\dots,g\}$ and $e_j\in\{-1,1\}$, we get a sequence of domains
\[G_n:=\gamma_{i_1}^{e_1}\circ\dots\circ\gamma_{i_n}^{e_n}(D'_{i_{n+1}}) \text{ for } e_{n+1}=1,\]
\[\text{respectively } G_n:=\gamma_{i_1}^{e_1}\circ\dots\circ\gamma_{i_n}^{e_n}(D_{i_{n+1}}) \text{ for } e_{n+1}=-1\]
with $G_{n+1}\subset G_n \text{ and } z\in G_n,~\forall n\in\mathds N$, and we get (as boundary curves of the domains $G_n$) a sequence of Jordan curves
\[K_n:=\gamma_{i_1}^{e_1}\circ\dots\circ\gamma_{i_n}^{e_n}(C'_{i_{n+1}}) \text{ for } e_{n+1}=1,\]
\[\text{respectively } K_n:=\gamma_{i_1}^{e_1}\circ\dots\circ\gamma_{i_n}^{e_n}(C_{i_{n+1}}) \text{ for } e_{n+1}=-1;~~e_i\in\{\pm 1\}.\]
Let $G:=\bigcap_{n=1}^\infty\overline{G_n}$. We want to show that $G=\{z\}$, which obviously implies $z\in\Lambda$.\\
$G$ is bounded and closed and, therefore, compact. So there exists an
\[r:=\max\{\vert x-y\vert: x, y\in G\}.\]
Assume that $r\neq 0$:
\begin{figure}[h]
\begin{center}
\includegraphics[scale=0.45]{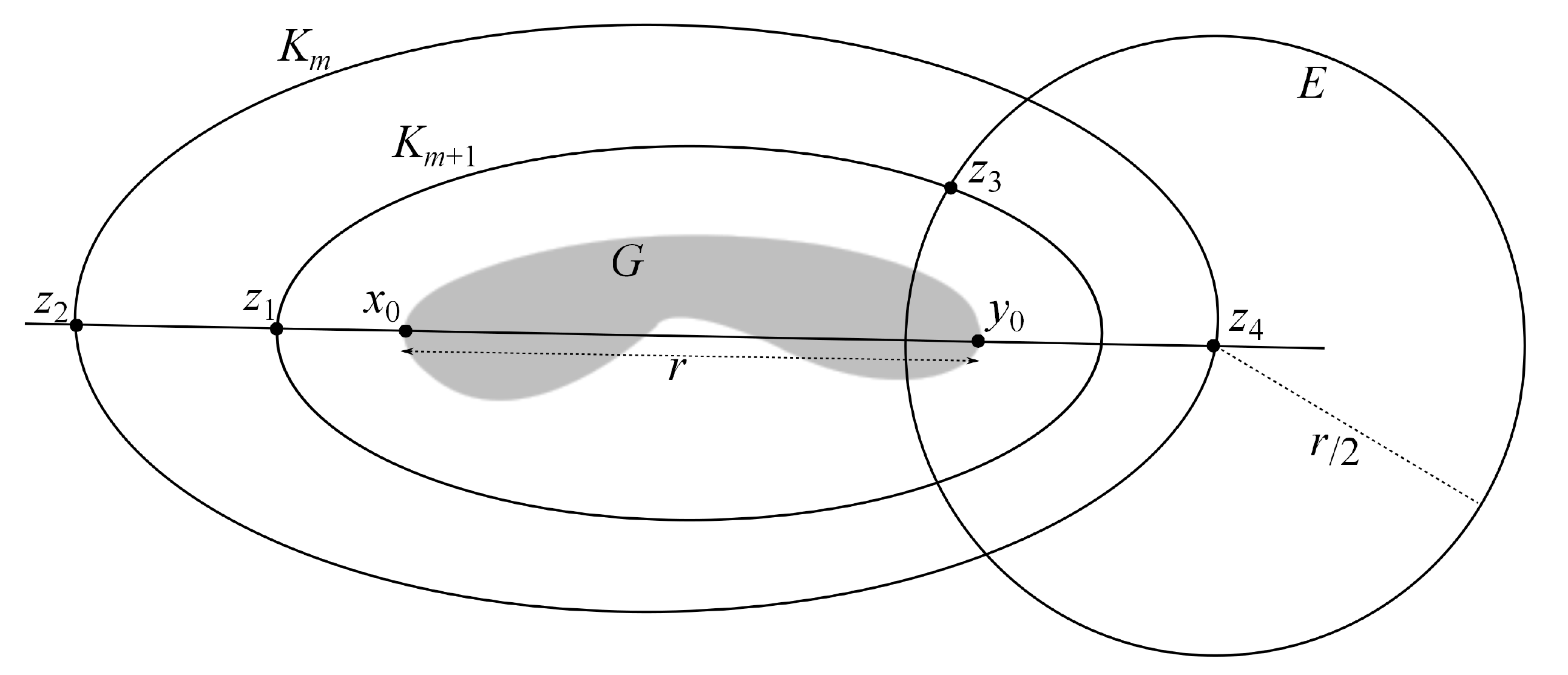}
\end{center}
\end{figure}\\
Let $x_0,y_0\in G$ with $\vert x_0 - y_0\vert =r$. Then for every $0<\varepsilon<\frac{r}{2}$ there exist two consecutive elements $K_m$ and $K_{m+1}$ of the sequence and points $z_1, z_2$ on the ray $x_0y_0^-$ and $z_4$ on the ray $y_0x_0^-$ with $z_1\in K_{m+1}$ and $z_2,z_4\in K_m$ with the following properties:
$$\vert z_1-z_2\vert\leq\varepsilon,~~~r\leq\vert z_1-z_4\vert\leq r+\varepsilon$$
Let us additionally consider the disc $E:=\{z\in\mathds C: \vert z-z_4\vert \leq\frac{r}{2}\}$. Since $y_0\in E$ and $x_0\notin E$, $\partial E$ intersects $K_{m+1}$. So there exists $z_3\in K_{m+1}$ with 
$$\vert z_3-z_4\vert=\frac{r}{2}.$$
It follows that
$$\vert z_1-z_3\vert\geq\vert z_1-z_4\vert-\vert z_3-z_4\vert\geq r-\frac{r}{2},$$
and this implies
$$\vert z_2-z_3\vert\geq\vert z_1-z_3\vert-\vert z_1-z_2\vert\geq \frac{r}{2}-\varepsilon.$$
So we conclude
\begin{equation}
\label{dv}
\vert DV(z_1,z_2,z_3,z_4)\vert=\frac{\vert z_1-z_2\vert\cdot\vert z_3-z_4\vert}{\vert z_2-z_3\vert\cdot\vert z_1-z_4\vert}\leq\frac{\varepsilon\cdot\frac{r}{2}}{(\frac{r}{2}-\varepsilon)\cdot r}\stackrel{\varepsilon\rightarrow 0}{\rightarrow}0.
\end{equation}
Now, let $\gamma:=(\gamma_{i_1}^{e_1}\circ\dots\circ\gamma_{i_m}^{e_m})^{-1}$. Then we have on the one hand
$$\vert DV(z_1,z_2,z_3,z_4)\vert=\vert DV(\gamma(z_1),\gamma(z_2),\gamma(z_3),\gamma(z_4))\vert,$$
since $\gamma$ is a Möbius transformation and, therefore, does not change the cross-ratio. On the other hand, we assumed $\infty\in A$, and, therefore,
$$0<\inf\{\vert w_i-w_j\vert : w_i\in B_i\neq B_j\ni w_j; B_i\in\{C_i,C'_i\}, B_j\in\{C_j, C'_j\}\}<\sup\{\dots\}<\infty,$$
and so we have
$$d:=\inf\{\vert DV(w_1,w_2,w_3,w_4)\vert:w_1,w_3\in B_i\neq B_j\ni w_2,w_4;\}>0.$$
Since $\gamma(K_{m+1})=C_{i_{m+2}}$ or $\gamma(K_{m+1})=C'_{i_{m+2}}$ and since $\gamma(K_m)$ is a component of $\partial A$ that is disjoint to $\gamma(K_{m+1})$ (either $C'_{i_{m+1}}$ or $C_{i_{m+1}}$), we have
\[\vert DV(\gamma(z_1),\gamma(z_2),\gamma(z_3),\gamma(z_4))\vert\geq d,\]
a contradiction to (\ref{dv}). It follows that $r=0$, and so $z$ is the only point in $G$. Therefore, $z\in\Lambda$.
\end{enumerate}
\end{Bew}

\begin{bem}
Thm.\ \ref{schottkygruppe} c) is a special case of the so called Ping-Pong-Lemma, see for example \hbox{\cite[Lemma 3.15]{farbmarg}.}
\end{bem}

\begin{bem}
Let $\Omega:=\bigcup_{\gamma\in\Gamma}\gamma(\overline{A})$ be the region of discontinuity of $\Gamma$. Then $X=\Omega/\Gamma$ is a compact Riemann surface of genus $g$.
\end{bem}

\begin{Bew}
The fundamental domain $A$ of $\Gamma$ is a $2g$-connected domain in $\mathds P^1(\C)$. The boundary curves $C_i$ and $C_i'$ are identified by $\Gamma$. So $X$ is topologically a sphere with $g$ handles. Since the group $\Gamma$ acts discontinuously on $\mathds P^1(\C)$ via Möbius transformations, $X$ inherits the complex structure of $\Omega$.
\end{Bew}

\begin{definition}
The quotient map $\Omega\overset{/\Gamma}{\rightarrow} X$ is called \begriff{Schottky covering}.
\end{definition}

\begin{satz}
\label{sub}
Every compact Riemann surface $X$ of genus $g\geq 1$ has such a Schottky covering.
\end{satz}
The following sketch of proof of this classical result is taken from \cite[5.2]{frankgabi}.
\begin{Bew}
Take disjoint, simply closed curves $c_1,\dots,c_g$ in $X$ such that $A:=X\setminus\bigcup_{i=1}^g c_i$ is connected.
Then $A$ is conformal equivalent to a domain in $\mathds{P}^1(\mathds{C})$ that is bounded by $2g$ closed curves, see \hbox{\cite[IV, 19F-G]{ahl}.}
Let $C_i$ and $C'_i$ for $i\in\{1,\dots,g\}$ be the boundary curves obtained by cutting the surface along the curves $c_i$. Now, let $F_g$ be the free group with generators $\gamma_1,\dots,\gamma_g$. Take a copy $A_w$ of $A$ for every reduced word $w\in F_g$. Now, we glue them together as follows: Identify $C_i$ on $A_{w'}$ with $C'_i$ on $A_w$ whenever $w=w'\gamma_i$. We get a domain $\Omega\subseteq\mathds{C}$, on which $F_g$ acts. This action extends to $\mathds{P}^1(\mathds{C})$, and $F_g$, thus, acts via Möbius transformations, see \hbox{\cite[IV, 19F]{ahl}.}
\end{Bew}

\section{The Schottky space and its complex structure}
\begin{definition}
Let $\tilde{\mathcal{S}}_g$ be the set of all $(\gamma_1,\dots,\gamma_g) \in\PSL_2(\mathds{C})^g$ that generate a Schottky group $\Gamma$ and are a Schottky base.
Take the equivalence classes under simultaneous conjugation:
$$\mathcal{S}_g:=\{(\gamma_1,\dots,\gamma_g)\}/\sim$$
$\mathcal{S}_g$ is called the \begriff{Schottky space} for Riemann surfaces of genus $g$.
\end{definition}
Let $s\in\tilde{\mathcal{S}}_g$, and let $\Gamma(s)$, $\Omega(s)$ and $X(s):= \Omega(s)/\Gamma(s)$ be the corresponding Schottky group, region of discontinuity and Riemann surface. This leads to an alternative definition of the Schottky space:
\begin{definition}
Take the pairs $(X,\sigma)$ such that
\begin{itemize}
\item[i)] $\sigma\colon F_g\rightarrow\PSL_2(\mathds{C})$ is an injective homomorphism,
\item[ii)] $\Gamma:=\sigma(F_g)$ is a Schottky group,
\item[iii)] $\Omega(\Gamma)/\Gamma\cong X.$
\end{itemize}
Then let $\mathcal{S}_g:=\{(X,\sigma)\}/\sim$ with
\[(X,\sigma)\sim (X',\sigma') \Leftrightarrow \exists A \in\PSL_2(\mathds{C})\colon \sigma'(\gamma)=A\sigma(\gamma)A^{-1}~~\forall \gamma\in F_g.\]
\end{definition}
Note that $X$ and $X'$ are isomorphic. That means that the \begriff{forgetful map}
\[\begin{array}{rrl}
\mu:&\mathcal{S}_g&\rightarrow\mathcal{M}_g\\
&[(X,\sigma)]&\mapsto[X]
\end{array}\]
is well defined.

\begin{definition}
Let $\gamma\in\PSL_2(\mathds{C})$ be loxodromic, i.e.\ conjugate to a Möbius transformation of the form $z\mapsto\lambda\cdot z$ with $0<\vert\lambda\vert<1$. Then we call $\lambda$ the \begriff{multiplier} of $\gamma$. It is well defined.
\end{definition}

\begin{satz}
\label{sraumoffeninc}
\begin{enumerate}
\item The map
$$\begin{array}{rcl}
\tilde{\mathcal{S}}_g&\hookrightarrow & (\mathds{P}^1(\mathds{C}))^{3g}\\
(\gamma_1,\dots,\gamma_g)&\mapsto &(z_1,w_1,\lambda_1,\dots,z_g,w_g,\lambda_g)
\end{array}$$
is an open embedding. Here we assume that $z_i$ is the attracting fixed point, $w_i$ is the repelling fixed point and $\lambda_i$ is the multiplier of $\gamma_i$ for $i\in\{1,\dots,g\}$. 
\item $S_1$ is a punctured disc.
\item For $g\geq2$, $\mathcal{S}_g$ is homeomorphic to an open subset of $\mathds{C}^{3g-3}$.
\end{enumerate}
\end{satz}
\begin{Bew}
\begin{enumerate}
\item The injectivity follows from the fact that a Möbius transformation is uniquely determined by its attracting fixed point, its repelling fixed point and its multiplier.\\
It remains to show that the image is open:\\
Let a triple $(z_i,w_i,\lambda_i)$ be given for one of the generators
\[\gamma_i(z)=\frac{a_iz+b_i}{c_iz+d_i} ~ \text{ with } ~ a_id_i-b_ic_i=1.\]
We can assume w.l.o.g. that $z_i,w_i\neq\infty$ and that $c_i\in \mathds H$ since these conditions are easily accomplished by conjugation. Note that   $(a_i,b_i,c_i,d_i)\in \mathds C^4$ is then unique, too. By Thm. \ref{schottkygruppe}a, we have $z_i\neq w_i$. Moreover, we have
$$(z_i,w_i,\lambda_i)\in E:=\{(z,w,\lambda):z,w\in\mathds C, z\neq w, \lambda\in\mathds D\setminus\{0\}\}~~(\text{open in }\mathds C^3).$$
For the fixed points $v_1$ and $v_2$ of $\gamma_i$, we have 
$$v_1= \frac{(a_i-d_i)+\sqrt{(a_i+d_i)^2-4}}{2c_i} ~ \text{ and } ~ v_2= \frac{(a_i-d_i)-\sqrt{(a_i+d_i)^2-4}}{2c_i}.$$ 
One of them is $z_i$, the other one is $w_i$. For the multiplier $\lambda_i$ of $\gamma_i$, we get
$$\begin{array}{ll}
&\frac{1}{\sqrt{\lambda_i}}+\sqrt{\lambda_i}=a_i+d_i\\
\Rightarrow &\lambda_i-\sqrt{\lambda_i}\cdot(a_i+d_i)+1=0\\
\Rightarrow&\lambda_i = \left(\frac{a_i+d_i}{2}+\sqrt{\left(\frac{a_i+d_i}{2}\right)^2-1}~\right)^2 \text{ and } \frac{1}{\lambda_i}=\left(\frac{a_i+d_i}{2}-\sqrt{\left(\frac{a_i+d_i}{2}\right)^2-1}~\right)^2 \text{ or vice versa.}
\end{array}$$
Here, $\lambda_i$ is the number, with absolute value less than 1. The terms under the square root are not zero since $\vert\lambda_i\vert\neq 1$.\\
Now, let
$$r_i:= z_i + w_i,~~~~ s_i:=z_i\cdot w_i~~\text{ and }~~ t_i:=\lambda_i+\frac{1}{\lambda_i}.$$
The values $r_i,s_i,t_i$ are uniquely determined by $z_i,w_i,\lambda_i$.
We get
\begin{align*}
r_i&=\frac{a_i-d_i}{c_i}\\
s_i&=\frac{(a_i-d_i)^2-(a_i+d_i)^2+4}{4c_i^2}\\
t_i&=(a_i+d_i)^2-2
\end{align*}
It follows that
\[a_i=\frac{1}{2}\left(r_i\sqrt{\frac{t_i-2}{r_i^2-4s_i}}+\sqrt{t_i+2}\right),\;c_i=\sqrt{\frac{t_i-2}{r_i^2-4s_i}},\;d_i=\frac{1}{2}\left(-r_i\sqrt{\frac{t_i-2}{r_i^2-4s_i}}+\sqrt{t_i+2}\right).\]
The only remainig indeterminacy in the choice of $(a_i,c_i,d_i)$ is controlled by the two terms under the square root. But the choice of $\sqrt{\frac{t_i-2}{r_i^2-4s_i}}$ is determined by $c_i\in\mathds H$, and the choice of $\sqrt{t_i+2}$ corresponds to the choice between $\gamma_i$ and $\gamma_i^{-1}$. However, this is nothing but the choice of which fixed point should be the attracting one and which the repelling one and is, therefore, determined by $z_i$ and $w_i$. Because of $c_i\neq 0$, also $b_i$ is determined by $a_id_i-b_ic_i=1$. Since $\lambda_i\neq\pm1$, we have $t_i\pm2\neq 0$. Furthermore, we have
$$r_i^2-4s_i=(v_1+v_2)^2-4v_1v_2=(v_1-v_2)^2\neq 0.$$
Therefore, there exist open sets $U\subset E$ and $V\subset\mathds C\times\mathds H\times\mathds C$ with $(z_i,w_i,\lambda_i)\in U$ and $(a_i,c_i,d_i)\in V$ such that the map
$$\begin{array}{llll}
f:& U& \rightarrow& V\\
&(z,w,\lambda)&\mapsto&\left(\frac{1}{2}\left(r\sqrt{\frac{t-2}{r^2-4s}}+\sqrt{t+2}\right),\sqrt{\frac{t-2}{r^2-4s}},\frac{1}{2}\left(-r\sqrt{\frac{t-2}{r^2-4s}}+\sqrt{t+2}\right)\right)
\end{array}$$
is well-defined and continuous and such that $f(z_i,w_i,\lambda_i)=(a_i,c_i,d_i)$. Here, $r,s$ and $t$ are defined analogously to $r_i,s_i$ and $t_i$.\\
Now, take $D_i$ and $D'_i$ as in the definition of Schottky groups, i.e.\ $\gamma_i(D_i)=\mathds P^1(\mathds C)\backslash\overline{D'_i}$.\\
The map
$$\begin{array}{lrcl}
\hat{\gamma_i}:&\mathds C\times\mathds H\times\mathds C\times\mathds P^1(\mathds C)&\rightarrow & \mathds{P}^1(\mathds{C})\\
&(a,c,d,x)&\mapsto &\frac{ax+\frac{ad-1}{c}}{cx+d}
\end{array}$$
is also continuous, and, therefore, $\hat{f}:=\hat{\gamma_i}\circ (f\times\id)$ is also continuous in $U\times\mathds P^1(\mathds C)$. There is, hence, for all $\delta>0$ an $\varepsilon:=\varepsilon(z_i,w_i,\lambda_i,x)>0$ such that
$$\vert (z_i,w_i,\lambda_i,x)-(z'_i,w'_i,\lambda'_i,x')\vert<\varepsilon\Rightarrow \vert\hat{f}(z_i,w_i,\lambda_i,x)-\hat{f}(z'_i,w'_i,\lambda'_i,x')\vert<\delta.$$
By the compactness of $\overline{D_i}$, it follows that for all $\delta$ there is some $\varepsilon:=\varepsilon(z_i,w_i,\lambda_i,D_i)>0$ such that
\begin{equation}
\label{stetig}
\forall\delta>0,\exists\varepsilon>0:\vert (z_i,w_i,\lambda_i)-(z'_i,w'_i,\lambda'_i)\vert<\varepsilon\Rightarrow \hat{f}(z'_i,w'_i,\lambda'_i,\overline{D_i}) =\mathds P^1(\mathds C)\backslash D''_i,
\end{equation}
where $D''_i$ is a domain such that
$$z'\in D''_i\Rightarrow d(z', D'_i)<\delta.$$
Since $\overline{D_1},\overline{D'_1}, \dots, \overline{D_g},\overline{D'_g}$ are pairwise disjoint, also $\overline{D_1}, \overline{D'_1}, \dots, \overline{D_i}, \overline{D''_i}, \dots, \overline{D_g}, \overline{D'_g}$ are pairwise disjoint as long as $\delta$ is small enough.
That means that also $(\gamma_1,\dots, \gamma'_i, \dots, \gamma_g)$ with
\[\gamma'_i\colon z\mapsto \frac{a'z+\frac{a'd'-1}{c'}}{c'z+d'}\qquad\text{ for } (a',b',c'):=f(z_i',w_i',\lambda_i'),\quad z_i',w_i',\lambda_i'\text{ as in }(\ref{stetig})\]
is a Schottky base. The assertion follows.
\item Every equivalence class has exactly one Möbius transformation of the form $z\mapsto \lambda\cdot z$ with $\lambda\in\mathds{C},~0<|\lambda|<1$.
\item There is exactly one element in the equivalence class of $s\in\tilde{\mathcal{S}}_g$ such that for this element, $\gamma_1$ has attracting fixed point $0$ and repelling fixed point $\infty$ and $\gamma_2$ has attracting fixed point $1$. Since $(\gamma_1,\dots,\gamma_g)$ is a Schottky base, $\gamma_2,\dots,\gamma_g$ do not fix the point $\infty$.
\end{enumerate}
\end{Bew}

\section{Coarse and fine moduli space}
\begin{definition}
Take two analytic spaces $\mathcal{H}$ and $\mathcal{B}$ and a proper, analytic map
\[\pi \colon \mathcal{H} \rightarrow \mathcal{B}.\]
Then we call $\mathcal{H}$ a \begriff{family of Riemann surfaces} or simply a ``family`` of genus $g$ over $\mathcal{B}$ if the fiber $\pi^{-1}(b)$ is a compact Riemann surface of genus $g$ for all $b\in \mathcal{B}$.
\end{definition}

\begin{definition}
\label{fms}
$\mathcal F_g$ is called \begriff{fine moduli space} for Riemann surfaces of genus $g$ if the following conditions are satisfied:
\begin{enumerate}
\item There is a bijection
\[\mathcal F_g \overset{\sim}{\longrightarrow} \{\text{isomorphism classes of compact Riemann surfaces of genus } g\}\]
that we'll use to identify $\mathcal F_g$ with the set of this isomorphism classes.
\item There is a family $\mathcal{U} \rightarrow \mathcal F_g$ (called \begriff{universal family}) such that for every family $\mathcal{H}\rightarrow \mathcal{B}$, the map $s\colon \mathcal{B}\rightarrow \mathcal F_g,~~b\mapsto [\pi^{-1}(b)]$ is analytic and such that $\mathcal{H}\cong\mathcal{U}\times_{\mathcal F_g}\mathcal{B}$, i.e. that
\[
\begin{xy}
\xymatrix{\mathcal{U}\times_{\mathcal F_g}\mathcal{B}\cong\hspace{-0.9cm}&\mathcal H \ar[d]^\pi\ar[r]& \mathcal U\ar[d]&&\\
&\mathcal{B}\ar[r]^s& \mathcal F_g&&
}
\end{xy}
\]
is a pullback square.
\end{enumerate}
\end{definition}

\begin{definition}
$\mathcal G_g$ is called \begriff{coarse moduli space} for Riemann surfaces of genus $g$ if the following conditions are satisfied:
\begin{enumerate}
\item As above in Def. \ref{fms}a.
\item For every family $\pi\colon \mathcal{H} \rightarrow \mathcal{B}$, the map $s\colon \mathcal{B} \rightarrow \mathcal G_g, b\mapsto [\pi^{-1}(b)]$ is analytic.
\item For every $\mathcal G'_g$ that satisfies a) and b), there is an analytic map $\Psi\colon \mathcal G_g \rightarrow \mathcal G'_g$ such that the following diagram commutes:
\[
\begin{xy}
\xymatrix{& \mathcal{H}\ar[d]&\\
& \mathcal{B}\ar[rd]^s\ar[ld]_{s'}&\\
\mathcal G'_g && \mathcal G_g\ar[ll]^{\Psi }}
\end{xy}
\]
\end{enumerate}
\end{definition}

\begin{satz}
$\mathcal{M}_g$ is a coarse moduli space for Riemann surfaces of genus $g$, and $\mathcal{T}_g$ is a fine moduli space for marked Riemann surfaces of genus $g$.
\end{satz}

\begin{Bew}
See \cite[2A, 2C]{harris}.
\end{Bew}
The next two theorems and the definition follow \cite[5.6, 5.7, 5.8]{frankgabi}.
\begin{satz}
\label{vergissanalytisch}
The forgetful map $\mu:\mathcal{S}_g\rightarrow \mathcal{M}_g$ is analytic and surjective.
\end{satz}
\begin{Bew} 
The surjectivity follows from the fact that every compact Riemann surface has a Schottky covering.
The only thing that remains to show is that $\mu$ is analytic. Since $\mathcal{M}_g$ is a coarse moduli space for Riemann surfaces, it suffices to find a family $\pi:\mathcal{C}_g\rightarrow \mathcal{S}_g$ of Riemann surfaces over $\mathcal{S}_g$ that induces $\mu$ in the following sense: $\mu (s)$ is the equivalence class of the Riemann surface 
C$_g :=\pi^{-1}(s)\subset\mathcal{C}_g$.\\
Let
$$\Omega_g := \lbrace (s,z) \in \mathcal{S}_g \times \mathds{P}^1(\mathds{C}): z \in \Omega (s) \rbrace.$$
That is a complex manifold, on which $F_g$ acts holomorphically by
$\varphi(s,z) := (s,\sigma(\varphi)(z))$ where $\sigma\colon F_g\rightarrow\PSL_2(\mathds C)$ is the injective homomorphism that is well defined by a normalization condition as in Thm. \ref{sraumoffeninc}c.\\
For $s\in \mathcal{S}_g$, we identify $\Omega(s)$ with $\{s\}\times\Omega(s)\subset\Omega_g$. Then C$_g:=\Omega(s)/F_g$ is a compact Riemann surface of genus $g$. Since $F_g$ acts trivially on the first component of $\Omega_g$, the projection to the first component $\hat{\pi}_1\colon \Omega_g\rightarrow \mathcal{S}_g$ factors through the quotient space $\mathcal{C}_g:=\Omega_g/F_g$. The induced map $\pi:\mathcal C_g\rightarrow \mathcal{S}_g$ is the family we are looking for.
\end{Bew}
\begin{definition}
\begin{enumerate}
\item Let $\mathcal{U}\rightarrow S$ be an analytic map between complex manifolds, and let \linebreak $\Gamma\leq\Aut(\mathcal{U}/S)$ be a proper discontinuous subgroup. Then the quotient map $\mathcal{U}\rightarrow\mathcal{U}/\Gamma = \mathcal{C}$ is called a \begriff{Schottky covering} if the induced map $\mathcal{C}\rightarrow S$ is a family of Riemann surfaces and if the restriction of the quotient map $U_s\rightarrow C_s$ is a Schottky covering for every $s\in S$.
\item A Schottky covering $\mathcal{U}\rightarrow\mathcal{U}/\Gamma$ together with an equivalence class of isomorphisms \hbox{$\sigma\colon F_g \rightarrow\Gamma$} is called \begriff{Schottky structure} on $\mathcal{U}$. Here, we will consider $\sigma$ and $\sigma'$ as equivalent if they differ only by an inner automorphism of $F_g$.
\end{enumerate}
\end{definition}

\begin{satz}
$\mathcal{S}_g$ is a fine moduli space for Riemann surfaces of genus $g$ with Schottky structure.
\end{satz}
\begin{Bew} (sketch)
Let $\mathcal{C}/\mathcal{S}$ be a family of Riemann surfaces and ($\mathcal{U}\rightarrow\mathcal{U}/\Gamma=\mathcal{C}, \sigma\colon F_g \overset{\sim}{\rightarrow}\Gamma$) a Schottky structure on $\mathcal{C}$.
Then there is a map $f\colon \mathcal{S}\rightarrow \mathcal{S}_g$ that maps a point $x$ to the equivalence class of the Schottky covering $U_x\rightarrow C_x$, see the following diagram:
\[
\begin{xy}
\xymatrix{\mathcal{U}\ar[r]\ar[d]^{/\Gamma}&\Omega_g\ar[d]^{/F_g}\\
\mathcal{C}\ar[r]\ar[d]&\mathcal{C}_g\ar[d]\\
\mathcal{S}\ar[r]^f & \mathcal{S}_g\ar[d]^{\mu}\\
& \mathcal{M}_g
}
\end{xy}
\]
Furthermore, we have $\mathcal{C}=\mathcal{C}_g\times_{\mathcal{S}_g}\mathcal{S}$ and $\mathcal{U}=\Omega_g\times_{\mathcal{C}_g}\mathcal{C}=\Omega_g\times_{\mathcal{S}_g}\mathcal{S}$.\\
It remains to show that $f$ is analytic, see \cite[5.8]{frankgabi} and \cite[§3]{frankger}.
\end{Bew}

\section{Relation between Teichmüller and Schottky space}

We now construct a family of Riemann surfaces over $\mathcal T_g$. For this purpose, we work with the following alternative definition of Teichmüller space:

\begin{definition}
\label{teichneu}
For $g\geq2$, take the (abstract) fundamental group of a surface of genus $g$ given by
\[\pi_g:=\left\langle a_1,b_1, \dots ,a_g,b_g\mid\prod_{i=1}^g[a_i,b_i]\right\rangle.\]
Moreover, let $\widetilde{\mathcal T_g}$ be the set of pairs $(X,\tau)$ where
\begin{itemize}
\item $\tau\colon \pi_g\hookrightarrow\PSL_2(\R)$ is an injective group homomorphism,
\item $\tau(\pi_g)$ is a cocompact Fuchsian group, and
\item $X$ is the Riemann surface $\FakRaum{\mathds H}{\tau(\pi_g)}$.
\end{itemize}
On $\widetilde{\mathcal T_g}$ we define the equivalence relation:
\[(X,\tau)\sim(X',\tau'):\Leftrightarrow\exists\gamma\in\PSL_2(\R)\colon c_\gamma\circ\tau=\tau'\]
where $c_\gamma$ denotes conjugation with $\gamma$:
\begin{align*}
c_\gamma\colon \PSL_2(\R)&\ra\PSL_2(\R)\\
A&\mapsto\gamma A\gamma^{-1}
\end{align*}
Then we define $\mathcal T_g:= \FakRaum{\widetilde{\mathcal T_g}}{\sim}$.
\end{definition}

\begin{bem}
\label{normteich}
In every equivalence class $\FakRaum{(\widetilde{X},\widetilde{\tau})}{\sim}$ we find exactly one representative $(X,\tau)$ such that $\tau(b_1)$ has attracting fixed point $0$ and repelling fixed point $\infty$ and $\tau(b_2)$ has attracting fixed point $1$.\\
From now on, we will identify $\mathcal T_g$ with the set of such representatives.
\end{bem}
This definition of $\mathcal T_g$ is equivalent to Definition \ref{teich}. We can see this by giving another definition of Teichmüller space, where we take a set of symplectic generators of the fundamental group up to an isotopic shift of the base point as markings (in particular, we look at the set symplectic generators up to an inner automorphism). Two marked riemann surfaces define the same point if there is a biholomorphic map between them that respects these markings, see \cite[§ 1.3]{imtan}. We'll denote this space by $\mathcal T'_g$. We show $\mathcal T'_g=\mathcal T_g$.\\
Let $[(X,\Sigma)]\in\mathcal T'_g$ where $\Sigma$ is a symplectic set of generators of $\pi_1(X)$.
Since we assume $g\geq2$, we can, by the uniformisation theorem, write $X$ as $\FakRaum{\mathds H}{\Gamma}$ for a Fuchsian group $\Gamma$ that is isomorphic to $\pi_g$. $\Gamma$ is uniquely determined by $X$ up to conjugation in $\PSL_2(\R)$. The choice of $\Sigma$ up to equivalence corresponds to the choice of a set of generators of the deck transformation group $\Gamma$ up to an inner automorphism, i.e.\ the choice of some $\tau\colon \pi_g\stackrel{\sim}{\ra}\Gamma$ up to conjugation. Thus, we have a bijection between the set of points $[(X,\Sigma)]\in\mathcal T'_g$ and the set of points $[(X,\tau)]\in\mathcal T_g$. For a bijection between $\mathcal T'_g$ and $\mathcal T_g$ from Definition \ref{teich}, see \cite[Thm. 1.4]{imtan}. An important ingredient to the proof in \cite{imtan} is the theorem of Dehn-Nielsen-Baer that will be of further interest for us.

\begin{satz}[Dehn-Nielsen-Baer]~\\
Let $X$ be a Riemann surface of genus $g$ and $p\in X$. Then we have:
\begin{align*}
&&\Mod_g&\ra\Out^+(\pi_1(X,p))\\
&&[f]&\mapsto[f_*],\qquad\qquad\qquad f\in\Homo^+(X)
\end{align*}
\end{satz}
\begin{Bew}
See \cite[Thm. 8.1]{farbmarg}.
\end{Bew}
Now, $\pi_g$ acts on $\mathcal T_g\times\mathds H$ via
$$\gamma\cdot((X,\tau), z):= ((X,\tau),\tau(\gamma)(z)), ~~~\text{ for } \gamma\in\pi_g.$$
We define $\mathcal C_{g,0}:= \FakRaum{\mathcal T_g\times\mathds H}{\pi_g}$.\\
Thus, $\mathcal C_{g,0}\ra\mathcal T_g$ is a family of Riemann surfaces over $\mathcal T_g$, and we get an isomorphism
\[\theta\colon \pi_g\ra\Aut((\mathcal T_g\times\mathds H)/\mathcal C_{g,0}).\]
Its restriction to the fiber $\{(X,\tau)\}\times\mathds H$ is exactly the isomorphism $$\theta_{(X,\tau)}:=\tau:\pi_g\ra\Aut(\mathds H/X),$$
and we, hence, get a so-called \begriff{Teichmüller structure} over $\mathcal{T}_g$:
$$\mathcal T_g\times\mathds H\overset{/\pi_g}{\longrightarrow}\mathcal{C}_{g,0} \longrightarrow \mathcal{T}_g$$
We now endow this family of Riemann surfaces with a Schottky structure:\\
First, we construct for every point $(X,\tau)\in\mathcal T_g$ a point in $\mathcal S_g$:\\
For this purpose, consider some $X = \mathds H/\tau(\pi_g)$ together with the isomorphism $\tau\colon \pi_g\ra\Deck(\mathds H/X)$ that is fixed by the normalization condition in Rem. \ref{normteich}. We can, therefore, view the elements of the image of $\tau$ as elements in $\pi_1(X)$, well defined up to an inner automorphism. So we get a symplectic representation of the fundamental group
\[\pi_1(X)=\left\langle\tau(a_1),\tau(b_1),\dots,\tau(a_g),\tau(b_g)\mid\prod_{i=1}^g[\tau(a_i),\tau(b_i)]\right\rangle.\]
As in Thm. \ref{sub}, the elements $\tau(a_1), \dots, \tau(a_g)$ give rise to a Schottky covering with region of discontinuity $\Omega_\tau$ and Schottky group $\Gamma_\tau(X)\leq\PSL_2(\C)$ with Schottky base (up to conjugation). We now fix an ordered set $(\gamma_1, \dots, \gamma_g)$ of free generators of $F_g$ and let $\sigma_{\alpha,\tau}\colon F_g\stackrel{\sim}{\ra}\Gamma_\tau(X)$ be the isomorphism that maps $(\gamma_1, \dots, \gamma_g)$ to the Schottky base $(\sigma_{\alpha,\tau}(\gamma_1),\dots,\sigma_{\alpha,\tau}(\gamma_g))$. Let the latter be normalized in such a way that $\sigma_{\alpha,\tau}(\gamma_1)$ has $0$ as attracting fixed point and $\infty$ as repelling fixed point and $\sigma_{\alpha,\tau}(\gamma_2)$ has $1$ as attracting fixed point. This Schottky base gives us the point $(X,\sigma_{\alpha,\tau})\in \mathcal S_g$.\\
Now, we want to construct Schottky covering with Schottky structure $\sigma_\alpha$ on a family of Riemann surfaces such that the restrictions of $\sigma_\alpha$ to the fibres is again the $\sigma_{\alpha,\tau}$ we started with. We define the following group homomorphism:
\begin{align*}
\alpha\colon\pi_g&\ra F_g\\
a_i&\mapsto 1,b_i\mapsto \gamma_i
\end{align*}
and denote its kernel by $N_\alpha$. Then $(\mathcal T_g\times\mathds H)/N_{\alpha}\ra\mathcal{C}_{g,0}$ is a Schottky covering and the universal Teichmüller structure $\theta$ induces a Schottky structure
\[\sigma_{\alpha}\colon \underbrace{\pi_g/N_{\alpha}}_{\cong F_g}\overset{\sim}{\longrightarrow}\Aut(((\mathcal T_g\times\mathds H)/N_{\alpha})/\mathcal{C}_{g,0})\]
on $\mathcal{C}_{g,0}$. The restriction to the fibre over $(X,\tau)$ is then
\[
\sigma_{\alpha,(X,\tau)}:\pi_g/N_{\alpha}\overset{\sim}{\rightarrow}\Aut((\mathds H/\tau(N_{\alpha}))/X)
\]
We may now choose the isomorphism
\[\widetilde{\sigma}_{\alpha,\tau}\colon \Aut((\mathds H/\tau(N_{\alpha}))/X)\ra\Gamma_\tau(X)\]
such that $\sigma_{\alpha,\tau}=\widetilde{\sigma}_{\alpha,\tau}\circ\sigma_{\alpha,(X,\tau)}$.\\
Since $\mathcal{S}_g$ is a fine moduli space, the following map is analytic:
\begin{align*}
s_{\alpha}:\quad\;\;\mathcal{T}_g&\rightarrow \mathcal{S}_g\\
(X,\tau)&\mapsto(X,\sigma_{\alpha,\tau})
\end{align*}

\begin{satz}
\label{teichquotabb}
\begin{enumerate}
\item
$s_{\alpha}$ is the quotient map for the following subgroup of $\Mod_g$:
\[\Mod_g(\alpha)=\{\varphi\in\Mod_g\mid\alpha\circ\tilde{\varphi}\equiv\alpha \operatorname{mod} \Inn(F_g) \text{ with } \tilde{\varphi}\in\Aut^+(\pi_g)\colon [\tilde{\varphi}]=\varphi\in\Out^+(\pi_g)\}\]
\item $s_{\alpha}\colon \mathcal{T}_g\rightarrow \mathcal{S}_g$ is the universal covering of the Schottky space.
\item $s_{\alpha}$ can be lifted to maps $\tilde{s}_{\alpha}$ and $\omega_{\alpha}$ such that the following diagram commutes:
\[
\begin{xy}
\xymatrix{\mathcal T_g\times \mathds H\ar[dd]_{/\pi_g}\ar[rd]^{/N_{\alpha}}&&\\
&\FakRaum{\mathcal T_g\times \mathds H}{N_{\alpha}}\ar[r]^{~~~~\omega_{\alpha}}\ar[dl]&\Omega_g\ar[d]^{/F_g}\\
\mathcal{C}_{g,0}\ar[d]\ar[rr]^{\tilde{s}_{\alpha}}&&\mathcal{C}_g\ar[d]\\
\mathcal{T}_g\ar[rr]^{s_{\alpha}}\ar[rd]_{/\Mod_g}&&\mathcal{S}_g\ar[dl]^{\mu}\\
&\mathcal{M}_g&}
\end{xy}
\]
\end{enumerate}
\end{satz}
\begin{Bew}
\begin{enumerate}
\item The group $\Mod_g(\alpha)$ is well defined:
Let $\tilde{\varphi}_1$ and $\tilde{\varphi}_2$ be representatives of $\varphi\in\Out^+(\pi_g)$. That means
\[\tilde{\varphi}_1=\gamma\tilde{\varphi}_2\gamma^{-1} \text{ for some } \gamma\in\pi_g.\]
Then we have
\[\alpha\circ\tilde{\varphi}_1= \alpha\circ(\gamma\tilde{\varphi}_2\gamma^{-1})\stackrel{\alpha \text{ hom.}}{=}\alpha(\gamma)(\alpha\circ\tilde{\varphi}_2)\alpha(\gamma)^{-1}.\]
It follows that
\[\alpha\circ\tilde{\varphi}_1\equiv \alpha\circ\tilde{\varphi}_2\operatorname{mod} \Inn(\pi_g),\]
and $\Mod_g(\alpha)$ is, hence, well-defined.\\\\
In order to prove the theorem, we first look at how $\Mod_g\cong\Out^+(\pi_g)$ acts on $\mathcal T_g$. Let $\varphi\in\Out^+(\pi_g)$ and $(X,\tau)\in\mathcal T_g$. We have
\[\varphi\cdot(X,\tau):=(X,\tau\circ\widetilde{\varphi}^{-1}),\]
where $\widetilde{\varphi}\in\Aut^+(\pi_g)$ is a representative of $\varphi$ such that $\tau\circ\widetilde{\varphi}^{-1}$ satisfies the normalization condition from definition \ref{teichneu} of $\mathcal T_g$.\\
We show that $s_{\alpha}((X,\tau))=s_{\alpha}(\varphi\cdot(X,\tau))$ if and only if $\varphi\in\Mod_g(\alpha)$.\\
Here, we have \[s_{\alpha}((X,\tau))=(X,\sigma_{\alpha,\tau})\text{ and }s_{\alpha}(\varphi\cdot(X,\tau))=(X,\sigma_{\alpha, \tau\circ\widetilde{\varphi}^{-1}}).\]
So we have
$$s_{\alpha}((X,\tau))=s_{\alpha}(\varphi\cdot(X,\tau))\Leftrightarrow(X,\sigma_{\alpha,\tau})\sim(X,\sigma_{\alpha, \tau\circ\widetilde{\varphi}^{-1}})\Leftrightarrow\sigma_{\alpha,\tau}=\sigma_{\alpha, \tau\circ\widetilde{\varphi}^{-1}}.$$\\
"'$\underline{\Rightarrow}$"': From the arguments above, one deduces that the images of $\sigma_{\alpha,\tau}$ and $\sigma_{\alpha, \tau\circ\widetilde{\varphi}^{-1}}$ are the same, and thus $\Gamma_\tau(X)=\Gamma_{\tau\circ\widetilde{\varphi}^{-1}}(X)$. For given $\tau$ and $\alpha$ (and the $\sigma_{\alpha,\tau}$ they induce), let us define $\beta_\tau\colon \Aut(\mathds H/X)\ra\Gamma_\tau(X)$ as the homomorphism renderin the following commutative diagram:
\[
\begin{xy}
\xymatrix{\pi_g\ar[d]_\alpha\ar[rr]^\tau&&\Aut(\mathds H/X)\ar[d]^{\beta_\tau}\\
F_g\ar[rr]^{\sigma_{\alpha,\tau}}&&\Gamma_\tau(X)}
\end{xy}
\]
We now have the following diagram:
\[
\begin{xy}
\xymatrix{&\pi_g\ar[dl]_{\tau}\ar[d]^{\alpha}\ar[r]^{\widetilde{\varphi}^{-1}}&\pi_g\ar@.[dl]_{\widetilde{\alpha}}\ar[d]^{\tau}\\
\Aut(\mathds H/X)\ar[d]_{\beta_\tau}&\pi_g/N_{\alpha}\ar[dl]^{\sigma_{\alpha,\tau}}\ar[rd]_{\sigma_{\alpha,\tau\circ\widetilde{\varphi}^{-1}}}&\Aut(\mathds H/X)\ar[d]^{\beta_{\tau\circ\widetilde{\varphi}^{-1}}}\\
\Gamma_\tau(X)\ar[rr]^\id&&\Gamma_{\tau\circ\widetilde{\varphi}^{-1}}(X)}
\end{xy}
\]
The diagram commutes without the dotted arrow because the parts on the top on the left and on the right commute by definition of $\beta_\tau$ and $\beta_{\tau\circ\widetilde{\varphi}^{-1}}$ and because the bottom part commutes by assumption. Because of the normalizing condition for $\tau$, $\tau\circ\widetilde{\varphi}^{-1}$, $\sigma_{\alpha,\tau}$ and $\sigma_{\alpha,\tau\circ\widetilde{\varphi}^{-1}}$ and because of $\Gamma_\tau(X)=\Gamma_{\tau\circ\widetilde{\varphi}^{-1}}(X)$, we also have $\beta_\tau=\beta_{\tau\circ\widetilde{\varphi}^{-1}}$.\\
We now define
\[\widetilde{\alpha}:=\sigma_{\alpha,\tau\circ\widetilde{\varphi}^{-1}}^{-1}\circ\beta_{\tau\circ\widetilde{\varphi}^{-1}}\circ\tau\]
and it follows that $\widetilde{\alpha}\circ\widetilde{\varphi}^{-1}=\alpha$, that is, the whole diagram commutes. Moreover, we have
\[\alpha=\sigma_{\alpha,\tau\circ\widetilde{\varphi}^{-1}}^{-1}\circ\id\circ\beta_\tau\circ\tau=\widetilde{\alpha}.\]
We may, thus, conclude $\alpha\circ\widetilde{\varphi}=\alpha$.\\\\
"'$\underline{\Leftarrow}$"' We assume that $\varphi\in\Mod_g(\alpha)$.\\
So for some $\widetilde{\varphi}\in\Aut^+(\pi_g)$ such that $[\tilde{\varphi}]=\varphi$, we have $\alpha\circ\widetilde{\varphi}\equiv\alpha\mod\Inn(F_g)$.\\
W.l.o.g. we can choose $\widetilde{\varphi}$ such that $\alpha\circ\widetilde{\varphi}=\alpha$. 
Now, $\tau$ gives us the Schottky cut system $(\tau(a_1),\dots,\tau(a_g))$ on $X$, and as in Thm. \ref{sub} we get a fundametal domain $\overline{A^\tau}\subset\mathds P^1(\C)$ of the Schottky group $\Gamma_\tau(X)$ that we indentify with the free group $F_g=\langle\gamma_1,\dots,\gamma_g\rangle$. As in \hbox{Thm. \ref{sub},} $F_g$ acts via Möbius transformations on the region of discontinuity $\Omega$ by sending a copy $\overline{A_w^\tau}$ of $\overline{A^\tau}$ to $\overline{A_{\gamma_iw}^\tau}$ with $\gamma_i$ (in the same way, $\tau\circ\widetilde{\varphi}^{-1}$ gives us a fundamental domain $\overline{A^{\tau\circ\widetilde{\varphi}^{-1}}}$). The (abstract) fundamental group $\pi_g$, therefore, also acts on $\Omega$ via $\delta\cdot\overline{A_w^\tau}=\overline{A_{\alpha(\delta)w}^\tau}$. This action coincides with the action of $\pi_g$ on $\Omega$ that is induced by $\tau\circ\widetilde{\varphi}^{-1}$, i.e.\ $\delta\cdot\overline{A_w^\tau}=\overline{A_{\alpha\circ\widetilde{\varphi}^{-1}(\delta)w}^\tau}$ since $\alpha\circ\widetilde{\varphi}^{-1}=\alpha$, by assumption. But this means nothing but that $\sigma_{\alpha,\tau\circ\widetilde{\varphi}^{-1}}\circ\alpha=\sigma_{\alpha,\tau}\circ\alpha$ 
and, therefore, $\sigma_{\alpha,\tau\circ\widetilde{\varphi}^{-1}}=\sigma_{\alpha,\tau}$.
\item This claim follows from the fact that $\mathcal{T}_g$ is simply connected and that $\Mod_g({\alpha})$ is torsion free, as we will see in Thm. \ref{torsionsfrei}. That means that $s_{\alpha}$ is unramified. An independent proof is in \cite[Thm. A and Rem. 5.1]{hejhal}.
\item Since $\FakRaum{\mathcal T_g\times \mathds H}{N_{\alpha}}\rightarrow \mathcal{C}_{g,0}$ is a Schottky covering and since $\mathcal{S}_g$ is a fine moduli space, we have: $\mathcal{C}_{g,0}=\mathcal{T}_g\times_{\mathcal{S}_g}\mathcal{C}_g$, and $\tilde{s}_{\alpha}$ is the projection to $\mathcal{C}_g$. Analogously, we construct $\omega_{\alpha}$.
\end{enumerate}
\end{Bew}

\section{Teichmüller discs in Schottky space}
We now look at Teichmüller discs whose image in $\mathcal{M}_g$ is an algebraic curve. We are interested in their image in Schottky space. An important fact in this context is Thm. \ref{affgamma}.
\begin{definition}
\label{sy}
A surjective homomorphism $\alpha\colon \pi_g\rightarrow F_g$ is called \begriff{symplectic} if there is a set of symplectic generators $(a_1, b_1,\dots, a_g, b_g)$ of $\pi_g$ with $\alpha(a_i)=1$ for $i\in\{1,\dots,g\}$. (With $\alpha(b_i)=:\gamma_i$, we then have that $(\gamma_1,\dots, \gamma_g)$ is a set of generators of the free group $F_g$.)
\end{definition}
\begin{bem}
\begin{enumerate}
\item Equivalently, we could require that $\alpha(b_i)=1$.
\item Let $M\subseteq\{1,\dots, g\}$. Then $\alpha$ is also symplectic if we require $\alpha(a_j)=1$ for $j\in M$ and $\alpha(b_k)=1$ for $k\notin M$, instead of what we required in Def. \ref{sy}.
\end{enumerate}
\end{bem}

\begin{Bew}
\begin{enumerate}
\item[(b)] Let $\hat{a_j}:=a_j, \hat{b_j}:=b_j, \hat{a_k}:=b_k^{-1}, \hat{b_k}:=b_ka_k$ for $j\in M$ and $k\notin M$. This is also a set of symplectic generators, and we have $\alpha(\hat{a_i})=1$ for all $i\in\{1,\dots,g\}$.
\item[(a)] This follows from (b).
\end{enumerate}
 \end{Bew}

\begin{satz}[Herrlich-Schmithüsen]~\\
\label{affgamma}
Let $\iota\colon \mathds{H} \rightarrow \mathcal{T}_g$ be a Teichmüller embedding such that the corresponding projective Veech group $\overline{\Gamma}_{\iota}\leq\PSL_2(\mathds{R})$ is a lattice and let $\Delta:=\iota(\mathds H)$. Then there is a symplectic homomorphism $\alpha:\pi_g \rightarrow F_g$ such that
$$\Stab(\Delta)\cap\Mod_g(\alpha)\neq\{1\}.$$
\end{satz}
\begin{Bew}
See \cite[Prop. 5.21]{frankgabi}. A constructive proof of a more general statement is given in \hbox{chapter \ref{kapori}} for origamis and in chapter \ref{kaphalbtrans} for flat surfaces in general.
\end{Bew}
Thm. \ref{affgamma} implies that the intersection is infinite since $\Mod_g(\alpha)$ is torsion free, see Thm. \ref{torsionsfrei}.\\
The restriction of $s_{\alpha}$ to $\Delta$ factorizes as follows:
\[\Delta\rightarrow\FakRaum{\Delta}{\Stab(\Delta)\cap\Mod_g(\alpha)}\rightarrow\Bild(\Delta)\subset \mathcal{S}_g,\]
and the image of the Teichmüller disc in the Schottky space, thus, cannot be isomorphic to a disc.\\
With remark \ref{affstab}, it follows that $\Stab(\Delta)\cap\Mod_g(\alpha)\cong\Aff^+(X)\cap\Mod_g(\alpha)$.

\begin{definition}
\label{mphi}
According to Dehn-Nielsen's theorem, an element $\varphi\in\Mod_g$ can be seen as an outer automorphism of the fundamental group $\pi_g:=\langle a_1,b_1,\dots,a_g,b_g\mid\Pi_{i=1}^g[a_i,b_i]=1\rangle$ and $\varphi$, hence, defines an automorphism of the abelianized fundamental group
\[\pi_g^{ab}:=\FakRaum{\pi_g}{[\pi_g,\pi_g]}\cong H_1(X,\Z)\cong\Z^{2g},\]
which is nothing but a matrix in $\Z^{2g\times2g}$:
\[M_\varphi:=\begin{pmatrix}
                         \sharp_{a_1}(\varphi(a_1))&\dots&\sharp_{a_1}(\varphi(a_g))&\sharp_{a_1}(\varphi(b_1))&\dots&\sharp_{a_1}(\varphi(b_g))\\
                         \vdots&&\vdots&\vdots&&\vdots\\                         
                         \sharp_{a_g}(\varphi(a_1))&\dots&\sharp_{a_g}(\varphi(a_g))&\sharp_{a_g}(\varphi(b_1))&\dots&\sharp_{a_g}(\varphi(b_g))\\
                         \sharp_{b_1}(\varphi(a_1))&\dots&\sharp_{b_1}(\varphi(a_g))&\sharp_{b_1}(\varphi(b_1))&\dots&\sharp_{b_1}(\varphi(b_g))\\
                         \vdots&&\vdots&\vdots&&\vdots\\
                         \sharp_{b_g}(\varphi(a_1))&\dots&\sharp_{b_g}(\varphi(a_g))&\sharp_{b_g}(\varphi(b_1))&\dots&\sharp_{b_g}(\varphi(b_g))\\
                        \end{pmatrix}.\]
\end{definition}
Here, we denote by $\sharp_{a_i}(\varphi(a_j))$ the sum of the exponents of $a_i$ in the word $\varphi(a_j)$. Indeed, $\varphi(a_j)$ is only well defined up to conjugation and modulo $\prod_{i=1}^g[a_i,b_i]$, but $\sharp_{a_i}(\varphi(a_j))$ is well-defined nevertheless. When we have $b_i$'s instead of $a_i$'s, the definition is completely analogous.
\begin{lemma}
\label{EW1}
Let $\alpha\colon \pi_g\ra F_g$ be the symplectic homomorphism with $\alpha(a_i)=1$ and $\alpha(b_i)=\gamma_i$ for $i\in\{1,\dots,g\}$ and $\varphi\in\Mod_g(\alpha)$. Then we have for $i,j\in\{1,\dots,g\}$:
\begin{enumerate}
\item $\sharp_{b_i}(\varphi(b_j)) = \delta_{ij}$,
\item $\sharp_{b_i}(\varphi(a_j)) = 0$,
\item $\sharp_{a_j}(\varphi(a_i)) = \delta_{ij}$,
\item $M_\varphi=\begin{pmatrix}I_g & A\\O_g&I_g\end{pmatrix}$,
with $I_g,O_g,A\in\Z^{g\times g}$. Here, $I_g$ and $O_g$ are the identity respectively the zero-matrix. In particular, $M_\varphi$ has eigenvalue 1 with (algebraic) multiplicity $2g$, and we have $\det(M_\varphi)=1$.
\end{enumerate}
\end{lemma}

\begin{Bew}
\begin{enumerate}
\item $\delta_{ij}=\sharp_{\gamma_i}(\alpha(b_j))\stackrel{\varphi\in\Mod_g(\alpha)}{=}\sharp_{\gamma_i}(\alpha\circ\varphi(b_j))\stackrel{\text{def. of }\alpha}{=}\sharp_{b_i}(\varphi(b_j))$.
\item $0=\sharp_{\gamma_i}(\alpha(a_j))=\sharp_{\gamma_i}(\alpha\circ\varphi(a_j))=\sharp_{b_i}(\varphi(a_j))$.
\item For $x,y\in\pi_g$, denote by $\hat{\text{i}}(x,y)$ the algebraic intersection number of $x$ and $y$, which remains invariant under isotopy and under homeomorphisms. Therefore, we have:
\begin{align*}
\delta_{ij}&=\hat{\text{i}}(a_i,b_j)\\
&=\hat{\text{i}}(\varphi(a_i),\varphi(b_j))\\
&=\sum_{k=1}^g\sharp_{a_k}(\varphi(a_i))\cdot\underbrace{\sharp_{b_k}(\varphi(b_j))}_{=\delta_{kj}}+\underbrace{\sharp_{b_k}(\varphi(a_i))}_{=0}\cdot\sharp_{a_k}(\varphi(b_j))\\
&=\sharp_{a_j}(\varphi(a_i)).
\end{align*}
\item This is only a)-c) in matrix notation.
\end{enumerate}
\end{Bew}

\begin{satz}
\label{torsionsfrei}
The group $\Mod_g(\alpha)$ is torsion free.
\end{satz}

\begin{Bew}
Let $\varphi\in\Mod_g(\alpha)\setminus\{\id\}$ with $\varphi^n=\id\;$ for some $n\in\N$. Then we also have $M_\varphi^n=I_{2g}$. But since $M_\varphi=\begin{pmatrix}I_g & A\\O_g&I_g\end{pmatrix}$, we have $M_\varphi^n:=\begin{pmatrix}I_g & nA\\O_g&I_g\end{pmatrix}$. So $M_\varphi$ can only be the identity. Therefore, $\varphi$ acts trivially on the homology. But the subgroup of $\Mod_g$ whose elements act trivially on the homology (the so called ``Torelli group'') is torsion free, see \cite[6.12]{farbmarg}. Therefore, we have $\varphi=\id$ and $\Mod_g(\alpha)$ is torsion free.
\end{Bew}

\chapter{Origamis}
\label{kapori}
In this chapter we look for ways to calculate the groups $\Aff^+(X)\cap\Mod_g(\alpha)$ for origamis.
We first give a brief summary of known facts about origamis that we'll need later. These facts are taken from \cite{gabidiss}, which also contains the proofs.
\section{A short summary of known facts about origamis}
\begin{definition}
An \begriff{origami} is given by a finite number of euclidean unity squares that are glued according to the following rules:
\begin{enumerate}
\item[i)] Every right side of a square is glued via a translation to a left side;
\item[ii)] every upper side of a square is glued via a translation to a lower side;
\item[iii)] the resulting closed surface is connected.
\end{enumerate}
\end{definition}

\begin{bem}
\label{Orieinf}
\begin{enumerate}
\item The simplest example for an origami is a unit square, whose upper side is glued with the lower side and whose right side is glued with the left side via a translation. The resulting surface is a torus that we'll denote by $E$. With $\infty$ we denote the image of the four vertices of the square in $E$.
\item Let $O$ be an origami and let $X$ the surface defined by $O$. Then $X$ covers the torus $E$, where the covering map is given by sending the copies of the unit squares in $X$ to the torus $E$. This covering $p\colon X\rightarrow E$ is finite and can only be ramified over the point $\infty$.
\begin{figure}[h]
\begin{center}
\includegraphics[scale=0.5]{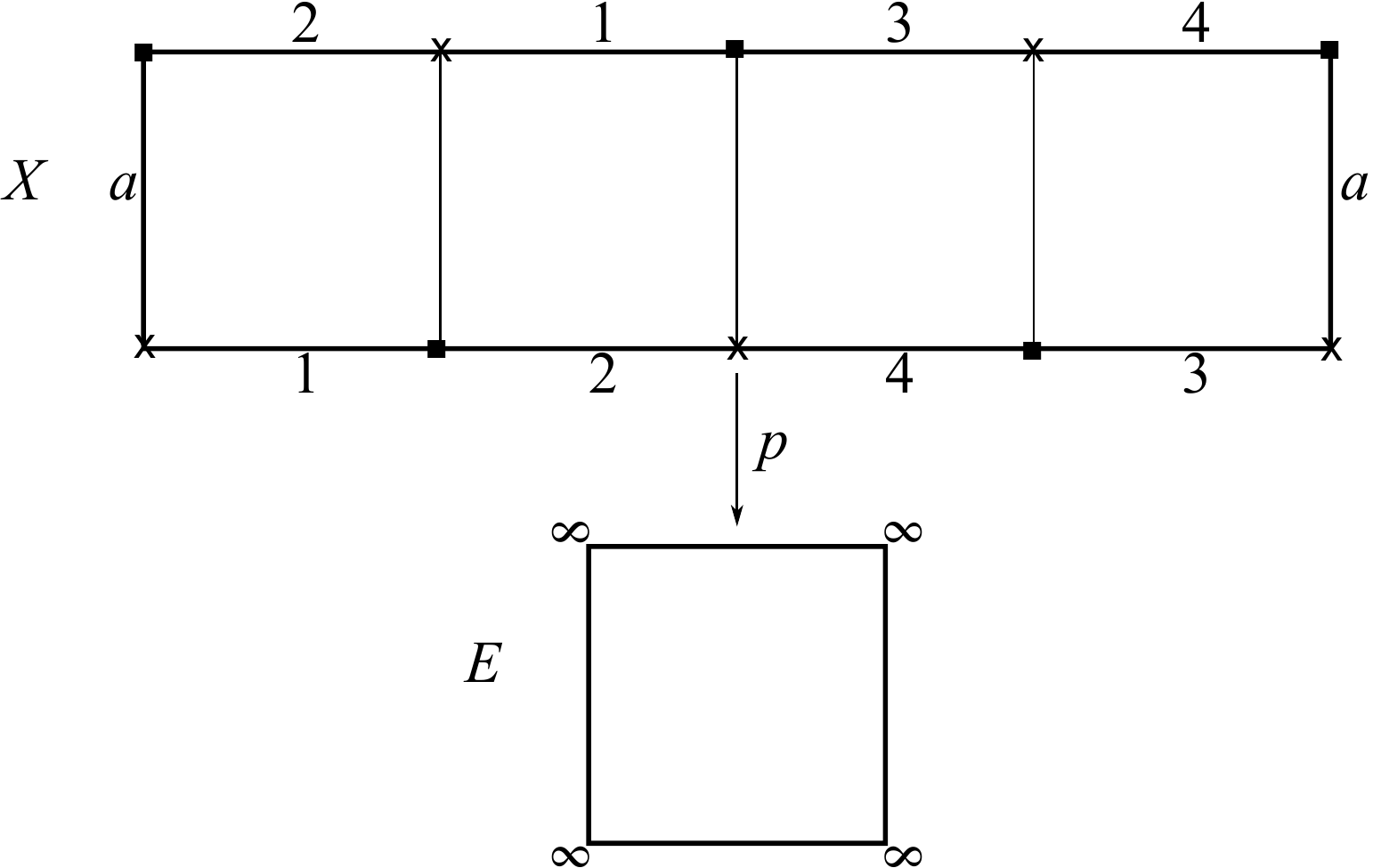}
\end{center}
\caption{Sides with the same labelling are identified}
\label{OrigamialsUberlagerung}
\end{figure}
\item An origami $O$ is uniquely determined by the unramified covering $p\colon X^*\rightarrow E^*$, where $E^*:=E\backslash\{\infty\}$ and  $X^*:=X\backslash p^{-1}(\infty)$. The points in $X\backslash X^*$ are called \begriff{vertices} of $O$.
\item The surface $X^*$, therefore, carries the structure of a translation surface that is induced by the corresponding structure on the torus. The translation structure defines a complex structure on $X$ and $X$ becomes a Riemann surface. From now on we'll take the translation structure on $E$ that comes from the lattice $\Lambda_I:=\{a+bi:a,b\in\mathds Z\}$, i.e.\ $E=\C/\Lambda_I$.
\item The fundamental group of $E^*$ is isomorphic to $F_2$. An isomorphism is given by sending the homotopy classes of the simple closed horizontal path and the simple closed vertical path that start and end in the center of the square to free generators $x$ and $y$. We fix this isomorphism $\pi_1(E)\cong F_2$ once and for all.
\begin{figure}[h]
\begin{center}
\includegraphics[scale=0.65]{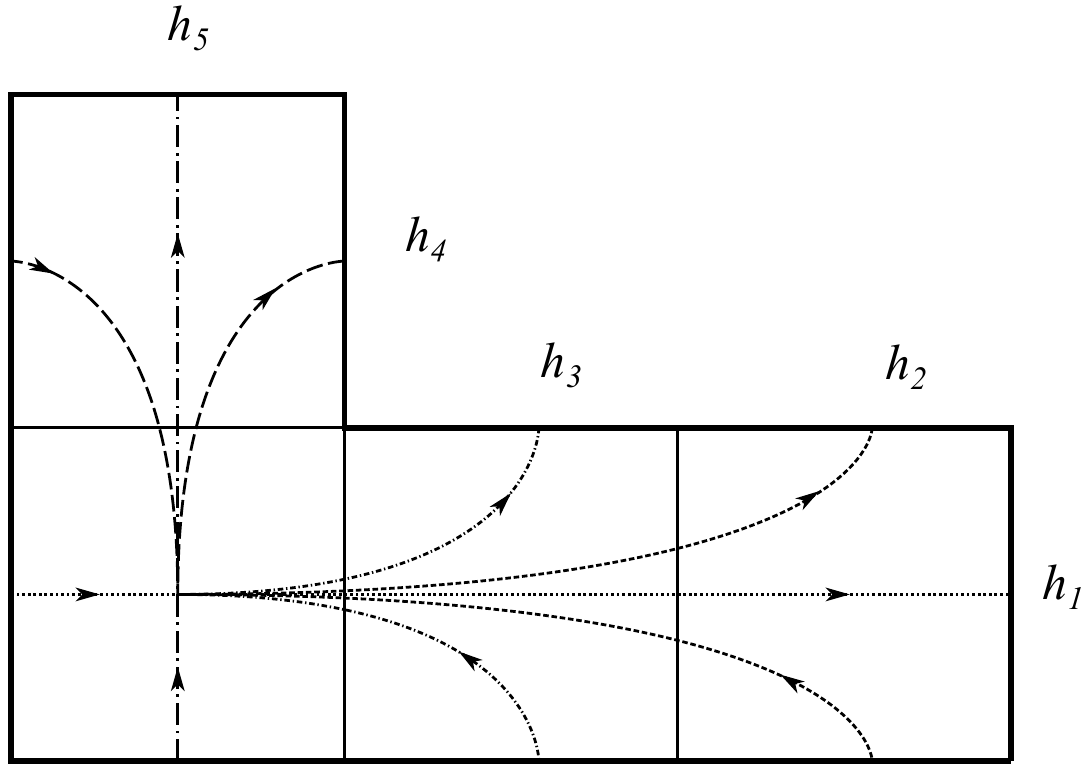}
$h_1=x^3,~h_2=x^2yx^{-2},~h_3=xyx^{-1},~h_4=yxy^{-1},~h_5=y^2$
\caption{Opposite sides are glued}
\end{center}
\label{L23}
\end{figure}
\item The fundamental group of $X^*$ is isomorphic to a subgroup of $F_2$ of finite index because $\pi_1(X^*)\cong\Deck(\HH/X^*)\leq\Deck(\HH/E^*)\cong F_2$.\\
One might also think of $\pi_1(X^*)$ as follows: Cut the punctured origami along sides of the squares until you get a simply connected surface and choose the center of one of the squares as base point $x_0$ of the fundamental group. For every side $i$ that has been cut, there is (up to inverse) a unique element $h_i\in F_2$ that belongs to a simple closed path in $X^*$ that starts in the base point and intersects the cut surface once in the side $i$. The subgroup $H$ of $F_2$ generated by the $h_i$'s is isomorphic to $\pi_1(X^*,x_0)$. Note that $H$ depends on the choice of the base point, another choice leads to a conjugation in $F_2$.
We, thus, have a bijection between origamis and conjugation classes of subgroups of $F_2$ with finite index.

\end{enumerate}
\end{bem}
For most theorems about origamis, we need the following theorem about automorphisms of free groups.
\begin{satz}
\label{hatbeta}
\begin{enumerate}
\item
The map
\[\begin{array}{lcl}
\hat{\beta}\colon &\Aut F_2 &\longrightarrow \GL_2(\mathds{Z})\\
&\varphi& \longmapsto
\begin{pmatrix}
\sharp_x(\varphi(x))&\sharp_x(\varphi(y))\\
\sharp_y(\varphi(x))&\sharp_y(\varphi(y))
\end{pmatrix}
\end{array}\]
is a surjective homomorphism with kernel $\Inn(F_2)$, where $(x,y)$ denotes a fixed pair of free generators and $\sharp_x(\varphi(y))$ denotes the sum of the exponents of $x$ that appear in the word $\varphi(y)$ (in the same way we define $\sharp_x(\varphi(x))$, $\sharp_y(\varphi(x))$ and $\sharp_y(\varphi(y))$).
\item Let $\varphi\in\Aut F_2$, $A:=\hat{\beta}(\varphi)$ and
\[\begin{array}{lcl}
\proj\colon & F_2&\rightarrow \mathds Z^2\\
& W(x,y)&\mapsto (\sharp_x(W(x,y)),\sharp_y(W(x,y)).
\end{array}\]
Then the following diagram commutes
\[
\begin{xy}
\xymatrix{
F_2 \ar[r]^\varphi\ar[d]_\proj & F_2\ar[d]^\proj\\
\mathds Z^2\ar[r]^{z\mapsto A\cdot z}&\mathds Z^2
}
\end{xy}
\]
\end{enumerate}
\end{satz}

\begin{Bew}
\begin{enumerate}
\item \underline{Homomorphism}: Let $\phi$ and $\varphi\in\Aut F_2$. Then we have:
\begin{align*}
\hat{\beta}(\phi\circ\varphi)&=
\begin{pmatrix}
\sharp_x(\phi\circ\varphi(x))&\sharp_x(\phi\circ\varphi(y))\\
\sharp_y(\phi\circ\varphi(x))&\sharp_y(\phi\circ\varphi(y))
\end{pmatrix} \\
&=\begin{pmatrix}
\sharp_x(\phi(x))\sharp_x(\varphi(x))+\sharp_x(\phi(y))\sharp_y(\varphi(x))&\sharp_x(\phi(x))\sharp_x(\varphi(y))+\sharp_x(\phi(y))\sharp_y(\varphi(y))\\
\sharp_y(\phi(x))\sharp_x(\varphi(x))+\sharp_y(\phi(y))\sharp_y(\varphi(x))&\sharp_y(\phi(x))\sharp_x(\varphi(y))+\sharp_y(\phi(y))\sharp_y(\varphi(y))
\end{pmatrix}\\
&=\begin{pmatrix}
\sharp_x(\phi(x))&\sharp_x(\phi(y))\\
\sharp_y(\phi(x))&\sharp_y(\phi(y))
\end{pmatrix} \cdot
\begin{pmatrix}
\sharp_x(\varphi(x))&\sharp_x(\varphi(y))\\
\sharp_y(\varphi(x))&\sharp_y(\varphi(y))
\end{pmatrix} = \hat{\beta}(\phi)\cdot\hat{\beta}(\varphi)
\end{align*}
Since $\hat{\beta}(\varphi)$ is invertible for all $\varphi\in\Aut F_2$, we have $\Bild(\hat{\beta})\subseteq\GL_2(\mathds{Z})$.\\\\
\underline{Surjective}: 
$\GL_2(\mathds{Z})$ is generated by the images of the automorphisms
\[(x\mapsto y, y\mapsto x), (x\mapsto y^{-1}, y\mapsto xy)\text{ and }(x\mapsto y^{-1}, y\mapsto x),\]
\[\text{i.e }
\begin{pmatrix}
0&1\\
1&0
\end{pmatrix},
\begin{pmatrix}
0&1\\
-1&1
\end{pmatrix}\text{ and }
\begin{pmatrix}
0&1\\
-1&0
\end{pmatrix}.\]
\underline{$\Inn(F_2)= \Kern\hat{\beta}$}:
For ``$\supseteq$'', see \cite{Chang}. The inclusion ``$\subseteq$'' is obvious.
\item Obvious.
\end{enumerate}
\end{Bew}

\begin{satz}
The Veech group $\Gamma(O)$ of an origami is given by:
$$\hat{\beta}(\Stab_{\Aut^+F_2}(H))~~\text{ with }~~\Stab_{\Aut^+F_2}(H):=\{\gamma\in\Aut^+F_2\colon \gamma(H)=H\},$$
where $H$ is a subgroup of $F_2$ associated with $O$ as in Rem. \ref{Orieinf}(f).
\end{satz}
\begin{Bew}
See \cite[Thm 1]{gabidiss}.
\end{Bew}
Note that $\hat{\beta}(\Stab_{\Aut^+F_2}(H))$ depends only on the conjugation class of $H$: Let $H'=aHa^{-1}$ and $\gamma\in\Stab_{\Aut^+F_2}(H)$. Then $$(c_{a\gamma(a)^{-1}}\circ\gamma)(H')=a\gamma(a)^{-1}\gamma(H')\gamma(a)a^{-1}=a\gamma(a^{-1}H'a)a^{-1}=a\gamma(H)a^{-1}=aHa^{-1}=H',$$
where $c_{h}\in\Inn F_2$ denotes conjugation with $h\in F_2$.

\section{Horizontal Schottky cut systems (HSS)}

\begin{definition} Let $F_2$ be the free group with a fixed pair of generators $(x,y)$.
\begin{enumerate}
\item $w\in F_2$ is called \begriff{horizontal} if it is of the following form:
$$w=\prod_{i=0}^nx^{c_i}yx^{d_i}y^{-1}~~~\text{or}~~~w=\prod_{i=0}^nx^{c_i}y^{-1}x^{d_i}y, ~~~\text{ for }c_i,d_i\in\Z, i\in\{1,\dots,n\}.$$
\item A word that is conjugate to a horizontal word is called \begriff{conjugate horizontal}.
\end{enumerate}
\end{definition}

\begin{bem}
Let $\varphi\in\Aut(F_2)$ be defined by $\varphi(x)=x$ and $\varphi(y)=x^ky$ for $k\geq0$, and let $w$ be a conjugate horizontal word. Then $\varphi(w)$ is also conjugate horizontal.
\end{bem}

\begin{defbem}
Let $O = (p\colon X^*\rightarrow E^*)$ be an origami and $x_0\in X^*$. Since the fundamental group $\pi_1(X^*,x_0)$ can be identified with a subgroup of $F_2$, we can talk about horizontal or conjugate horizontal elements in $\pi_1(X^*,x_0)$. The latter does not depend on the choice of $x_0$ since for $x_0'\in X^*$, the fundamental groups $\pi_1(X^*,x_0)$ and $\pi_1(X^*,x_0')$ are conjugate to each other in $F_2$.\\
Elements in $\pi_1(X,x_0)$ are called \begriff{(conjugate) horizontal} if they have a (conjugate) horizontal lift in $\pi_1(X^*,x_0)$. The same applies to free homotopy classes of paths in $X^*$ and $X$. For free homotopy classes, ``horizontal'' means the same as ``(conjugate) horizontal''\\
The kernel of the map $\pi_1(X^*,x_0)\ra\pi_1(X,x_0)$ is generated by paths around the punctures, which are conjugate horizontal.
\end{defbem}

\begin{definition}
We call a Schottky cut system on an origami that consists only of horizontal free homotopy classes of paths a \begriff{horizontal Schottky cut system (HSS)}.
\end{definition}

\begin{satz}
\label{hsss}
In an origami there is always a HSS.
\end{satz}

\begin{Bew}
Let $O = (p\colon X^*\rightarrow E^*)$ be an origami of genus $g$ and let $Z_1,\dots,Z_m$ be the horizontal cylinders of $O$ and let $c_1,\dots,c_n$ be a non-separating set (i.e.\ such that $X\setminus(\bigcup_{i=1}^nc_i)$ is connected) of pairwise disjoint horizontal curves on $X$ . We show by induction on $n$ that there is such a set for all $n\leq g$.\\
For $n=0$, there is nothing to show. Let $n<g$ and let $c_1,\dots,c_n$ be such a set. There are two cases:\\
\underline{Case 1}: There is a horizontal curve $\gamma'$ around a cylinder such that $\gamma'\in X\setminus(\bigcup_{i=1}^nc_i)$ and such that $\gamma'$ does not separate $X\setminus(\bigcup_{i=1}^nc_i)$.
Then $c_{n+1}:=\gamma'$ satisfies our requirements.\\
\underline{Case 2}: Every path around a cylinder which is not in $\{c_1,\dots,c_n\}$ separates $X\setminus(\bigcup_{i=1}^nc_i)$. Let $\gamma$ be some curve around a cylinder. We construct a finite sequence of simply closed curves $(\gamma_j)$ such that $\gamma_j$ does not separate $X\setminus(\bigcup_{i=1}^nc_i)$ for all $j$ and such that $\vert\gamma_{j+1}\cap\gamma\vert<\vert\gamma_j\cap\gamma\vert$.\\
$X\setminus(\bigcup_{i=1}^nc_i)$ is homeomorphic so the surface $S_{g-n}^{2n}$. Since $g-n\geq1$, there is a simple closed, non-separating curve $\gamma_1$ in $X\setminus(\bigcup_{i=1}^nc_i)$, i.e.\ $X\setminus(\gamma_1\cup\bigcup_{i=1}^nc_i)$ is connected. If $\gamma_1$ is horizontal, we are finished. Otherwise, a word representation of $\gamma_1$ contains a word of the form $yx^ay$ or $y^{-1}x^ay^{-1}$ ($a\in\Z$). That means that there is a cylinder $Z$ such that $\gamma_1$ enters $Z$ on the lower side and exits on the upper side (or the other way round). Let $\gamma$ be the curve that goes horizontally around $Z$. Suppose $\vert\gamma_1\cap\gamma\vert=1$. Then $\gamma_1\setminus\gamma$ is a curve from one side of $\gamma$ to the other and that means that $\gamma$ doesn't separate $X\setminus\bigcup_{i=1}^nc_i$ which contradicts our assumptions. So we have $\vert\gamma_1\cap\gamma\vert\geq2$.
\begin{figure}[h]
\begin{center}
\includegraphics[scale=0.3]{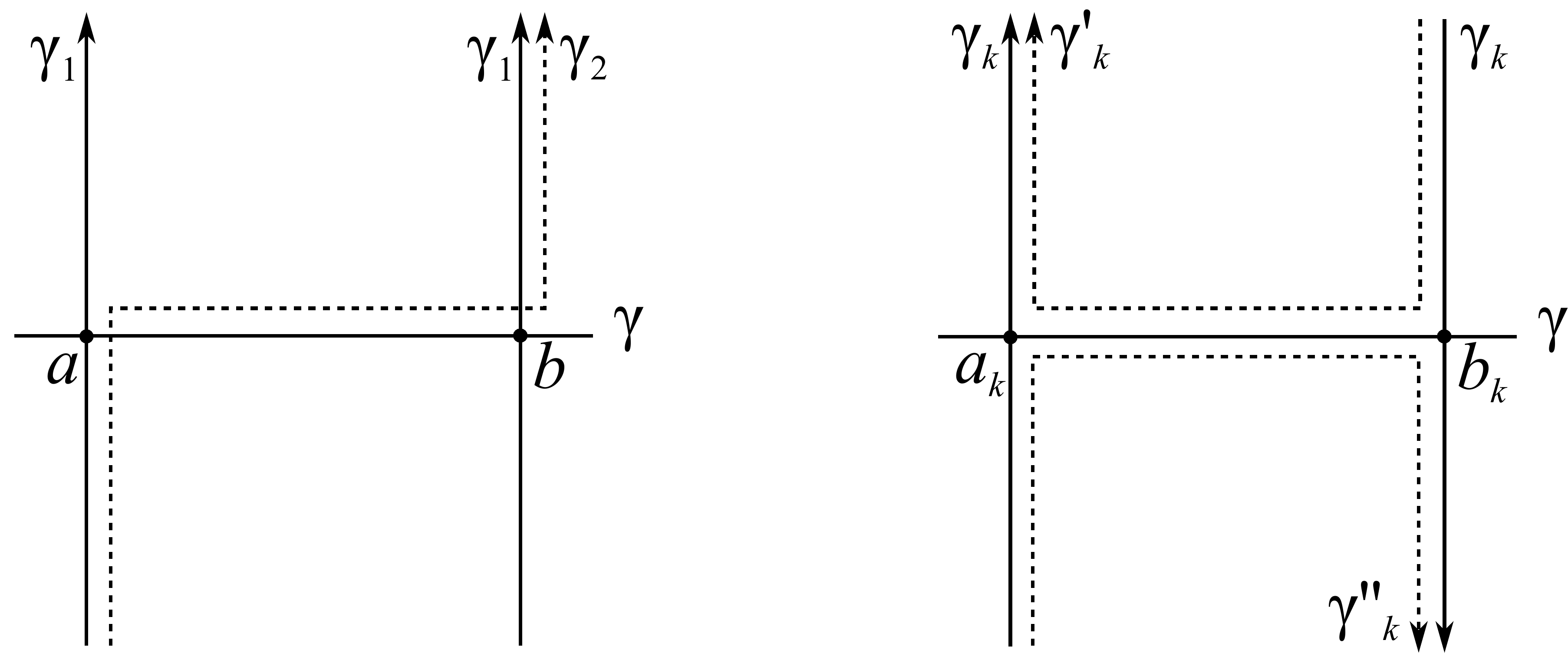}
\end{center}
\end{figure}\\
\underline{Step 1}: Let $a$ and $b$ be two intersection points of $\gamma_1$ and $\gamma$ that are consecutive with respect to $\gamma$ and with the same orientation (if such points do not exist, proceed with step 2). Let $\gamma_2$ be the curve that we obtain by substitution of the piece of $\gamma_1$ between $a$ and $b$ by the the piece of $\gamma$ between $a$ and $b$ (see the left figure above). This curve does not separate $X\setminus(\bigcup_{i=1}^nc_i)$ because the piece of $\gamma_1$ that has been substituted is a curve from one side of $\gamma_2$ to the other. Repeat step 1 until we obtain a curve $\gamma_k$ such that there are no more such consecutive points with the same orientation. This happens after a finite number of steps since $\vert\gamma_{j+1}\cap\gamma\vert\leq\vert\gamma_j\cap\gamma\vert -1$.\\
\underline{Step 2}: If $\vert\gamma_k\cap\gamma\vert\neq0$, there are consecutive (with respect to $\gamma$) intersection points $a_k$ and $b_k$ with opposite orientation. We obtain two curves $\gamma_k'$ and $\gamma_k''$ in the following way: Let $\gamma_k'$ be the curve that starts in $a_k$, follows $\gamma_k$ until $b_k$ and then follows $\gamma$ until $a_k$. Further, let $\gamma_k''$ be the curve that starts in $b_k$, follows $\gamma_k$ until $a_k$ and then follows $\gamma$ until $b_k$ such that $\gamma_k'$ and $\gamma_k''$ share the same part of the way with $\gamma$ (the part between $a_k$ and $b_k$, where there are no other points of $\vert\gamma\cap\gamma_k\vert$ apart from $a_k$ and $b_k$), see the right side of the figure above.\\
It cannot be that both $\gamma_k'$ and $\gamma_k''$ separate $X\setminus(\bigcup_{i=1}^nc_i)$ because then we could glue the surface with boundary $\gamma_k'$ to the surface with boundary $\gamma_k''$ on their common piece on $\gamma$. This would be a surface with $\gamma_k$ as boundary, and so $\gamma_k$ would be separating, a contradiction. Define $\gamma_{k+1}$ as one of the curves $\gamma_k'$ and $\gamma_k''$ that doesn't separate $X\setminus(\bigcup_{i=1}^nc_i)$.\\
We repeat step 1 and 2, until we get some curve $\gamma_r$ that is disjoint to $\gamma$. Then we apply these steps to the intersection points with the horizontal ways around the other horizontal cylinders of the origami and at the end we obtain a horizontal simple closed curve $\gamma_s$ that doesn't separate $X\setminus(\bigcup_{i=1}^nc_i)$. Let $c_{n+1}:=\gamma_s$.
\end{Bew}

\begin{bem}
In an origami there is always a HSS with at least one curve around a horizontal cylinder, i.e.\ a curve of the form $x^m$ (after proper choice of the base point), where $m$ denotes the length of that cylinder.
\end{bem}

\begin{Bew}
The horizontal cut $c_1$ around an arbitrary cylinder $Z$ is not separating. This follows from the fact that a path that starts in a square in $Z$ (in the point $x_0$) and goes in $y$-direction, comes back to $x_0$ after finitely many $y$-steps (there is a $t\in\N$ with $y^t\in\pi_1(X^*,x_0)$ because the fundamental group $\pi_1(X^*,x_0)$ has finite index in $F_2$). Now, at the latest, the path enters $Z$ from the lower side. That means that the upper and lower half of $Z$ lie in the same path component of $X\setminus c_1$. (This is equivalent to the fact that the intersection graph of a stable curve that is a limit point of an origamicurve has no bridges, see \cite[Lemma 2]{michimaier}. Also the other direction of the implication of the cited lemma is true, see \cite[Satz 2]{michimaier}.) So  $c_1$ can be taken as the first path for the construction of a HSS like in the proof of Thm. \ref{hsss}.
\end{Bew}

\begin{ko}
\label{ko}
The fundamental group of an origami always has a symplectic representation \linebreak$\langle a_1,b_1,\dots,a_g,b_g\mid\Pi_{i=1}^g[a_i,b_i]\rangle$ such that $a_1,\dots,a_g$ are conjugate horizontal. (A Schottky cut system can always be extended to a symplectic set of generators, as follows from the classification of surfaces, see \cite[§ 4.5]{kinsey}.)

\end{ko}
Every subgroup $H\leq F_2$ of finite index is isomorphic to the fundamental group of a punctured origami and we get the fundamental group of the corresponding non-punctured origami by factoring out the normal subgroup that is generated by the set $R$ of paths around the punctures. Therefore, the last result can be stated algebraically: 
\begin{ko}
Choose a base $(x,y)$ of the free group $F_2$ and let $H\leq F_2$ be a finite index subgroup. We define 
$$\overline{H}:=\FakRaum{H}{\langle\langle R\rangle\rangle},~~\text{ where } ~~ R:=\{w\in H\mid w=\tilde{w}(xyx^{-1}y^{-1})^n\tilde{w}^{-1}, \tilde{w}\in F_2, n\in\N\}.$$
Then there are $a_1,b_1,\dots,a_g,b_g\in\overline{H}$ such that the $a_1,\dots,a_g$ have conjugate horizontal representatives in $H$ and
\[\overline{H}=\langle a_1,b_1,\dots,a_g,b_g\mid\Pi_{i=1}^g[a_i,b_i]=1\rangle.\]
\end{ko}

\begin{satz}
\label{zylhorimkern}
Let $O = (p\colon X^*\rightarrow E^*)$ be an origami and let $(a_1,b_1,\dots,a_g,b_g)$ be a set of symplectic generators of $\pi_1(X,x_0)$ such that $a_1,..,a_g$ are conjugate horizontal. Further, let $\gamma:=\beta x^{\ell(Z)}\beta^-$ be a path, where
\begin{itemize}
\item $Z$ is a horizontal cylinder of length $\ell(Z)$ and
\item $\beta$ is a path that starts and ends in a point $x_0$ in $Z$.
\end{itemize}
Then $[\gamma]$ is in the normal subgroup generated by $a_1,\dots,a_g$.
\end{satz}

\begin{Bew}
We will call a homotopy that fixes the base point $x_0$ a $x_0$-homotopy (respectively we'll talk about $x_0$-homotopy classes) to distinguish them from free homotopies.\\
Let $Z$ be a cylinder and let $\gamma'$ be the horizontal path around $Z$, i.e.\ the path such that $x^{\ell(Z)}$ is the image of $\gamma'$ in $E^*$. Further, let $\alpha_1,\dots,\alpha_g$ be pairwise disjoint representatives of the free homotopy classes $a_1,\dots,a_g$. Then i$(\gamma',\alpha_i)=0$ for all $i\in\{1,\dots,g\}$, where i$(\cdot,\cdot)$ is the geometric intersection number between two curves. Therefore, we can choose $\alpha_1,\dots,\alpha_g$ that don't intersect $\gamma'$. Since $\widehat{X}:=\overline{X\setminus\bigcup_{i=1}^g\alpha_i}$ is connected, there is a path $\beta$ from the base point $x_0$ of the fundamental group to $\gamma'$ that does not intersect the curves $\alpha_i$. Let $\gamma:= \beta\gamma'\beta^-$. We have $[\gamma]\in\pi_1(\widehat{X},x_0)$. But $\pi_1(\widehat{X},x_0)$ is generated by the $x_0$-homotopy classes of the boundary curves around the $2g$ boundary components of $\widehat{X}$. In $X$, these boundary curves lie in homotopy classes conjugate to $a_1,\dots,a_g$.\\
\begin{figure}[h]
\begin{center}
\includegraphics[scale=0.55]{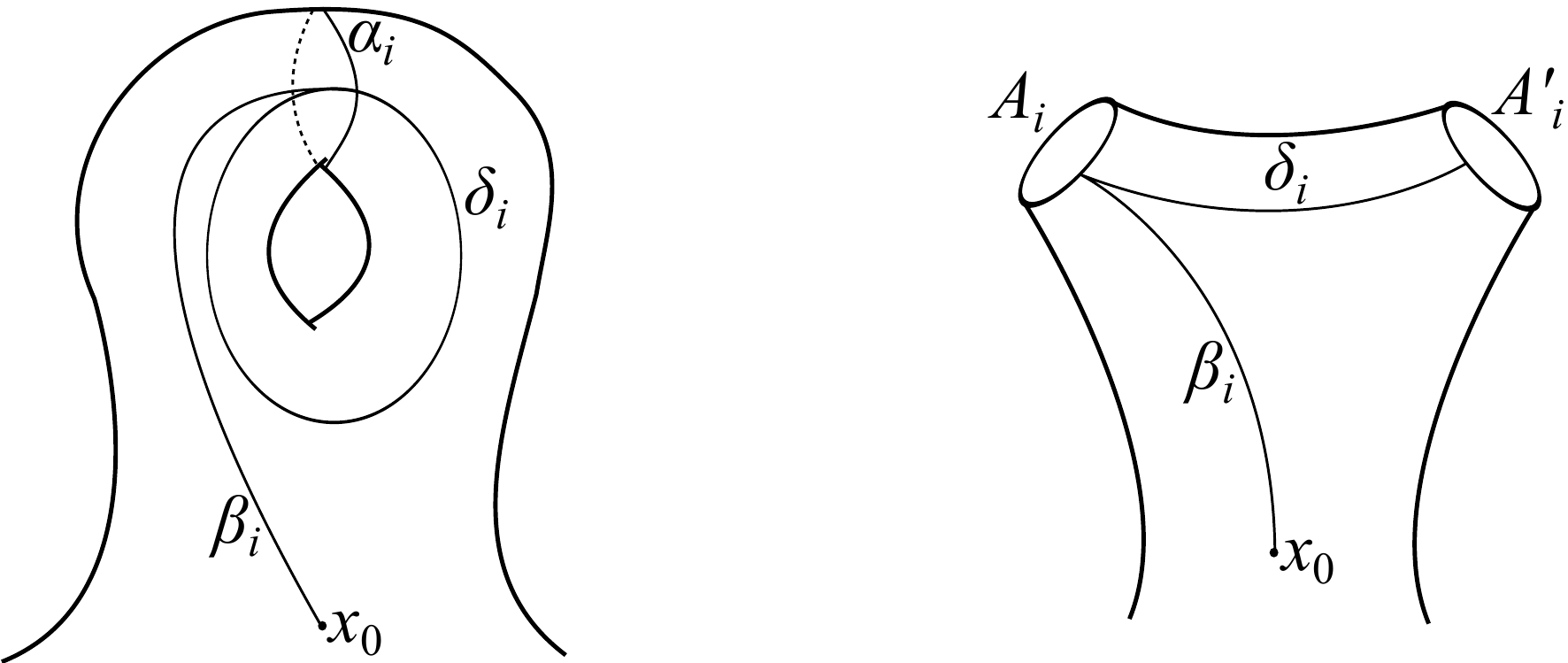}
\end{center}
\end{figure}\\
More precisely: Let $A_1,\dots,A_g,A_1',\dots,A_g'$ be the boundary curves of $\widehat{X}$, where for $i\in\{1,\dots,g\}$ $A_i,A_i'$ are the curves that we get when we cut $X$ along $\alpha_i$. Let $\beta_i$ be a path from $x_0$ to $\alpha_i$ such that $\alpha_i':=\beta_i\alpha_i\beta_i^-$ is a representative of the $x_0$-homotopy class $a_i$ (we assume w.l.o.g. that $\beta_1,\dots,\beta_g$ are also connected in $\widehat{X}$). W.l.o.g., we can choose $\beta_i$ for $i\in\{1,\dots,g\}$ to be a curve that goes from $x_0$ to $A_i$, after cutting the surface $X$, that is, the curve $\alpha_i'$ in $\widehat{X}$ lies in the $x_0$-homotopy class of a simple closed curve around $A_i$.\\
For each $i\in\{1,\dots,g\}$, we now choose a curve $\delta_i$ from $A_i$ to $A_i'$ that starts at the end point of $\beta_i$ on $A_i$ and ends at the corresponding point on $A_i'$; $\delta_i$ is a closed curve on $X$ and so is
$$\alpha_i'':=\beta_i\delta_i\alpha_i\delta_i^-\beta_i^-.$$
For their $x_0$-homotopy class in $X$, we have:
$$[\alpha_i'']=\underbrace{[\beta_i\delta_i\beta_i^-]}_{:=d_i\in\pi_1(X,x_0)}\underbrace{[\beta_i\alpha_i\beta_i^-]}_{=a_i}\underbrace{[\beta_i\delta_i^-\beta_i^-]}_{=d_i^{-1}\in\pi_1(X,x_0)},$$
i.e.\ $[\alpha_i'']$ is conjugate to $a_i$. In $\widehat{X}$ the curve $\alpha_i''$ goes around $A_i'$. This means that $\pi_1(\widehat{X},x_0)$ is generated by the $x_0$-homotopy classes of $\alpha_1',\dots,\alpha_g',\alpha_1'',\dots,\alpha_g''$. But these correspond to
\[a_1,\dots,a_g,d_1a_1d_1^{-1},\dots,d_ga_gd_g^{-1}\]
in $\pi_1(X,x_0)$. So we have $\gamma\in\langle\langle a_1,\dots,a_g\rangle\rangle$. The choice of $\beta$ is not a restriction because for every path $\hat{\gamma}:=\hat{\beta}x^{\ell(Z)}\hat{\beta}^-$ that satisfies the assumptions, we have that $[\hat{\gamma}]$ is conjugate to $[\gamma]$ in $\pi_1(X,x_0)$.
\end{Bew}
We can formulate also this theorem in algebraic terms:
\begin{ko}
Let $H\leq F_2, R$ and $\overline{H}=\langle a_1,b_1,\dots,a_g,b_g\mid\Pi_{i=1}^g[a_i,b_i]=1\rangle$ be defined as above. Then all $\overline{wx^nw^{-1}}$ 
that are in $\overline{H}$ are also in the normal subgroup generated by $a_1,\dots,a_g$.
\end{ko}

\section{The group $\Aff^+(X,\mu)\cap\Mod_g(\alpha)$ for a given $\alpha$}

\begin{satz}
\label{Oridiagramm}
Let $O = (p\colon X^*\rightarrow E^*)$ be an origami, $x_0\in X^*$ the base point of the fundamental group and $H\cong\pi_1(X^*,x_0)\cong\Deck(\HH/X^*)$ the corresponding subgroup of $F_2$. Moreover, let $f\in\Aff^+(X)$ and $\f\in\Aut^+(\pi_1(X,x_0))$ be a representative of the outer automorphism of $\pi_1(X,x_0)$ induced by $f$ and $i_*:\pi_1(X^*,x_0)\ra\pi_1(X,x_0)$ the homomorphism induced by the inclusion \hbox{$i\colon X^*\hookrightarrow X$.} Then there is a lift $\hat{f}_\star\colon F_2\ra F_2$ of $f$ such that the following diagram commutes for $A:= \hat{\beta}(\hat{f}_\star)$:
\[
\begin{xy}
\xymatrix{F_2\ar[ddd]_{\proj}\ar[rrr]^{\hat{f}_\star}&&& F_2\ar[ddd]^{\proj}\\
&\pi_1(X^*,x_0)\ar[ddl]_{\proj}\ar[r]^{\hat{f}_\star\mid_H}\ar[d]_{(1)}^{i_*}\ar@{_(->}[ul]&\pi_1(X^*,x_0)\ar[ddr]^{\proj}\ar[d]^{(2)}_{i_*}\ar@{^(->}[ur]&\\
&\pi_1(X,x_0)\ar[dl]\ar[r]^{\overline{f}}&\pi_1(X,x_0)\ar[dr]&\\
\Z^2 \ar[rrr]^{z\mapsto Az}&&&\Z^2
}
\end{xy}
\]
\end{satz}
\begin{Bew}
$F_2$ is isomorphic to the group $\Deck(\HH/E^*)$, where $E^*$ is the once-punctured torus. More precisely: 
$$F_2=\{h\in\Aff^+(\HH)\mid \der(h)=I\}.$$
We defined the homomorphism ``$\der$'' in §\ref{teichveech}. With this interpretation of $F_2$ we can define the following homomorphism, which is known to be an isomorphism, see \cite[Prop 3.5]{gabidiss}:
\[\begin{array}{rrl}
\star:&\Aff^+(\HH)&\ra\Aut^+(F_2)\\
&h&\mapsto(h_\star:\sigma\mapsto h\circ\sigma\circ h^{-1})
\end{array}\]
Now, every $f\in\Aff^+(X^*)$ can be lifted to some $\hat{f}\in\Aff^+(\HH)$, see \cite[Prop 3.3]{gabidiss}.\\
We have $\der(f)=\der(\hat{f})$ and by the proof of \cite[Prop 3.5]{gabidiss} we have $\der({\hat{f}})=\hat{\beta}(\hat{f}_\star)=A$. Therefore, the outer diagram commutes by Thm. \ref{hatbeta}b.\\
In order to prove that the diagram in the middle commutes we must show that the restriction of $\hat{f}_{\star}:F_2\ra F_2$ to the subgroup $H\cong\pi_1(X^*,x_0)$ is induced by the restriction $f\colon X^*\ra X^*$. For this purpose, we first recall the construction of the automorphism of $\pi_1(X^*,x_0)$ induced by $f$ and then remind the reader of the canonical isomorphism between $\pi_1(X^*,x_0)$ and the deck transformation group. With these prerequisites at our disposal, we then can compare $\hat{f}_{\star}$ with the automorphism induced by $f$, simply by interpreting $\hat{f}_{\star}|_H$ as automorphism of $\Deck(\mathds H/X^*)$ since $H\cong\Deck(\mathds H/X^*)$.\\
Let $p:\mathds H\ra X^*$ be the covering map, and choose a preimage $\tilde{x}_0\in\mathds H$ of $x_0$ and a path $\tilde{s}\colon [0,1]\ra\mathds H$ that starts in $\tilde{s}(0)=\tilde{x}_0$ and ends in $\tilde{s}(1)=\hat{f}(\tilde{x}_0)$. Then $s:=p\circ\tilde{s}$ is a path from $x_0$ to $f(x_0)$. The map $f$ induces an $f_*\in\Aut^+(\pi_1(X^*,x_0))$ as follows: Let $a:[0,1]\ra X^*$ be a closed path in $X^*$ with $a(0)=a(1)=x_0$ and let $[a]\in\pi_1(X^*,x_0)$ be the homotopy class of $a$. Then define $f_*([a]):=[s\cdot(f\circ a)\cdot s^-]$. (Note that even if $f_*$ is usually only well-defined up to conjugation, here it is well-defined since $s$ is defined up to homotopy by $\tilde{x}_0$ and the lift $\hat{f}$.)\\
Moreover, we identify $\pi_1(X^*,x_0)$ with $\Deck(\mathds H/X^*)$ as follows: Take $a$ and \hbox{$[a]\in\pi_1(X^*,x_0)$} as above; let $\tilde{a}$ be the lift of $a$ that starts in $\tilde{x}_0$ and let $h_a$ be the deck transformation that maps $\tilde{x}_0$ to the end point of $\tilde{a}$. Then we associate $[a]\in\pi_1(X^*,x_0)$ to the deck transformation $h_a$. This is well defined because lifts of homotopic paths are homotopic and $h_a$ exist since the covering $\mathds H/X^*$ is galois.\\
Now take $a$ as above. We view $\hat{f}_{\star}|_H$ as an automorphism of $\Deck(\mathds H/X^*)$. We have to show that $\hat{f}_{\star}(h_a)=h_{s\cdot(f\circ a)\cdot s^-}$.\\
To this end, take $\tilde{a}$ as above and let $\tilde{y}_0$ be the end point of $\tilde{a}$.\ The deck transformation $\hat{f}_{\star}(h_a)=\hat{f}\circ h_a\circ\hat{f}^{-1}$ is the one that maps $\hat{f}(\tilde{x}_0)$ to $\hat{f}(\tilde{y}_0)$.
To understand the deck transformation $h_{s\cdot(f\circ a)\cdot s^-}$, we look at the lift of the path $s\cdot(f\circ a)\cdot s^-$ that starts in $\tilde{x}_0$. The first part of it is the lift $\tilde{s}$ of $s$, the second part is the lift $\hat{f}\circ \tilde{a}$ of $f\circ a$ that starts in $\hat{f}(\tilde{x}_0)$ and ends in $\hat{f}(\tilde{y}_0)$, the third part is the lift $\hat{s}^-$ of $s^-$ that starts in $\hat{f}(\tilde{y}_0)$. Let $\tilde{z}_0$ be its end point. Then $h_{s\cdot(f\circ a)\cdot s^-}$ is the deck transformation that maps $\tilde{x}_0$ to $\tilde{z}_0$, and this is exactly the deck transformation that maps $\tilde{s}$ to $\hat{s}$ since it is the deck transformation that maps the starting point $\tilde{x}_0$ of $\tilde{s}$ to the starting point $\tilde{z}_0$ of $\hat{s}$ and, therefore, also maps the end point of $\tilde{s}$ to the end point of $\hat{s}$, i.e.\ maps $\hat{f}(\tilde{x}_0)$ to $\hat{f}(\tilde{y}_0)$. Therefore, it is the deck transformation $\hat{f}_{\star}(h_a)$. That means, the diagram in the middle commutes.\\
The diagrams (1) and (2): Let $h_1,\dots,h_d$ be free generators of $H\cong\pi_1(X^*,x_0)$ and $r_1,\dots,r_l\in H$ the homotopy classes of simple closed paths (with a fixed orientation) around the punctures. Then we have
\[\pi_1(X,x_0)=\langle h_1,\dots,h_d\mid r_1,\dots,r_l\rangle.\]
In $F_2$ the relation $r_i$ is conjugate to $(xyx^{-1}y^{-1})^{v_i}$, where $v_i$ denotes the ramification index of $p$ in the corresponding puncture for $i\in\{1,\dots,l\}$. 
So we have:
\[\proj(r_i)=(0,0)^T\quad i\in\{1,\dots,l\}\]
and the homomorphism $\proj\colon \pi_1(X^*,x_0)\ra\mathds Z^2$ factorises over $\pi_1(X,x_0)$ and the diagrams (1) and (2) commute. 
We also denote the corresponding (well defined) homomorphism $\pi_1(X,x_0)\ra\Z^2$ by $\proj$.\\
Since the homomorphism $\pi_1(X^*,x_0)\ra\pi_1(X,x_0)$ is surjective, also the lower quadrangle commutes. For the rest, there's nothing to show.
\end{Bew}

\begin{satz}
\label{Formelprojkern}
Let $O=(p:X^*\rightarrow E^*)$ be an origami and $(a_1,b_1,\dots,a_g,b_g)$ a set of symplectic generators of $\pi_1(X,x_0)$, $\alpha$ a symplectic homomorphism
with $\alpha(a_i)=1$ and $\alpha(b_i)=\gamma_i$ for $i\in\{1,\dots,g\}$ and let $f\in\Aff^+(X)\cap\Mod_g(\alpha)$. Then we have:
\begin{enumerate}
\item $\der(f)\neq\pm I$. That means, $f$ is neither a translation nor an involution,
\item $\der(f)\cdot\proj(a_i)=\proj(a_i),~~i\in\{1,\dots,g\}$,
\item $\der(f)\cdot\proj(b_i)=\proj(b_i)+\sum_{j=1}^gh_{ij}\cdot\proj(a_j),~~i\in\{1,\dots,g\},~h_{ij}\in\mathds Z$.
\end{enumerate}
\end{satz}

\begin{Bew}
\begin{enumerate}
\item A translation of an origami corresponds to a permutation of the squares that respects the glueings. In particular, translations have finite order. The same holds for involutions since they have the form $z\mapsto-z+b$ and thier square is, thus, a translation. Since $\Mod_g(\alpha)$ is torsion free by Thm. \ref{torsionsfrei}, elements in $\Aff^+(X)\cap\Mod_g(\alpha)$ cannot have finite order.
\item Because the lower quadrangle in Thm. \ref{Oridiagramm} commutes, we have
$$\der(f)\cdot\proj(a_i)=\proj(f_*(a_i)).$$
Because of Lemma \ref{EW} b) and c), this equals $\proj(a_i)$.
\item $\Kern\alpha$ is the normal subgroup generated by $a_1,\dots,a_g$. Therefore, $\proj(\Kern\alpha)$ is generated by $\proj(a_1),\dots,\proj(a_g)$. Since $f\in\Mod_g(\alpha)$, we have $f_*(\Kern\alpha)=\Kern\alpha$ and $f_*(b_i)\equiv b_i\operatorname{mod}\Kern\alpha$. Again it follows from the fact that the lower quadrangle in Thm. \ref{Oridiagramm} commutes that $\der(f)\cdot\proj(b_i)=\proj(f_*(b_i))$. The assertion now follows from
\[\der(f)\cdot\proj(b_i)\equiv\proj(b_i)\operatorname{mod}\proj(\Kern\alpha).\]
\end{enumerate}
\end{Bew}

\begin{lemma}
\label{EW}
For a field $\K$ and $n\geq 2$, let $\phi:\K^2\ra\K^2$ and $\hat{\phi}:\K^n\ra\K^n$ be bijective linear maps and $p:\K^n\ra\K^2$ a surjective linear map. Moreover, assume that $\phi$ is diagonalizable with eigenvalues $\lambda_1$ and $\lambda_2$ and that the diagram
\[
\begin{xy}
\xymatrix{\K^n\ar[r]^{\hat{\phi}}\ar[d]_p&\K^n\ar[d]^p\\
\K^2\ar[r]^{\phi}&\K^2\\
}
\end{xy}
\]
commutes. Then $\lambda_1$ and $\lambda_2$ are also eigenvalues of $\hat{\phi}$.
\end{lemma}

\begin{Bew}
Let $v_1$ and $v_2$ be linear independent eigenvectors for the eigenvalues $\lambda_1$ and $\lambda_2$ and take the basis $B:=\{v_1,v_2\}$ of $\K^2$ and let $w_1,w_2\in\K^n$ be preimages of $v_1$ and $v_2$ under $p$. Moreover, let $\{w_3,\dots,w_n\}$ be a basis of $\Kern(p)$. This can be extended to a basis $B'$ of $\K^n$ with $w_1$ and $w_2$ because
\[
\sum_{i=1}^na_iw_i=0\Rightarrow p(\sum_{i=1}^na_iw_i)=0\Rightarrow a_1p(w_1)+a_2p(w_2)=0\Rightarrow a_1v_1+a_2v_2=0\Rightarrow a_1=a_2=0
\]
But since $\{w_3,\dots,w_n\}$ is a basis of $\Kern(p)$, it also follows $a_3=\dots=a_n=0$.\\
Since $\hat{\phi}$ is bijective and the diagram commutes, we have $\hat{\phi}(\Kern(p))=\Kern(p)$. The commutativity of the diagram implies $\hat{\phi}(w_i)\in\lambda_iw_i+\Kern(p)$ for $i\in\{1,2\}$.
Then the transformation matrix of $\hat{\phi}$ with respect to the basis $B'$ has the form
\[
\begin{pmatrix}
\lambda_1 & 0 & 0 & \cdots & 0\\
0 & \lambda_2 & 0 & \cdots & 0\\
* & & \cdots & & *\\
\vdots & & & & \vdots\\
* & & \cdots & & *\\
\end{pmatrix}
\]
and this matrix has eigenvalues $\lambda_1$ and $\lambda_2$.
\end{Bew}

\begin{satz}
\label{nichthyp}
Let $O = (p\colon X^*\rightarrow E^*)$ be an origami and let $f\in\Aff^+(X)$ such that $A:=\der(f)$ is hyperbolic. Further, let $M_f\in\Z^{2g\times2g}$ be the transformation matrix for the map on $H_1(X,\Z)$ induced by $f$ (see definition \ref{mphi}). Then the two eigenvalues of $A$ are also eigenvalues of $M_f$, viewed as matrix in $\R^{2g\times2g}$.
\end{satz}

\begin{Bew}
Denote by $\pi_1^{ab}(X,x_0)$ the abelianised fundamental group and let $\overline{f}$ be the automorphism of $\pi_1(X,x_0)$ induced by $f$ (which is unique up to an inner automorphism) as in Thm. \ref{Oridiagramm}. We have the following diagram:
\[
\begin{xy}
\xymatrix{\pi_1(X,x_0)\ar[dd]_{\proj}\ar[rrrr]^{\overline{f}}_{\stackrel{\stackrel{~}{\stackrel{~}{\stackrel{~}{~}}}}{\text{ (I) }}}\ar[dr]^h&&&&\pi_1(X,x_0)\ar[dl]_h\ar[dd]^{\proj}&\\
&\pi_1^{ab}(X,x_0)\ar[rr]^{x\mapsto M_fx}_{\stackrel{\stackrel{~}{\stackrel{~}{\stackrel{~}{~}}}}{\text{ (II) }}}\ar[dl]&&\pi_1^{ab}(X,x_0)\ar[dr]&\\
\Z^2 \ar[rrrr]^{z\mapsto Az}&&&&\Z^2
}
\end{xy}
\] 
Let $h$ be the projection of $\pi_1(X,x_0)$ to its abelianisation. Note that the kernel of $h$ is contained in the kernel of $\proj$. The homomorphism $\proj$ hence factorizes over $\pi_1^{ab}(X,x_0)$ and we have a homomorphism from $\pi_1^{ab}(X,x_0)$ to $\Z^2$. The two triangular parts of the diagram, therefore, commute.\\
We now show that part (II) commutes: The outer quadrangle commutes by Thm. \ref{Oridiagramm} and part (I) commutes by definition of $M_f$. Since the map $h$ is surjective, part (II) commutes, too.\\
We have $\pi_1^{ab}(X,x_0)\cong\Z^{2g}$ and the group homomorphisms defined by $M_f$ and $A$ can be extended uniquely to linear maps form $\R^{2g}$ to $\R^{2g}$ and from $\R^2$ to $\R^2$, respectively. We are, therefore, in the situation of Lemma \ref{EW} and the assertion follows.
\end{Bew}

\begin{satz}
\label{isomzuz}
Let $O=(p\colon X^*\rightarrow E^*)$ be an origami of genus $g$ and $\alpha:\;\pi_1(X,x_0)\ra F_g$ a symplectic homomorphism. Then the group $\Gamma:=\der(\Aff^+(X)\cap\Mod_g(\alpha))$ has no hyperbolic elements and $\Aff^+(X)\cap\Mod_g(\alpha)$ consists of the identity only or is isomorphic to $\Z$.
\end{satz}

\begin{Bew}
Let $f\in\Aff^+(X)\cap\Mod_g(\alpha)$ and define $A:=\der(f)$. According to Lemma \ref{EW1}, the matrix $M_f$ has eigenvalue 1 with multiplicity $2g$. Theorem \ref{nichthyp} implies that $A$ cannot be hyperbolic, and since $\Mod_g(\alpha)$ is torsion free, $A$ cannot be elliptic. Therefore, $\Gamma$ can only have parabolic elements and the identity.\\
We show that the group $\overline{\Gamma}:=\Gamma/\{\pm I\}$ is abelian. We can view $\overline{\Gamma}\leq\PSL_2(\Z)$ as group of Möbius transformations in $\Isom^+(\mathds H)$. Let $B,C\in\overline{\Gamma}$. If $B$ and $C$ have the same fixed point (w.l.o.g. let $\infty$ be this fixed point), then $B$ and $C$ are of the form \[B=\begin{bmatrix}1&b\\0&1\end{bmatrix}\quad\text{ and }\quad C=\begin{bmatrix}1&c\\0&1\end{bmatrix},\quad b,c\in\R\]
and commute. If they don't have the same fixed point, let w.l.o.g. be $\infty$ the fixed point of $B$ and $0$ the fixed point of $C$. $B$ and $C$ are then of the form
\[B=\begin{bmatrix}1&b\\0&1\end{bmatrix}\quad\text{ and }\quad C=\begin{bmatrix}1&0\\c&1\end{bmatrix},\quad b,c\in\R,\]
where we assume w.l.o.g. that $b,c>0$. But
\[BCB^{-1}C^{-1}=\begin{bmatrix}1&b\\0&1\end{bmatrix}\cdot\begin{bmatrix}1&0\\c&1\end{bmatrix}\cdot\begin{bmatrix}1&-b\\0&1\end{bmatrix}\cdot\begin{bmatrix}1&0\\-c&1\end{bmatrix}=\begin{bmatrix}(1+bc)^2-bc&-b^2c\\bc^2&1-bc\end{bmatrix}\]
has trace $(1+bc)^2-bc+1-bc=2+(bc)^2>2$ and is, therefore, hyperbolic, a contradiction.\\
So $\overline{\Gamma}$ is abelian and, thus, isomorphic to $\Z$ or trivial because a Fuchsian group is abelian if and only if it is cyclic (see \cite[Thm. 2.3.6]{katok}).\\
$\Aff^+(X)\cap\Mod_g(\alpha)$ is also isomorphic to $\Z$ or trivial since the group homomorphism
\[/\{\pm I\}\circ\der\mid_{\Aff^+(X)\cap\Mod_g(\alpha)}:\;\Aff^+(X)\cap\Mod_g(\alpha)\ra\overline{\Gamma}\]
is injective. This follows from the fact that the kernel of the above map can contain only the identity, translations and involutions. But there are no translations and involutions in $\Aff^+(X)\cap\Mod_g(\alpha)$, by Thm. \ref{Formelprojkern}a, so the kernel is trivial.
\end{Bew}

\begin{bem}
The statement of Thm. \ref{nichthyp} is a well known fact for translation surfaces. For this more general case, the proof is more complicated and will be done later, see Thm. \ref{nichthypallgemein}. But this new proof for the special case of origamis is elementary. That's why we included it.
\end{bem}

\begin{satz}
\label{schnittnichtleer}
Let $O=(p:X^*\rightarrow E^*)$ be an origami and $(a_1,b_1,\dots,a_g,b_g)$ a set of symplectic generators of $\pi_1(X,x_0)$ such that $a_1,\dots,a_g$ are conjugate horizontal, and let
$$\alpha:\pi_1(X,x_0)\ra F_g=\langle\gamma_1,\dots,\gamma_g\rangle$$
be the symplectic homomorphism induced by $\alpha(a_i)\mapsto1$ and $\alpha(b_i)\mapsto\gamma_i$. Then we have:
\[\Aff^+(X)\cap\Mod_g(\alpha)\cong\Z.\]
\end{satz}

\begin{Bew}
Let $Z_1,\dots,Z_r$ be the horizontal cylinder of $O$ and let
$$m:=\kgV(\ell(Z_1),\dots,\ell(Z_r)).$$
Then we have $A:=\begin{pmatrix}1&m\\0&1\end{pmatrix}\in\Gamma(O)$ and there exists an $f\in\Aff^+(X)$ with $\der(f)=A$.\\
Denote by $\hat{f}_\star$ the automorphism given by $\hat{f}_\star(x)=x$ and $\hat{f}_\star(y)=x^my$. This map $\hat{f}_\star$ is a lift of $f$ to $F_2$ as in Thm. \ref{Oridiagramm}.
Let $w\in\pi_1(X^*,x_0)$ be a reduced word in $x$ and $y$. We denote by $\overline{w}\in\pi_1(X,x_0)$ the equivalence class of $w$ and assume w.l.o.g. that $w$ contains at least one $y$ and is of the form $w=w_1y^{\varepsilon}w_2$ with $\varepsilon\in\{-1,1\}$.\\
We show $\alpha(\overline{\hat{f}_{\star}(w)})=\alpha(\overline{w_1y^{\varepsilon}\hat{f}_{\star}(w_2)})$ by induction on the (absolute) number $n$ of $y$ that appear in the part $w_1$ of $w$ (that means, we count $y$ and $y^{-1}$ both as occurences of $y$ in $w_1$).\\
In the case $n=0$, $w_1$ is of the form $x^k$. We first look at the case $\varepsilon=1$:
\begin{align*}
\hat{f}_{\star}(w)&=\hat{f}_{\star}(x^k)\hat{f}_{\star}(y)\hat{f}_{\star}(w_2)\\
&=x^kx^my\hat{f}_{\star}(w_2)\\
&=x^mx^ky\hat{f}_{\star}(w_2)
\end{align*}
By Thm. \ref{zylhorimkern}, we have $\overline{x^m}\in\Kern(\alpha)$. Therefore,
\[\alpha(\overline{\hat{f}_{\star}(w)})=\alpha(\overline{x^m}\cdot\overline{x^ky\hat{f}_{\star}(w_2)})=\alpha(\overline{x^ky\hat{f}_{\star}(w_2)}).\]
Now let $\varepsilon=-1$:
\begin{align*}
\hat{f}_{\star}(w)&=\hat{f}_{\star}(x^k)\hat{f}_{\star}(y^{-1})\hat{f}_{\star}(w_2)\\
&=x^ky^{-1}x^{-m}\hat{f}_{\star}(w_2)\\
&=x^ky^{-1}x^{-m}yx^{-k}x^ky^{-1}\hat{f}_{\star}(w_2)
\end{align*}
Here, too, Thm. \ref{zylhorimkern} implies $\overline{x^ky^{-1}x^{-m}yx^{-k}}\in\Kern(\alpha)$ and, therefore, $\alpha(\overline{\hat{f}_{\star}(w)})=\alpha(\overline{x^ky^{-1}\hat{f}_{\star}(w_2)})$.\\
Let us now consider the case where the number of $y$ appearing in $w_1$ is $n+1$:\\
Then we can write $w_1$ as $w_1=\tilde{w}_1y^{\delta}x^k$ with $\delta\in\{-1,1\}$ and we have 
\[\alpha(\overline{\hat{f}_{\star}(w)})=\alpha(\overline{\tilde{w}_1y^{\delta}\hat{f}_{\star}(x^ky^{\varepsilon}w_2)}),\]
by induction.\\
Again, we first look at the case $\varepsilon=1$:
\begin{align*}
\alpha(\overline{\hat{f}_{\star}(w)})&=\alpha(\overline{\tilde{w}_1y^{\delta}\hat{f}_{\star}(x^kyw_2)})\\
&=\alpha(\overline{\tilde{w}_1y^{\delta}x^kx^my\hat{f}_{\star}(w_2)})\\
&=\alpha(\overline{(\tilde{w}_1y^{\delta}x^k)x^m(\tilde{w}_1y^{\delta}x^k)^{-1}(\tilde{w}_1y^{\delta}x^k)y\hat{f}_{\star}(w_2)})\\
&=\alpha(\underbrace{\overline{(\tilde{w}_1y^{\delta}x^k)x^m(\tilde{w}_1y^{\delta}x^k)^{-1}}}_{\in\Kern(\alpha)}\cdot\overline{\tilde{w}_1y^{\delta}x^ky\hat{f}_{\star}(w_2)})\\
&=\alpha(\overline{\tilde{w}_1y^{\delta}x^ky\hat{f}_{\star}(w_2)})\\
&=\alpha(\overline{w_1y\hat{f}_{\star}(w_2)}).
\end{align*}
Now, let $\varepsilon=-1$: Then,
\begin{align*}
\alpha(\overline{\hat{f}_{\star}(w)})&=\alpha(\overline{\tilde{w}_1y^{\delta}\hat{f}_{\star}(x^ky^{-1}w_2)})\\
&=\alpha(\overline{\tilde{w}_1y^{\delta}x^ky^{-1}x^{-m}\hat{f}_{\star}(w_2)})\\
&=\alpha(\overline{(\tilde{w}_1y^{\delta}x^ky^{-1})x^{-m}(\tilde{w}_1y^{\delta}x^ky^{-1})^{-1}}\cdot\overline{\tilde{w}_1y^{\delta}x^ky^{-1}\hat{f}_{\star}(w_2)})\\
&=\alpha(\overline{w_1y^{-1}\hat{f}_{\star}(w_2)}).
\end{align*}
So the assertion $\alpha(\overline{\hat{f}_{\star}(w)})=\alpha(\overline{w_1y^{\varepsilon}\hat{f}_{\star}(w_2)})$ follows, and since we can also choose $w_2$ to be some part of the word $w$ after the last $y^{\pm1}$, i.e.\ $w_2=x^r$ for some $r\in\Z$, it follows that $\alpha(\overline{\hat{f}_{\star}(w)})=\alpha(\overline{w})$.\\
Since the diagram in Thm. \ref{Oridiagramm} commutes, we have $\alpha\circ f_*(\overline{w})=\alpha(\overline{w})$ and, therefore, $f\in\Mod_g(\alpha)$.\\
So we have $\Aff^+(X)\cap\Mod_g(\alpha)\neq\{1\}$, and by Thm. \ref{isomzuz}, it follows $\Aff^+(X)\cap\Mod_g(\alpha)\cong\Z$.
\end{Bew}

\begin{ko}
\label{schrg}
The statements of Thm. \ref{hsss}, \ref{zylhorimkern} and \ref{schnittnichtleer} apply also for vertical Schottky cut systems (which can be defined analogously to the horizontal ones). $\Aff^+(X)\cap\Mod_g(\alpha)$ is then generated by an element of the form $\begin{pmatrix}1&0\\c&1\end{pmatrix}$.
\end{ko}

\begin{satz}
\label{schief}
Let $O = (p\colon X^*\rightarrow E^*)$ be an origami of genus $g$, and let $v=(p,q)\in\Z^2$. Then there is a symplectic homomorphism $\alpha\colon \pi_1(X)\ra F_g$ and an $f\in\Aff^+(X)\cap\Mod_g(\alpha)$, $f\neq\id$ such that $A:=\der(f)$ has eigenvector $v$. Moreover, $A$ is parabolic.
\end{satz}
\begin{Bew}
If such an $f$ exists, then $A$ is parabolic by Thm. \ref{isomzuz}.\\
Let $v=(p,q)\in\Z^2$. W.l.o.g. choose $p$ and $q$ such that $\ggT(p,q)=1$ and $q\geq0$. For $q=0$ and $p=0$, we have already proven the assertion in Thm. \ref{schnittnichtleer} respectively Cor. \ref{schrg}. So we can assume that $q>0$ and $p\neq0$.\\
Let $G$ be the set of the straight lines in $\C$ with gradient $\frac{q}{p}$ that intersect the lattice $\Lambda_I:=\Z + i\Z$. Every unit square with vertices in the points of the lattice is intersected by $G$ in a way that its horizontal side is divided in $q$ equal segments and its vertical side is divided in $\vert p\vert$ equal segments. The same happens to the squares $Q_1,\dots,Q_d$ of $O$ if we look at $p^{-1}(h(G))$, where $p:X\ra E$ is the ramified covering from remark \ref{Orieinf}b and $h:\C\ra E$ is the covering of translation surfaces with $h(\Lambda_I)=\infty$:
$$\C\stackrel{h}{\longrightarrow}E\stackrel{p}{\longleftarrow}X$$
A horizontal cylinder that consists of $m$ squares is divided by $p^{-1}(h(G))$ in $mq$ parallelogramms. So we get a new partition of the translation surface $X$ in $dq$ parallelogramms with side length $\frac{1}{q}$ and $\frac{\sqrt{p^2+q^2}}{q}$ and height 1, where the lower/right side of a parallelogram is glued to the upper/left side of another parallelogram. The side of length $\frac{\sqrt{p^2+q^2}}{q}$ is parallel to the eigenvector $v$.\\
\begin{figure}[h]
\begin{center}
\includegraphics[scale=0.47]{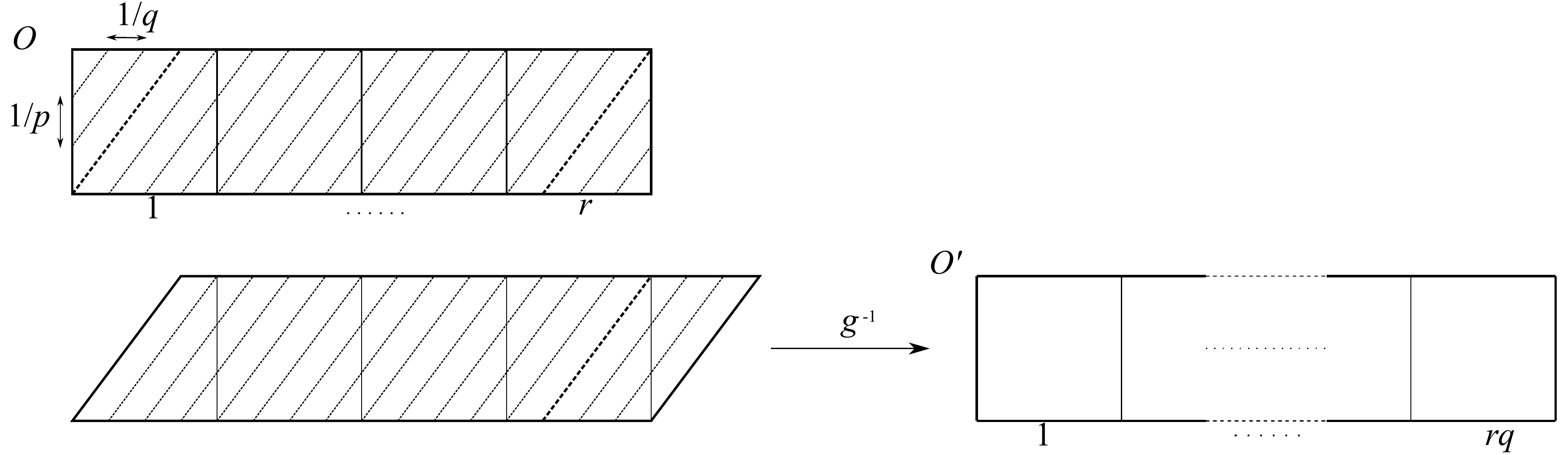}
\end{center}
\end{figure}\\
Now, let $O'$ be the origami that consist of $dq$ squares that are glued
like the parallelogramms in $O$. We denote by $X'$ the translation surface defined by $O'$ and look at the affine map
\begin{align*}
g:\; X'&\longrightarrow X\\
z&\longmapsto
\begin{pmatrix}
\frac{1}{q}&\frac{p}{q}\\
0&1
\end{pmatrix}\cdot z,
\end{align*}
that maps the unit squares of $O'$ to the parallelogramms in $O$. By Cor. \ref{schrg}, there is a set of symplectic generators $(a_1,b_1,\dots,a_g,b_g)$ of $\pi_1(X',x_0)$ such that $a_1,\dots,a_g$ are conjugate vertical, a symplectic homomorphism 
\begin{align*}
\alpha':\;\pi_1(X',x_0)&\longrightarrow F_g=\langle\gamma_1,\dots,\gamma_g\rangle\\
\quad a_i\mapsto1,\;\;&b_i\longmapsto\gamma_i
\end{align*}
and an $f'$ with $\langle f'\rangle=\Aff^+(X')\cap\Mod_g(\alpha')$ and $\der(f')=\begin{pmatrix}1&0\\c&1\end{pmatrix}$ for some $c\neq0$.\\
With $\alpha:=\alpha'\circ g^{-1}_*$, in the following diagram the outer quadrangle, the upper and the lower part commute. Since $f'\in\Mod_g(\alpha')$, also the left part commutes up to an inner automorphism of $F_g$. Therefore, also the right part commutes up to an inner automorphism of $F_g$, so that we have $g\circ f'\circ g^{-1}\in\Aff^+(X)\cap\Mod_g(\alpha)$. Since the diagram is symmetric in $X$ and $X'$, it also follows that $\Aff^+(X)\cap\Mod_g(\alpha)=\langle g\circ f'\circ g^{-1}\rangle$.
\[
\begin{xy}
\xymatrix{
\pi_1(X')\ar[rr]^{g_*}\ar[dd]_{f'_*}\ar[rd]_{\alpha'}&&\pi_1(X)\ar[dd]^{g_*\circ f'_*\circ g^{-1}_*}\ar[ld]^{\alpha}\\
&F_g&\\
\pi_1(X')\ar[rr]_{g_*}\ar[ru]^{\alpha'}&&\pi_1(X)\ar[lu]_{\alpha}}
\end{xy}
\]
Moreover, we have
\begin{align*}
\der(g\circ f'\circ g^{-1})\cdot\begin{pmatrix}p\\q\end{pmatrix}&=\der(g)\cdot\der(f')\cdot\der(g^{-1})\cdot\begin{pmatrix}p\\q\end{pmatrix}\\
&=\begin{pmatrix}\frac{1}{q}&\frac{p}{q}\\0&1\end{pmatrix}\cdot
\begin{pmatrix}1&0\\c&1\end{pmatrix}\cdot
\begin{pmatrix}q&-p\\0&1\end{pmatrix}\cdot
\begin{pmatrix}p\\q\end{pmatrix}\\
&=\begin{pmatrix}1+pc&-\frac{p^2c}{q}\\cq&-cp+1\end{pmatrix}\cdot
\begin{pmatrix}p\\q\end{pmatrix}=\begin{pmatrix}p\\q\end{pmatrix}.
\end{align*}
So $f:=g\circ f'\circ g^{-1}$ satisfies the assertion.
\end{Bew}

\begin{bem}
The proof of Thm. \ref{schief} also implies that the Veech group of an origami contains for every direction $(p,q)\in\Z^2$ a parabolic element $A$ with corresponding eigenvector. This is already known and follows from the fact that $\Gamma(O)$ has finite index in $\SL_2(\Z)$.
\end{bem}

\chapter{Flat surfaces}
\label{kaphalbtrans}
The statements of chapter \ref{kapori} about horizontal Schottky cut systems (in particular Thm. \ref{hsss} and \ref{zylhorimkern}) can be extended to flat surfaces (see definition \ref{flach}). The same holds for the statement that there is always a symplectic homomorphism $\alpha$ such that the group $\Aff^+(X,\mu)\cap\Mod_g(\alpha)$ is not trivial if $\overline{\Gamma}(X,\mu)$ contains a parabolic element. For translation surfaces, it can also be shown that in this case this group is isomorphic to $\Z$ (analogous to Thm. \ref{schnittnichtleer} for origamis). First of all, we'll transfer the concept of HSS to flat surfaces and we summarise some facts about differential forms and \v{C}ech cohomology that we'll need later.

\section{$\Aff^+(X,\mu)\cap\Mod_g(\alpha)\neq\{\id\}$ holds for an adequate $\alpha$}

\begin{definition}
Let $q$ be a holomorphic quadratic differential on a Riemann surface $X$ of genus $g\geq2$, and let $\mu$ the flat structure induced by $q$ on $X^*=X\setminus\{\text{zeros of }q\}$ (see remark \ref{hqm}).
\begin{enumerate}
\item We call a maximal real curve in $X^*$ that is locally mapped by a chart of $\mu$ to a horizontal segment (or a segment with direction $v\in\R^2$) a \begriff{horizontal trajectory} respectively \begriff{trajectory with direction $v$}.
\item A \begriff{critical trajectory} is one that ends in $X\setminus X^*$. Otherwise we call it a \begriff{regular trajectory}.
\item If all regular trajectories in one direction $v\in\R^2$ are closed, we call this direction a \begriff{Strebel direction}.
\end{enumerate}
\end{definition}

\begin{satz}
Let $q$ and $X$ be given as above and assume that the horizontal direction of $q$ is a Strebel direction. Then $(X,\mu)$ can be realized as patchwork of rectangles whose horizontal sides are in Strebel direction. These rectangles yield a decomposition of $X$ in horizontal cylinders.
\end{satz}

\begin{Bew}
See \cite[§4.1.1]{frankgabi}.
\end{Bew}

\begin{definition}
Let $X$ be a Riemann surface $X$ with flat structure $\mu$ realized as patchwork of rectangles. A \begriff{horizontal path} (with respect to the given realization) is a closed path $c$ with the property that for every $\varepsilon>0$ there exists a path $c_\varepsilon$ homotopic to $c$ that lies in a $\varepsilon$-neighbourhood of the horizontal sides of the rectangles. All the definitions of the ``horizontal notions'', like ``conjugate horizontal'', ``(conjugate) horizontal as element in $\pi_1(X)$'' and ``HSS'' can now be carried over from origami to arbitrary flat surfaces.
\end{definition}

\begin{bem}
\label{eckigesvertretersystem}
Let $(X,\mu)$ be a flat surface, realized as patchwork of rectangles. There is for each element $\gamma\in\pi_1(X)$ some representative that intersects the boundaries of the cylinders orthogonally, doesn't pass through critical points and follows vertical and horizontal segments in the cylinders. 
\end{bem}

\begin{satz}
\label{alleparab}
Let $(X,\mu)$ be a flat surface of genus $g$ whose projective Veech group $\overline{\Gamma}(X,\mu)$ contains a parabolic element $A$, and let $f\in\Aff^+(X,\mu)$ with $A=D(f)$. Then there is a symplectic homomorphism $\alpha:\pi_1(X)\ra F_g$ and an $n\in\N$ with
\[f^n\in\Mod_g(\alpha)\cap\Aff^+(X,\mu).\] 
\end{satz}

\begin{Bew} 
By \cite[Prop. 2.4]{veech}, the eigendirection of $A$ is a Strebel direction (let $v$ be a vector in this direction). Then we can realize $(X,\mu)$ as a patchwork of rectangles that have two sides parallel to $v$. We assume w.l.o.g. that these sides are horizontal, i.e.\ that $v=(1,0)^\top$. As in Thm. \ref{hsss}, it follows that there is a HSS with respect to the given realization of $X$ as patchwork of rectangles. Let $(a_1,b_1,\dots,a_g,b_g)$ be a set of symplectic generators of the fundamental group $\pi_1(X)$ such that $a_1,\dots,a_g$ are conjugate horizontal, and let $\alpha:\pi_1(X)\ra F_g=\langle\gamma_1,\dots,\gamma_g\rangle$ be defined by $\alpha(a_i)=1$ and $\alpha(b_i)=\gamma_i$ for $i\in\{1,\dots,g\}$. As in Thm. \ref{zylhorimkern} it follows that an element in $\pi_1(X)$ that is freely homotopic to a horizontal cut around a cylinder lies in the normal subgroup generated by $a_1,\dots,a_g$ because in the proofs of these theorems only the flat structure of origamis has been used. A power $f^n$ (for an $n\in\N$) of $f$ is a multiple Dehn twist on the middle of the cylinders. With remark \ref{eckigesvertretersystem} it is easy to see that $f^n$ acts on the fundamental group in a way that only adds horizontal loops going around the cylinders to a given element in $\pi_1(X)$. These horizontal curves are contained in the kernel of $\alpha$. The assertion follows.
\end{Bew}

\begin{ko}
\label{insbesgitter}
In particular, the assertion of the theorem holds if the projective Veech group $\overline{\Gamma}(X,\mu)$ is a lattice in $\PSL_2(\R)$.
\end{ko}

\begin{Bew}
If $\overline{\Gamma}(X,\mu)$ is a lattice in $\PSL_2(\R)$, there exists an $f\in\Aff^+(X,\mu)$ such that $A:=D(f)$ is parabolic, see the remark after \cite[Prop. 2.10]{veech}.
\end{Bew}

\section{Differential forms}
To prove a more general version of Thm. \ref{nichthyp} for translation surfaces, we first need some theorems about cohomology of compact Riemann surfaces that we state without proofs and which can be found in \cite{forster}. We also refer to that source for the definitions and theorems in this section and in section \ref{cech}.\\\\
Let $a$ be a point on a Riemann surface $X$, and denote by $\mathscr E$
the sheaf of real differentiable and complex valued functions on $X$. We
denote by $\mathscr E_a$ the stalk of $\mathscr E$ in the point $a$. So $\mathscr E_a$ consists of the differentiable germs of functions in $a$. Let $\textfrak m_a\subset\mathscr E_a$ be the maximal ideal of germs of functions that vanish in $a$, and let $\textfrak m^2_a\subset\textfrak m_a$ be the ideal of germs of functions that vanish in $a$ with order 2 at least. The \begriff{cotangent space} of $X$ at the point $a$ is defined as
\[T^{(1)}_a:=\textfrak m_a/\textfrak m^2_a.\]
For an open (chart-)neighbourhood $U$ of $a$ and $f\in\mathscr E(U)$, the differential $d_af\in T^{(1)}_a$ is defined as
\[d_af:=(f-f(a))\mod\mathfrak m_a^2.\]
It can be shown \cite[9.4]{forster} that
\[d_af=\frac{\partial f}{\partial x}(a)d_ax+\frac{\partial f}{\partial y}(a)d_ay=\frac{\partial f}{\partial z}(a)d_az+\frac{\partial f}{\partial\overline{z}}(a)d_a\overline{z}.\]
Therefore, $\{d_ax,d_ay\}$ and $\{d_az,d_a\overline{z}\}$ each form a base of $T^{(1)}_a$ as $\C$-vector space.\\
Since on a Riemann surface the transition maps are biholomorphic, the vector subspaces $T^{(1,0)}_a:=\C d_az$ and $T^{(0,1)}_a:=\C d_a\overline{z}$ do not depend on the chart. We then speak of differential forms of type (1,0) respectively type (0,1). By construction, we have
\[T^{(1)}_a=T^{(1,0)}_a\oplus T^{(0,1)}_a.\]

\begin{definition}
\begin{enumerate}
\item A \begriff{differential form of 1. order} (below we simply refer to them as ``differential forms'') on a Riemann surface $X$ is a map
\[\omega\colon X\ra\bigcup_{a\in X}T^{(1)}_a\]
with $\omega(a)\in T^{(1)}_a$ for all $a\in X$.
\item For $f\in\mathscr E(X)$, the differential forms $df, d'f, d''f$ are defined by:
\[(df)(a):=d_af;\quad\quad(d'f)(a):=d_a'f:=\frac{\partial f}{\partial z}(a)d_az;\quad\quad(d''f)(a):=d_a''f:=\frac{\partial f}{\partial\overline{z}}(a)d_a\overline{z}\]
\item A differential form $\omega$ is called \begriff{differentiable} or \begriff{holomorphic} if they can be realised with respect to every chart $(U,z)$ as
\[\omega=fdz+gd\overline{z}\quad\text{ in } U\text{ with } f,g\in\mathscr E(U)\]
respectively
\[\omega=fdz\quad\text{ in } U\text{ with } f\in\mathscr O(U),\]
where $\mathscr O$ is the sheaf of holomorphic functions.
\item Let $U\subseteq X$ be open. We denote by $\mathscr E^{(1)}(U)$ the $\mathscr E(U)$-module of differentiable differential forms on $U$, by $\mathscr E^{(1,0)}(U)$ and $\mathscr E^{(0,1)}(U)$ the submodules of type (1,0) respectively (0,1), and by $\Omega(U)$ the $\mathscr O(U)$-module of holomorphic differential forms on $U$.\\
We denote the sheafs of modules on $X$ that we get from $\mathscr E^{(1)}(U)$, $\mathscr E^{(1,0)}(U)$, $\mathscr E^{(0,1)}(U)$ and $\Omega(U)$ by
\[\mathscr E^{(1)},\quad\quad\mathscr E^{(1,0)},\quad\quad\mathscr E^{(0,1)}\quad\text{ and }\quad\Omega.\]
\item For $U\subseteq X$, let $\omega=fdz+gd\overline{z}\in\mathscr E^{(1)}(U)$ with $f,g\in\mathscr E(U)$. We define the differential form $\overline{\omega}:=\overline{f}d\overline{z}+\overline{g}dz\in\mathscr E^{(1)}(U)$ and 
\[\overline{\Omega}(U):= \{\omega\in\mathscr E^{(1)}(U)\mid\overline{\omega}\in\Omega(U)\}.\]
Let $\overline{\Omega}$ be the sheaf defined so. Differential forms in $\overline{\Omega}(U)$ are called \begriff{antiholomorphic}.
\end{enumerate}
\end{definition}

\begin{definition}
Let $X$ and $Y$ be two Riemann surfaces, $f\colon X\ra Y$ a differentiable map and $\omega\in\mathscr E^{(1)}(Y)$. For a local representation $\omega=gdz+hd\overline{z}$, the \begriff{pullback} of $\omega$ with respect to $f$ in $\mathscr E^{(1)}(X)$ is given by 
\[f^*\omega:= (g\circ f)d(z\circ f)+(h\circ f)d(\overline{z}\circ f).\]
These local representations can be glued on $X$ such that they define a differential form on $X$, see \cite[9.17]{forster}.
\end{definition}

\section{\v{C}ech cohomology}
\label{cech}
\begin{definition}
Let $X$ be a topological space and $\mathscr F$ a sheaf of abelian groups on it.
\begin{enumerate}
\item For a covering $\textfrak U=(U_i)_{i\in I}$, we define for $q\in\N_0$ the \begriff{$q$-th cochain group} of $\mathscr F$ with respect to $\textfrak U$ as
\[C^q(\textfrak U,\mathscr F):=\prod_{(i_0,\dots,i_q)\in I^{q+1}}\mathscr F(U_{i_0}\cap\dots\cap U_{i_q}).\]
Its elements are called \begriff{$q$-cochains}. The sum in $C^q(\textfrak U,\mathscr F)$ is defined component wise.
\item For $(f_i)_{i\in I}\in C^0(\textfrak U,\mathscr F)$ and $(f_{ij})_{i,j\in I}\in C^1(\textfrak U,\mathscr F)$, let the \begriff{coboundary operators}
\[\delta\colon C^0(\textfrak U,\mathscr F)\ra C^1(\textfrak U,\mathscr F)\quad\text{ and }\quad \delta\colon C^1(\textfrak U,\mathscr F)\ra C^2(\textfrak U,\mathscr F)\]
be the homomorphisms which are defined by
\[\delta((f_i)_{i\in I})=(g_{ij})_{i,j\in I}\;\text{ with }\;g_{ij}:=f_j\mid_{U_i\cap U_j}-f_i\mid_{U_i\cap U_j}\;\text{ respectively }\]
\[\delta((f_{ij})_{i,j\in I})=(g_{ijk})_{i,j,k\in I}\;\text{ with }\;g_{ijk}:=f_{jk}\mid_{U_i\cap U_j\cap U_k}-f_{ik}\mid_{U_i\cap U_j\cap U_k}+f_{ij}\mid_{U_i\cap U_j\cap U_k}.\]
Note that $\delta\circ\delta\colon C^0(\textfrak U,\mathscr F)\ra C^2(\textfrak U,\mathscr F)$ maps everything to zero.
\item The group of \begriff{1-cocycles} is defined as
\[Z^1(\textfrak U,\mathscr F):=\Kern(C^1(\textfrak U,\mathscr F)\stackrel{\delta}{\ra}C^2(\textfrak U,\mathscr F))\]
and the group of \begriff{1-coboundaries} is defined as
\[B^1(\textfrak U,\mathscr F):=\Bild(C^0(\textfrak U,\mathscr F)\stackrel{\delta}{\ra}C^1(\textfrak U,\mathscr F)).\]
\item The \begriff{1. cohomology group} of $\mathscr F$ with respect to $\textfrak U$ is defined as
\[H^1(\textfrak U,\mathscr F):=Z^1(\textfrak U,\mathscr F)/B^1(\textfrak U,\mathscr F).\]
The factor group is well defined because every coboundary is also cocycle.
\end{enumerate}
\end{definition}
It is possible to refine the coverings of $X$ more and more, where $\textfrak V$ is called a \begriff{refinement} of $\textfrak U$ (write $\textfrak V\leq\textfrak U$) if every open set of $\textfrak V$ is a subset of an open set of $\textfrak U$. If we choose for every $V\in\textfrak V$ some $U_V\in\textfrak U$ with $V\subseteq U_V$, the corresponding restriction morphisms of the sheaf $\mathscr F$ induce a homomorphism
\[t^{\textfrak U}_{\textfrak V}\colon H^1(\textfrak U,\mathscr F)\ra H^1(\textfrak V,\mathscr F).\]
It can be shown that $t^{\textfrak U}_{\textfrak V}$ is injective \cite[12.4]{forster} and does not depend on the choice of $U_V$ \linebreak\cite[12.3]{forster}. Moreover, we have
$\textfrak W\leq\textfrak V\leq\textfrak U$:
\[t^{\textfrak V}_{\textfrak W}\circ t^{\textfrak U}_{\textfrak V}=t^{\textfrak U}_{\textfrak W}.\]
Therefore, we have an equivalence relation $\sim$ on the disjoint union
\[\widetilde{H^1}(X,\mathscr F):=\coprod_{\textfrak U\text{ open covering of } X}H^1(\textfrak U,\mathscr F),\]
that is defined by 
\[\chi\sim\xi:\Leftrightarrow\exists\;\textfrak V\text{ with }\textfrak V\leq\textfrak U, \textfrak V\leq\textfrak U'\colon t^{\textfrak U}_{\textfrak V}(\chi)=t^{\textfrak U'}_{\textfrak V}(\xi)\]
for $\chi\in H^1(\textfrak U,\mathscr F), \xi\in H^1(\textfrak U',\mathscr F)$.
\begin{definition} Let $X$ be a topological space with a sheaf $\mathscr F$ on it. The \begriff{1. cohomology group} of $X$ with respect to $\mathscr F$ is given by
\[H^1(X,\mathscr F)=\FakRaum{\widetilde{H^1}(X,\mathscr F)}{\sim} = \varinjlim H^1(\textfrak U,\mathscr F).\]
\end{definition}
Like the 1. cohomology group one can also define the $n$-th cohomology group for $n\in\N_0$. For us, it is interesting to look at the 0. cohomology group:
\begin{definition} Let $X$ be a topological space with a sheaf $\mathscr F$ on it. For a covering $\textfrak U$, define
\begin{align*}
Z^0(\textfrak U,\mathscr F)&:=\Kern(C^0(\textfrak U,\mathscr F)\stackrel{\delta}{\ra}C^1(\textfrak U,\mathscr F))\\
B^0(\textfrak U,\mathscr F)&:=0\\
H^0(\textfrak U,\mathscr F)&:=Z^0(\textfrak U,\mathscr F)/B^0(\textfrak U,\mathscr F)=Z^0(\textfrak U,\mathscr F).
\end{align*}
It can easily be shown that $H^0(\textfrak U,\mathscr F)\cong\mathscr F(X)$, independently of the covering \cite[12.10]{forster}. Therefore, we can define
\[H^0(X,\mathscr F):=\mathscr F(X)\]
as the \begriff{0. cohomology group} of $X$ with values in $\mathscr F$.
\end{definition}
\begin{bem}
If the sheaf $\mathscr F$ is a sheaf of vector spaces, then $H^1(X,\mathscr F)$ and $H^0(X,\mathscr F)$ are vector spaces.
\end{bem}
Given the above definitions, we now cite some theorems that we will need below.
\begin{satz}
\label{einfzus}
Let $X$ be a simply connected Riemann surface. Then we have
\[H^1(X,\Z)=H^1(X,\C)=0\]
where $\Z$ and $\C$ are the sheafs of locally constant functions with values in $\Z$ respectively $\C$.
\end{satz}

\begin{Bew}
See \cite[12.7]{forster}.
\end{Bew}

\begin{satz}[Leray]~\\
\label{leray} Let $X$ be a topological space with sheaf $\mathscr F$ and let $\textfrak U$ be an open covering of $X$ with $H^1(U,\mathscr F)=0$ for all $U\in\textfrak U$. Then we have:
\[H^1(X,\mathscr F)\cong H^1(\textfrak U,\mathscr F).\]
The covering $\textfrak U$ is called \begriff{Leray covering}.
\end{satz}

\begin{Bew} See \cite[12.8]{forster}.
\end{Bew}

\begin{defbem}
Let $X$ be a topological space and $\mathscr F$, $\mathscr G$ and $\mathscr H$ sheafs on $X$.
\begin{enumerate}
\item A sheaf homomorphism $\alpha\colon \mathscr F\ra\mathscr G$ is a family of homomorphisms
\[\alpha_U\colon \mathscr F(U)\ra\mathscr G(U), U\subseteq X \text{ open},\]
that are compatible with the restriction morphisms, i.e.\ the following diagram commutes for all open sets with $V\subseteq U$:
\[\begin{xy}
\xymatrix{\mathscr F(U)\ar[r]^{\alpha_U}\ar[d]&\mathscr G(U)\ar[d]\\
\mathscr F(V)\ar[r]^{\alpha_V}&\mathscr G(V)}   
\end{xy}
\]
\item $\alpha$ induces a homomorphism of the stalks $\alpha_x\colon \mathscr F_x\ra\mathscr G_x$ for all $x\in X$.
\item A sequence $\mathscr F\stackrel{\alpha}{\ra}\mathscr G\stackrel{\beta}{\ra}\mathscr H$ of sheaf homomorphisms is called \begriff{exact} if 
\[\mathscr F_x\stackrel{\alpha_x}{\longrightarrow}\mathscr G_x\stackrel{\beta_x}{\longrightarrow}\mathscr H_x\]
is exact for all $x\in X$.
\end{enumerate}
\end{defbem}
A sheaf homomorphism $\alpha:\mathscr F\ra\mathscr G$ induces homomorphisms of the cohomology groups (see also \cite[15.10]{forster})
\begin{align*}
\alpha^0:\;& H^0(X,\mathscr F)\ra H^0(X,\mathscr G) \text{ and }\\
\alpha^1:\;& H^1(X,\mathscr F)\ra H^1(X,\mathscr G). 
\end{align*}
$\alpha^0$ is defined as the homomorphism $\alpha_X\colon \mathscr F(X)\ra\mathscr G(X)$ and $\alpha^1$ is defined as follows:\\
Let $\textfrak U$ be a covering of $X$. Then we get a homomorphism
\[\alpha_{\textfrak U}:\;C^1(\textfrak U,\mathscr F)\ra C^1(\textfrak U,\mathscr G),\]
by mapping a cochain $(f_{ij})$ to the cochain $(\alpha(f_{ij}))$. This induces a homomorphism of the cohomologies with respect to $\textfrak U$ and this in turn induces the homomorphism
\[\alpha^1:\;H^1(X,\mathscr F)\ra H^1(X,\mathscr G).\]
\begin{satz}
\label{exakt}
Let $X$ be a topological space with sheafs $\mathscr F$, $\mathscr G$ and $\mathscr H$ and sheaf homomorphisms $\alpha:\mathscr F\ra\mathscr G$ and $\beta:\mathscr G\ra\mathscr H$, where the sheaf sequence
\[0\ra\mathscr F\stackrel{\alpha}{\ra}\mathscr G\stackrel{\beta}{\ra}\mathscr H\ra0\]
is exact. Then the following induced sequence is also exact:
\[0\ra H^0(X,\mathscr F)\stackrel{\alpha^0}{\ra} H^0(X,\mathscr G)\stackrel{\beta^0}{\ra} H^0(X,\mathscr H)\stackrel{\delta^*}{\ra} H^1(X,\mathscr F)\stackrel{\alpha^1}{\ra} H^1(X,\mathscr G)\stackrel{\beta^1}{\ra} H^1(X,\mathscr H)\]
\end{satz}

\begin{beme}
The homomorphism $\delta^*$ is defined in \cite[15.11]{forster}. For us, it is enough to know that such a homomorphism exists.
\end{beme}

\begin{Bew}
See \cite[15.12]{forster}.
\end{Bew}
For what follows, we have to introduce the notion of harmonic differential forms. One possibility to define them is the following (for another possibility see \cite[19.2]{forster} and \cite[19.3]{forster}):

\begin{definition}
\begin{enumerate}
\item $\omega\in\mathscr E^{(1)}(X)$ is called \begriff{harmonic} if $\omega=\omega_1+\omega_2$ with $\omega_1\in\Omega(X)$ and $\omega_2\in\overline{\Omega}(X)$
\item We denote the vector space of harmonic differential forms with $\Harm^1(X)$.
\end{enumerate}
\end{definition}

\begin{bem}
$\Harm^1(X)$ is a $2g$ dimensional vectorspace over $\C$ and we have
\[\Harm^1(X)=\Omega(X)\oplus\overline{\Omega}(X).\]
\end{bem}
\begin{Bew}
This follows from \cite[19.6 and 19.11]{forster}
\end{Bew}

\begin{satz}[de Rahm and Hodge]~\\
\label{de Rahm and Hodge}
Let $X$ be a compact Riemann surface of genus $g$. We have
\[\Harm^1(X)\cong H^1(X,\C).\]
Therefore, $H^1(X,\C)$ is isomorphic to $\C^{2g}$ as $\C$-vector space,.
\end{satz}
\begin{Bew}
See \cite[19.14]{forster}.
\end{Bew}

\section{For translation surfaces, $\Aff^+(X,\mu)\cap\Mod_g(\alpha)$ is cyclic}

\begin{lemma}
\label{cohom}
Let $X$ be a compact Riemann surface of genus $g$ with sheafs $\Z$ and $\C$ and let $\iota\colon \Z\hookrightarrow\C$ be the sheaf homomorphism induced by the natural embedding.\\
The map $\iota^1\colon H^1(X,\Z)\ra H^1(X,\C)$ induced by $\iota$ is injective and maps every set of generators of $H^1(X,\Z)$ to a base of $H^1(X,\C)$.
\end{lemma}

\begin{Bew}
Let  $\C^*$ be the sheaf of locally constant functions to $\C^*$, and let $\exp\colon \C\ra\C^*$ be the sheaf homomorphism induced by $z\mapsto e^{2\pi iz}$. We have the following exact sequence of sheafs:
\[
0\longrightarrow\Z\stackrel{\iota}{\longrightarrow}\C\stackrel{\exp}{\longrightarrow}\C^*\longrightarrow0
\]
By Thm. \ref{exakt}, this induces a long exact sequence of groups:
\[
0\ra H^0(X,\Z)\stackrel{\iota^0}{\ra} H^0(X,\C)\stackrel{\exp^0}{\ra} H^0(X,\C^*)\stackrel{\delta^*}{\ra} H^1(X,\Z)\stackrel{\iota^1}{\ra} H^1(X,\C)\stackrel{\exp^1}{\ra} H^1(X,\C^*)
\]
$H^0(X,\Z)$, $H^0(X,\C)$ and $H^0(X,\C^*)$ are the groups of constant functions from $X$ to $\Z$, $\C$ and $\C^*$ with addition in $\Z$ and in $\C$ respectively multiplication in $\C^*$ as group operations. Therefore, they are isomorphic to the groups $\Z$, $\C$ and $\C^*$ in a canonical way; moreover we have $\iota^0=\iota$ and $\exp^0=\exp$. Therefore, $\exp^0$ is surjective and the exactness of the sequence implies that $\delta^*$ maps everything to zero and $\iota^1$ is injective.\\
Now, let w.l.o.g. $g\geq1$. We can construct a finite covering $\textfrak U$ of $X$ of simply connected charts with the property that every two elements in $\textfrak U$ have connected (possibly empty) intersection as follows:\\
On $X$ we have a metric induced by the universal covering. Let $\ell(X)$ be the length of a shortest closed geodesic on $X$. For every $p\in X$, let \[U_p:=B_{\frac{1}{4}\ell(X)}(p):=\{x\in X\mid d(x,p)<\frac{1}{4}\ell(X)\}.\]
The covering $\{U_p\mid p\in X\}$ of $X$ satisfies the desired property. Since $X$ is compact, we need only finitely many of these neighbourhoods $U_1,\dots,U_n$, to cover $X$. Let
\[\textfrak U:=\{U_i\mid i=1,\dots,n\}.\]
$\textfrak U$ is a Leray covering since $U_1,\dots,U_n$ are simply connected and, therefore, by Thm. \ref{einfzus} we have
\[H^1(U_i,\C)=H^1(U_i,\Z)=0,\quad\forall U_i\in\textfrak U.\]
By the Theorem of Leray (\ref{leray}), we have
\[H^1(\textfrak U,\C)\cong H^1(X,\C)\quad\text{ and }\quad H^1(\textfrak U,\Z)\cong H^1(X,\Z),\]
where the isomorphisms are canonical.\\
Because of condition (*) the locally constant maps from $U_i\cap U_j$ to $\C$ respectively $\Z$ are also globally constant. Therefore, the cochain groups $C^1(\textfrak U,\Z)\leq C^1(\textfrak U,\C)$ can be identified with $\Z^{n^2}$ and $\C^{n^2}:=(c_{ij})_{i,j\in\{1,\dots,n\}}$. $C^1(\textfrak U,\C)$ is a $\C$-vector space. The sub-vectorspace of cocycles $Z^1(\textfrak U,\C)$ consists of those cochains that satisfy the cocycle relations, i.e.\
\[\forall i,j,k\colon U_i\cap U_j\cap U_k\neq\emptyset\colon c_{jk}-c_{ik}+c_{ij}=0.\]
It follows that $Z^1(\textfrak U,\C)$ is the solution space of a system of linear equations with coefficients in $\Z$ and is, therefore, generated by elements in $C^1(\textfrak U,\Z)\cap Z^1(\textfrak U,\C)=Z^1(\textfrak U,\Z)$.\\
The embedding $C^1(\textfrak U,\Z)\hookrightarrow C^1(\textfrak U,\C)$ induces a map $\iota_{\textfrak U}^1\colon H^1(\textfrak U,\Z)\hookrightarrow H^1(\textfrak U,\C)$ that induces in turn the embedding $\iota^1$. Since $H^1(\textfrak U,\Z)\cong H^1(X,\Z)$, $\iota_{\textfrak U}^1$ is injective, too. By definition of $\iota_{\textfrak U}^1$, the following diagram commutes:
\[\begin{xy}
\xymatrix{
Z^1(\textfrak U,\Z)\ar@{->>}[d]\ar@{^(->}[r]&Z^1(\textfrak U,\C)\ar@{->>}[d]\\
H^1(\textfrak U,\Z)\ar@{^(->}[r]^{\iota_{\textfrak U}^1}&H^1(\textfrak U,\C)
}
\end{xy}\]
Therefore, also $H^1(\textfrak U,\C)$ has the embedding of a set of generators of $H^1(\textfrak U,\Z)$ as base. Since $\iota^1$ is induced by $\iota_{\textfrak U}^1$, the corresponding assertion applies also to $H^1(X,\Z)$ and $H^1(X,\C)$.
\end{Bew}

\begin{lemma}
\label{z2g}
Let $X$ be a Riemann surface of genus $g$. Then the group $H^1(X,\Z)$ is isomorphic to $\Z^{2g}$. In other words: $H^1(X,\Z)\cong\pi_1^{ab}(X)$.
\end{lemma}

\begin{Bew}
This follows directly from Lemma \ref{cohom} and Thm. \ref{de Rahm and Hodge}.
\end{Bew}

\begin{satz}
\label{nichthypallgemein}
Let $(X,\mu)$ be a translation surface and $f\in\Aff^+(X,\mu)$ with $A := \der(f)$. Then the eigenvalues of $A$ are also eigenvalues of $M_f$.
\end{satz}

\begin{Bew}
By Lemma \ref{cohom} and Lemma \ref{z2g}, we can extend $M_f\in\Aut(\pi_1^{ab}(X))$ uniquely to a (bijective) linear map 
\[M_f\colon H^1(X,\C)\ra H^1(X,\C).\]
By the Theorem of de-Rahm and Hodge, we have $H^1(X,\C)\cong\Omega(X)\oplus\overline{\Omega}(X)$. The differentials of the form $dz$ on the charts can be glued because of the translation structure such that we can talk about a differential $\omega:= dz$ on $X$. The same holds for $\overline{\omega}:=d\overline{z}$. We have $\omega\in\Omega(X)$ and $\overline{\omega}\in\overline{\Omega}(X)$. Moreover, $\omega$ and $\overline{\omega}$ are linearly independent and span a 2-dimensional $\C$-vector space $[\omega,\overline{\omega}]\subseteq H^1(X,\C)$. We show that $M_f$ maps this vector space to itself and that $M_f\vert_{[\omega, \overline{\omega}]}$ has the same eigenvalues as $A$. These are, hence, eigenvalues of $M_f$, too.\\
To get $M_f\vert_{[\omega, \overline{\omega}]}\in\SL_2(\C)$, we have to take a look at $f^*(dz)$ and $f^*(d\overline{z})$ with
\begin{align*}
f\;:\;(X,\mu)&\ra(X,\mu)\\
\begin{pmatrix}x\\y\end{pmatrix}&\mapsto\begin{pmatrix}\alpha&\beta\\\gamma&\delta\end{pmatrix}\cdot\begin{pmatrix}x\\y\end{pmatrix}=\begin{pmatrix}\alpha x+\beta y\\\gamma x+\delta y\end{pmatrix},
\end{align*}
and $z=x+iy$. With $2x=z+\overline{z}$ and $2y=-i(z-\overline{z})$ we can write $f$ as follows:
\[z\mapsto\frac{1}{2}[\alpha(z+\overline{z})-i\beta(z-\overline{z})+i\gamma(z+\overline{z})+\delta(z-\overline{z})]\]
Therefore, we have
\[f^*(dz)=d(f(z))=\frac{1}{2}[(\alpha-i\beta+i\gamma+\delta)dz+(\alpha+i\beta+i\gamma-\delta)d\overline{z}]\]
\[f^*(d\overline{z})=d(\overline{f(z)})=\frac{1}{2}[(\alpha-i\beta-i\gamma-\delta)dz+(\alpha+i\beta-i\gamma+\delta)d\overline{z}]\]
Therefore, we have
\[M_f\vert_{[\omega,\overline{\omega}]}=\frac{1}{2}\begin{pmatrix}
\alpha-i\beta+i\gamma+\delta&\alpha+i\beta+i\gamma-\delta\\
\alpha-i\beta-i\gamma-\delta&\alpha+i\beta-i\gamma+\delta
\end{pmatrix}.
\]
The characteristic polynomial of $M_f\vert_{[\omega,\overline{\omega}]}$ is the same as the one of $A$ because
\[\frac{1}{4}[(\alpha+\delta-i(\beta-\gamma)-2x)(\alpha+\delta+i(\beta-\gamma)-2x)-(\alpha-\delta+i(\beta+\gamma))(\alpha-\delta-i(\beta+\gamma))]\]
\[=\alpha\delta-\beta\gamma-(\alpha+\delta)x+x^2=1-(\alpha+\delta)x+x^2.\]
Therefore, $A$ and $M_f\vert_{[\omega,\overline{\omega}]}$ have the same eigenvalues.
\end{Bew}

\begin{ko}
Let $(X,\mu)$ be a translation surface, $\alpha:\pi_1(X)\ra F_g$ a symplectic homomorphism and $f\in\Aff^+(X)\cap\Mod_g(\alpha)$. Then $\der(f)$ is not hyperbolic.
\end{ko}
\begin{Bew}
This follows directly from Thm. \ref{nichthypallgemein} and Lemma \ref{EW1} d.
\end{Bew}
\begin{ko}
\label{isomzuzallg}
If $\Aff^+(X)\cap\Mod_g(\alpha)$ is not trivial, we have $\Aff^+(X)\cap\Mod_g(\alpha)\cong\Z$.
\end{ko}
\begin{ko}
Let $(X,\mu)$ be a translation surface, whose Veech group is a lattice, and let \linebreak$f\in\Aff^+(X,\mu)$. Then there is a symplectic homomorphism $\alpha$ such that
\[f^n\in\Aff^+(X,\mu)\cap\Mod_g(\alpha) \quad \text{ for some } n\in\N,\]
if and only if $\der(f)$ is not elliptic and the only eigenvalues of matrix $M_f$ are roots of unity.
\end{ko}
\begin{Bew}
Let $f\in\Aff^+(X)$, and assume that $\der(f)$ is not elliptic and that $M_f$ has only roots of unity as eigenvalues. Two of these eigenvalues, $\lambda_1$ and $\lambda_2$, are also eigenvalues of $\der(f)$, according to Thm. \ref{nichthypallgemein}. Since $\lambda_1$ and $\lambda_2$ are roots of unity, we have $|\lambda_1|=|\lambda_2|=1$ and since $\der(f)$ is not elliptic, $\lambda_1$ and $\lambda_2$ are real and, hence, either $-1$ or $1$. Then $\der(f^2)$ has only the eigenvalue $1$ and is, therefore, parabolic (or the identity). By Thm. \ref{alleparab}, $f$ satisfies the assertion.\\
Let us now assume that $f\in\Aff^+(X,\mu)$ and that there is a symplectic homomorphism $\alpha$ such that
\[f^n\in\Aff^+(X,\mu)\cap\Mod_g(\alpha) \quad \text{ for some } n\in\N.\]
Since $\Gamma(X,\mu)$ is a lattice, elliptic elements in $\Gamma(X,\mu)$ have finite order and can, therefore, not be in $\Mod_g(\alpha)$ since this group is torsion free by Thm. \ref{torsionsfrei}. By Lemma \ref{EW1}, $f^n\in\Mod_g(\alpha)$ implies that $(M_f)^n=M_{f^n}$ has $1$ as its only eigenvalue. All the eigenvalues $\lambda_1,\dots\lambda_{2g}$ of $M_f$, therefore, satisfy $\lambda_1^n=\dots=\lambda_{2g}^n=1$. It follows that the eigenvalues of $M_f$ are roots of unity.
\end{Bew}
The following corollary states a well known fact, see e.g.\ \hbox{\cite[Lemma 2.3.17]{ollidiss}.}
\begin{ko}
\label{nichtintorelli}
Let $(X,\mu)$ be a translation surface of genus $g\geq2$ and $f\in\Aff^+(X)$ with $f\neq\id$. Then $f$ is not in the Torelli group (i.e.\ acts not trivially on the homology).
\end{ko}

\begin{Bew}
Let $f$ be in the Torelli group. Then we have, in particular, $M_f\vert_{[\omega,\overline{\omega}]}=I$ with the notations as in the proof of Thm.\ref{nichthypallgemein}. But it also follows from the proof of Thm. \ref{nichthypallgemein} that $M_f\vert_{[\omega,\overline{\omega}]}$ and $\der(f)$ have the same characteristic polynomial. Therefore, $\der(f)=I$, and $f$ is hence a translation. But translations have finite order, by the $84(g-1)$-Theorem (see \cite[Thm. 7.4]{farbmarg}). Since the Torelli group is torsion free (\cite[Thm. 6.12]{farbmarg}), we have a contradiction.
\end{Bew}

\section{The image of $\Delta$ in the Schottky space}
Let $(X,q)=(X,\mu)$ be a flat surface and let $\Delta_q$ be the Teichmüller disc that corresponds to the   quadratic holomorphic differential $q$. In summary, we can say that for every parabolic element $\tilde{\tau}\in\Aff^+(X,\mu)$, there exists a symplectic homomorphism $\alpha:\pi_1(X)\ra F_g$ such that
\[\Z\cong\langle\tilde{\tau}^n\rangle\leq\Aff^+(X,\mu)\cap\Mod_g(\alpha)\cong\Stab(\Delta_q)\cap\Mod_g(\alpha).\]
In particular, if $X$ is a translation surface then
\[\langle\tilde{\tau}^n\rangle\cong\Aff^+(X,\mu)\cap\Mod_g(\alpha)\cong\Stab(\Delta_q)\cap\Mod_g(\alpha).\]
The map $\pr\circ\iota \colon \mathds{H}\rightarrow \pr(\Delta_q)=:C \subseteq \mathcal{M}_g$ factorizes by Thm. \ref{rar} (as we saw in the proof of Thm. \ref{algkurveechgitter}) through $\overline{\Gamma}^*(X,\mu) = R\overline{\Gamma}(X,\mu)R^{-1}$. With the statements above we get the following commutative diagram (where $C$ is an algebraic curve if $\overline{\Gamma}(X,\mu)$ is a lattice in $\PSL_2(\R)$):
\[
\begin{xy}
\xymatrix{\mathds H\ar[rr]^{t \mapsto -\overline{t}}\ar[d] && \mathds H\ar[rr]^{\iota}\ar[d] && \Delta_q\ar[d]^{s_\alpha}\\
\mathds H/\langle\tau\rangle\ar[rr]^{\text{antihol.}}\ar[d]&&\mathds H/\langle\tau\rangle^*\ar[rr]^{\text{bihol.}}\ar[d]&& s_\alpha(\Delta_q)\ar[d]\\
\mathds H/\overline{\Gamma}(X,\mu)\ar[rr]^{\text{antihol.}}&&\mathds H/\overline{\Gamma}^*(X,\mu)\ar[rr]^{\text{birat.}} && C.}
\end{xy}
\]
\textbf{Summary}:\\
For translation surfaces, $s_\alpha(\Delta_q)$ is biholomorphic to $\mathds H$ or to $\FakRaum{\mathds H}{\langle\tau\rangle}$ for  $\tau=\begin{pmatrix}1&1\\0&1\end{pmatrix}.$

\chapter{An algorithm to find a HSS}
\label{kapalgorithmus}
\begin{definition}
Take a compact, orientable surface $X$ of genus $g$ and a finite set $\mathcal M$ of simply closed, pairwise disjoint curves on $X$.
A subset $\mathcal M'$ of $\mathcal M$ such that
\begin{itemize}
 \item $X\setminus\bigcup\limits_{C\in\mathcal M'}C$ is connected and
 \item $X\setminus(D\cup\bigcup\limits_{C\in\mathcal M'}C)$ is not connected for all $D\in\mathcal M\setminus\mathcal M'$,
\end{itemize}
is called \begriff{maximal non-separating system} in $\mathcal M$.
\end{definition}

\begin{satz}
\label{maxsys}
Let $X$ and $\mathcal M$ be given as above and let $\mathcal M'\subseteq\mathcal M$ and $\mathcal M''\subseteq\mathcal M$ be two maximal non-separating systems. Then we have $\vert\mathcal M'\vert=\vert\mathcal M''\vert$.
\end{satz}

\begin{Bew}
Let $\Upsilon:=\Upsilon(\mathcal M)$ be the graph defined as follows:
\begin{itemize}
\item $V(\Upsilon):=$ set of components of $X\setminus\bigcup\limits_{C\in\mathcal M}C$,
\item Let $e_1,e_2\in V(\Upsilon)$. For all $L\in \mathcal M$ such that $e_1\cup L\cup e_2$ is connected, add a geometric edge between $e_1$ and $e_2$.
\end{itemize}
Let $n:=\vert V(\Upsilon)\vert$ and $m:=\vert E(\Upsilon)\vert=\vert\mathcal M\vert$.\\
To cut the surface $X$ along a maximal non-separating system in $\mathcal M$ corresponds to delete as many edges from $\Upsilon$ such that the remaining graph $\Upsilon'$ is connected and that $\forall k\in K(\Upsilon')$ the graph $\Upsilon'\setminus k$ is not connected. Therefore, $\Upsilon'$ is also a spanning tree and has $n-1$ edges. So $m-(n-1)$ edges have been deleted and this number is $\vert\mathcal M'\vert=\vert\mathcal M''\vert$.
\end{Bew}

\begin{figure}[h]
\begin{center}
\includegraphics[scale=0.45]{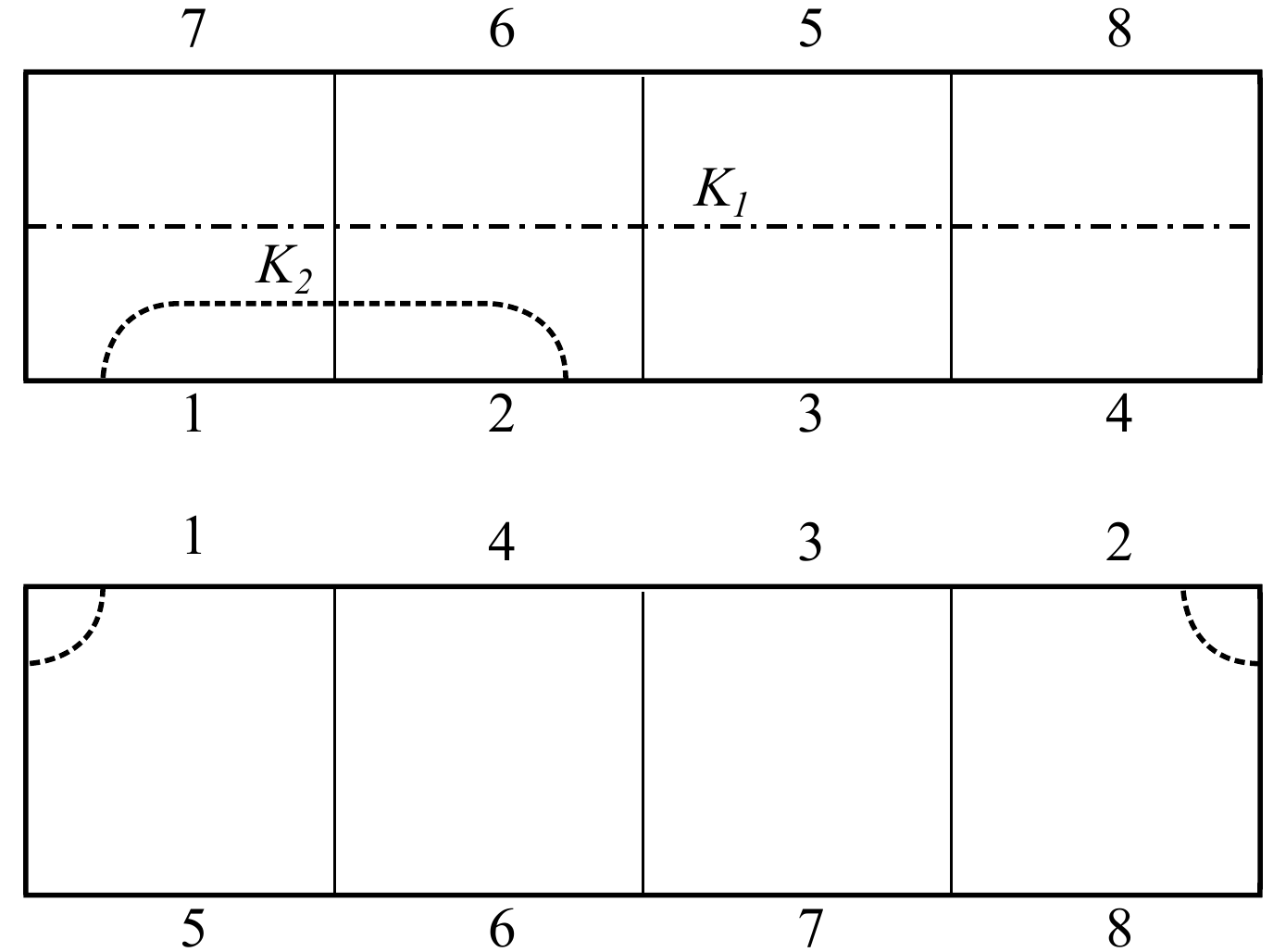}
\end{center}
\caption{Tajectory of curves on the cylinders}
\label{Origami3}
\end{figure}

\begin{figure}[h]
\begin{center}
\includegraphics[scale=0.45]{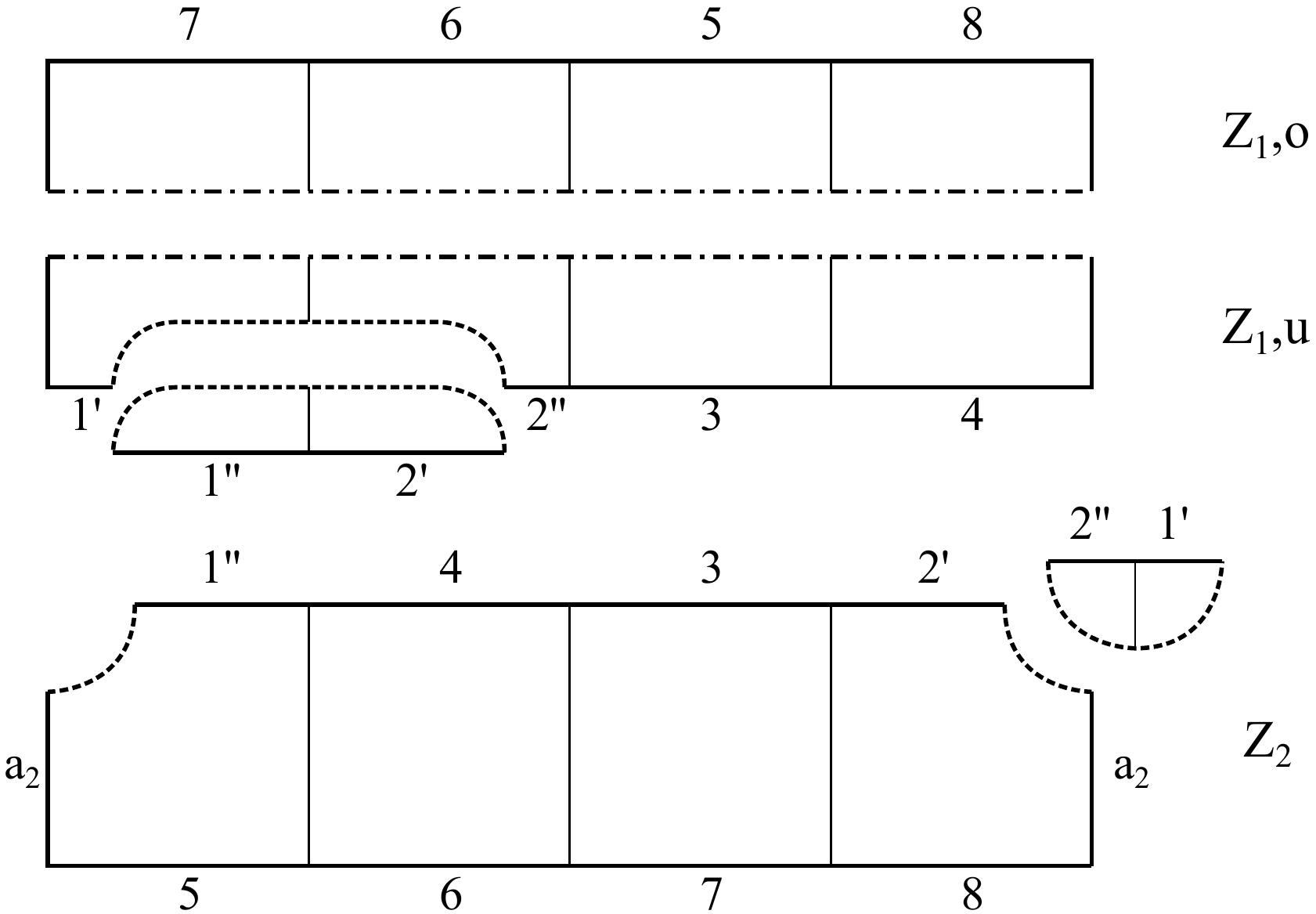}
\end{center}
\caption{Cut origami}
\label{Origamiaufgeschnitten3}
\end{figure}

\begin{figure}[h]
\begin{center}
\includegraphics[scale=0.45]{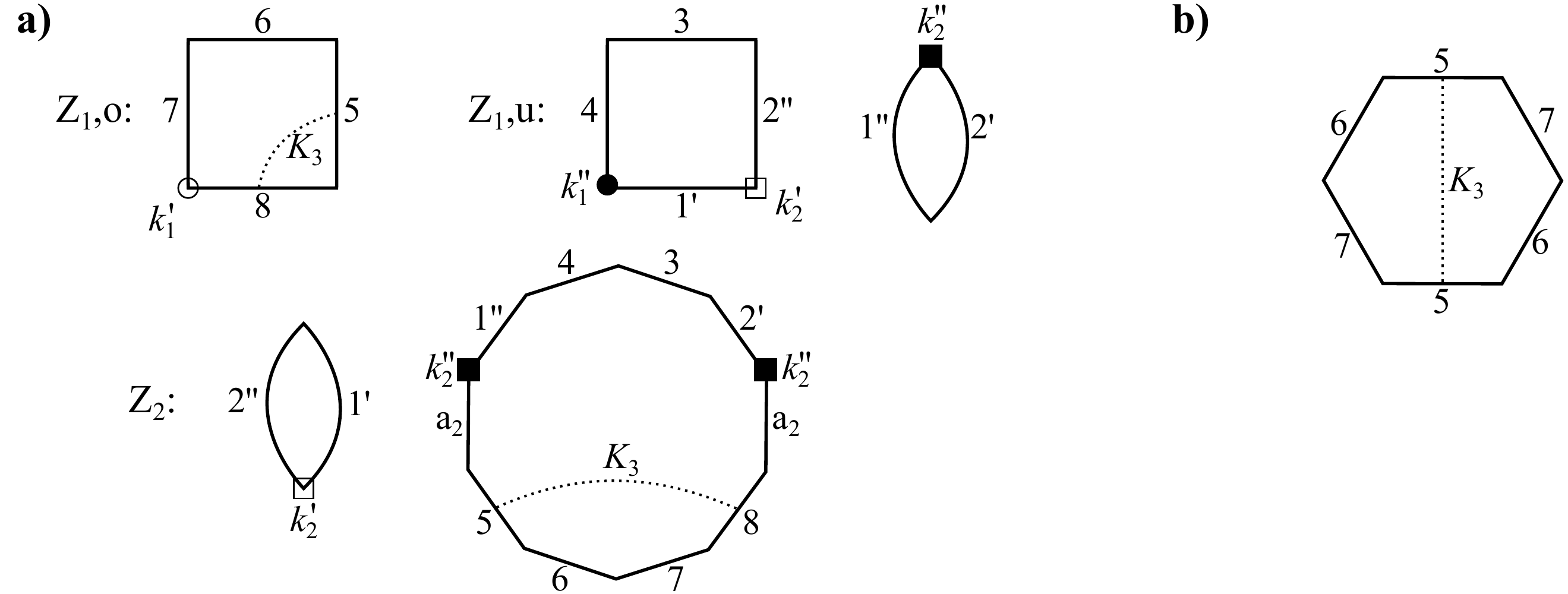}
\caption{The surface $X'$}
\label{Origamiaufgeschnittenzus3}
\end{center}
\end{figure}

\begin{bem}
To require that $\mathcal M$ is finite is not really a restriction because there are at most $3g-3$ pairwise disjoint, non-isotopic curves (pants decomposition). But in a maximal non-separating system $\mathcal M$ two curves cannot be isotopic, and a curve can be replaced by an isotopic one.
\end{bem}
The algorithm to find a HSS that will be formalized below shall be first explained with an example.\\
Take an origami $O=(p:X^*\rightarrow E^*)$. First we cut in the middle as much horizontal cylinders as possible such that the resulting surface is still connected. In other words: We look for a maximal non-separating system in the set of all horizontal cuts through the middle of the cylinders. The number of cylinders that we cut in this step does not depend on the choice of the particular cylinder, by Thm. \ref{maxsys}. In figure \ref{Origami3}, the horizontal cut through $Z_1$, namely $K_1$, is an example for such a maximal system.\\
When we have such a maximal system, we proceed by constructing a new non-separating horizontal curve (as $K_2$ in figure \ref{Origami3}) and repeat this step until we have $g$ such curves.
~\\
This step (that has been done already once in figure \ref{Origami3} to get $K_2$) works as follows: We cut $X$ along the $n$ curves that we already have constructed in the first step (figure \ref{Origamiaufgeschnitten3}) and identify each of the $2n$ boundary curves of the resulting surface to a point (figure \ref{Origamiaufgeschnittenzus3}a, here the points that come from $K_1$ and $K_2$ are denoted as $k_1'$, $k_1''$, $k_2'$ and $k_2''$). We now have a closed surface $X'$ of genus $g-n$. On this surface we find another non-separating curve if we glue the polygons from figure \ref{Origamiaufgeschnittenzus3} to one polygon and take a pair of identified sides that separates another such pair, and connect the centers of these sides with a straight line.\\
If we now follow back this straight line along the step, where we glued the polygons to one, we can draw the new curve in our origami (see figure \ref{neuekurve}). This curve is already horizontal as follows from Lemma \ref{schonhorizontal} below.\\
\\
\begin{figure}[h!]
\begin{center}
\includegraphics[scale=0.45]{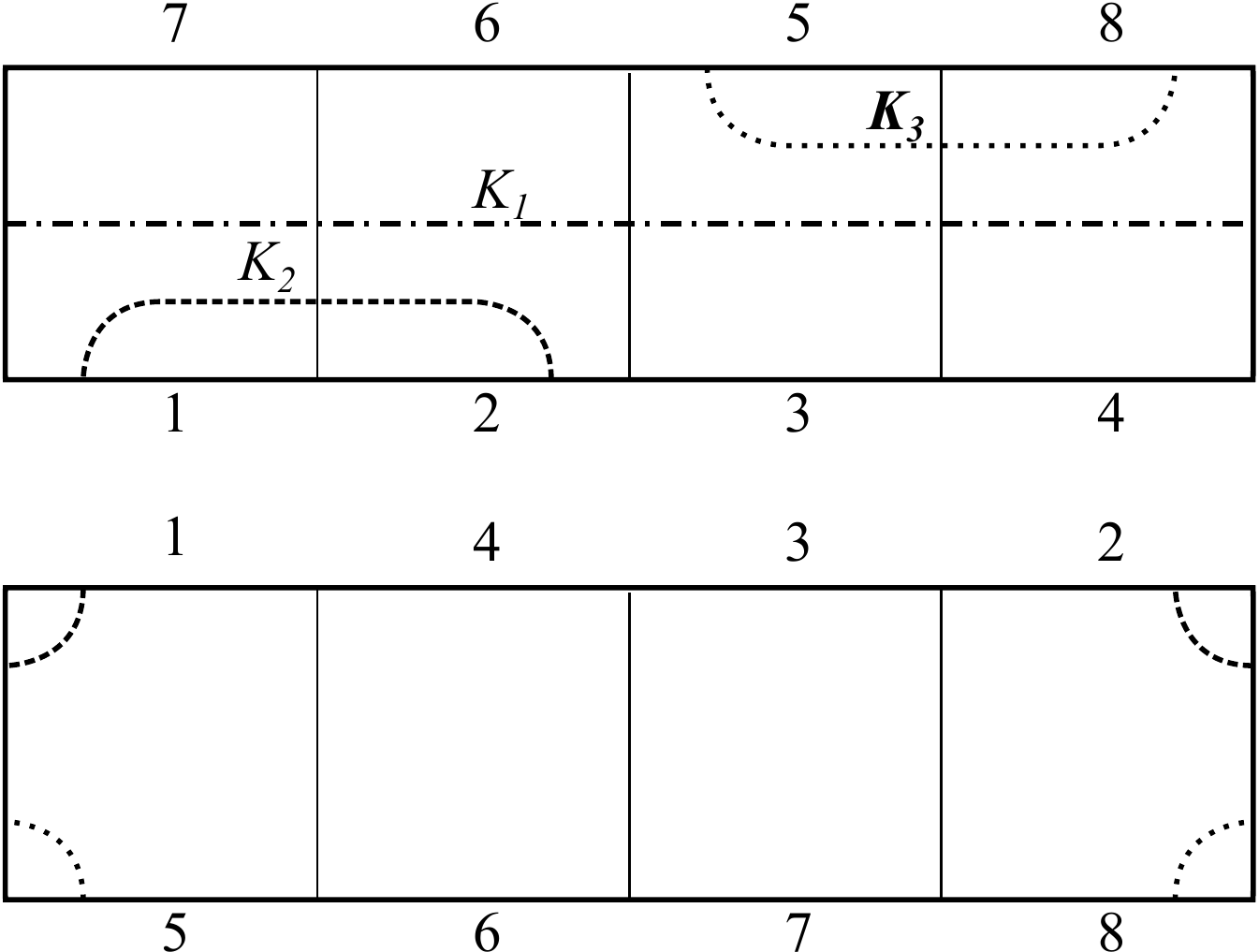}
\caption{The new curve $K_3$ on the surface $X$}
\label{neuekurve}
\end{center}
\end{figure}

\section{The algorithm}
From now on we view an origami $O = (p\colon X^*\rightarrow E^*)$ as set $Q:=\{1,\dots,d\}$ of squares, together with 2 permutations $p_1,p_2\in\mathfrak{S}_d$ such that $\langle p_1,p_2\rangle$ acts transitively on $Q$ and such that $p_1$ and $p_2$ are defined as follows:
\begin{center}
\begin{tabular}[c]{rlcrl}
$p_1:Q$&$\ra Q,$&\quad&$p_2:Q$&$\ra Q$\\
$x$&$\mapsto$ right neighbour of $x$ &&$x$&$\mapsto$ upper neighbour of $x$
\end{tabular}
\end{center}
The horizontal cylinders of $O$ are the orbits of $p_1$, the vertices of $O$ are the orbits of $p_2p_1p_2^{-1}p_1^{-1}$; here we identify a square with its left lower vertex and we use the mapping notation, i.e.\ $p_1p_2(x):=p_1(p_2(x))$. Identifying the squares of $O$ with their lower edge, $Q$ can also be seen as the set of horizontal edges. The genus of $X$ is then $g=\frac{d-\sharp\text{vertices}}{2}+1$.

\paragraph{Step 1: Find a maximal system of horizontal cuts at the middle of cylinders}~\\
First we need a maximal non-separating system in the set of horizontal cuts at the middle of all cylinders.
Therefore, we define a graph $\Upsilon(O)$ as follows:
\begin{itemize}
\item $V(\Upsilon(O))$: Take two vertices $z_i^o$ and $z_i^u$ for every horizontal cylinder $Z_i$.
\item $E(\Upsilon(O))$: Let $z_i^o$ and $z_j^u$ be connected with an edge if $p_2(Z_i)\cap Z_j\neq\emptyset.$
\end{itemize}
The set of components of $\Upsilon(O)$ corresponds to the set of components of the origami after cutting it at the middle of all horizontal cylinders. We add consecutively bridges connecting vertices $z_i^o$ and $z_i^u$ for some $i$ and thereby decrease the number of connected components of the graph by 1. We proceed until the graph is connected.\\
Adding a bridge connecting $z_i^o$ and $z_i^u$ corresponds to glueing the cut surface along the horizontal cut on $Z_i$.
Therefore, horizontal cuts along the cylinders, for which \textit{no} bridge has been added, form a maximal non-separating system in the set of horizontal cuts along all cylinders. If the number of these cuts is less than $g$, proceed with step 2.

\paragraph{Step 2: Add a new horizontal curve to a system of non-separating horizontal cuts at cylinders}~\\
First, we want to represent the polygons that add up to $X'$ (see figure \ref{Origamiaufgeschnittenzus3}a) in a proper way by cyclic lists (i.e.\ lists, where the order of the entries is only given up to cyclic permutation) that represent the order of the edges in counterclockwise direction.\\
Therefore, we define for every cylinder $Z=\{p_1(x),\dots,p_1^n(x)=x\}$ two cyclic lists
\[L(Z,u):=[p_1(x),\dots,p_1^n(x)]\;\text{ and }\;L(Z,o):=[p_2p_1^n(x),\dots,p_2p_1(x)].\]
These correspond to the lower edges of the cylinder from left to right respectively to the upper edges from right to left. If $Z$ has been cut in the first step, these lists correspond to two polygons of $X'$; if $Z$ has not been cut, let $a_Z$ be some further entry in the list, and we define additionally the cyclic list
\[L(Z):=[a_Z,L(Z,u),a_Z,L(Z,o)]=[a_Z,p_1(x),\dots,p_1^n(x), a_Z, p_2p_1^n(x),\dots,p_2p_1(x)].\]
$L(Z)$ corresponds to one of the polygons of $X'$, where $a_Z$ belongs to a vertical edge   (see cylinder $Z_2$ in figure \ref{Origamiaufgeschnitten3} respectively the corresponding polygon in figure \ref{Origamiaufgeschnittenzus3}a).\\
Now, let  $\mathcal Z:=\{Z_1,\dots,Z_b\}$ be the set of horizontal cylinders of $O$, $\overline{\mathcal Z}:=\{Z_1,\dots,Z_a\}\subseteq\mathcal Z$ the subset of the cylinders that have been cut in step 1 and $\mathcal Z':=\{Z_{a+1},\dots,Z_b\}\subseteq\mathcal Z$ the subset of the cylinders that have not been cut in step 1. We define three sets of lists:
$$\mathcal L_u:=\{L(Z_1,u),\dots,L(Z_b,u)\}\text{, }\;\;\;\;\mathcal L_o:=\{L(Z_1,o),\dots,L(Z_b,o)\}$$
and
$$\mathcal L:=\{L(Z_1,u),\dots,L(Z_a,u),L(Z_1,o),\dots,L(Z_a,o),L(Z_{a+1}),\dots,L(Z_b)\}$$
$\mathcal L$ corresponds to a representation of $X'$ by several polygons to be glued, as in figure \ref{Origamiaufgeschnittenzus3}a. Moreover, we define the following partitions of $\mathcal L_u$ and $\mathcal L_o$:
$$\mathcal L_u=\mathcal L_{Z_1,u}\cup\dots\cup\mathcal L_{Z_b,u}\text{ with } \mathcal L_{Z_i,u} := \{L(Z_i,u)\}\text{ for } i\in\{1,\dots,b\}$$
$$\mathcal L_o=\mathcal L_{Z_1,o}\cup\dots\cup\mathcal L_{Z_b,o}\text{ with } \mathcal L_{Z_i,o} := \{L(Z_i,o)\}\text{ for } i\in\{1,\dots,b\}$$
In this step the sets $\mathcal L_{Z_i,u}$ and $\mathcal L_{Z_i,o}$ still have only one element; in the next steps we shall collect there what we cut out of the upper respectively lower half of the corresponding cylinder (see figure \ref{Origamiaufgeschnitten3}). In general $\mathcal L_u$ and $\mathcal L_o$ are \textit{not} subsets of $\mathcal L$ because we account also for those cylinders that in step 1 were not cut in the middle what belongs to the ``upper half'', i.e.\ lies in $\mathcal L_o$, and what belongs to the ``lower half'', i.e.\ belongs to $\mathcal L_u$.\\
Just as we glued the polygons in figure \ref{Origamiaufgeschnittenzus3} from a) to b), we now merge the lists in $\mathcal L$. This happens step by step (for example by always concatenating the first list with the next possible) as follows:\\
Two lists $L$ and $M$ with common entry $x$ are concatenated at $x$, by replacing $x$ in $L$ by the list $M\setminus x$, beginning with the entry after $x$ and proceeding cyclically until the entry before $x$. More precisely: If
$$L = [a_1,\dots,a_n,x,b_1,\dots,b_m]\;\text{ and }\;M = [c_1,\dots,c_p,x,d_1,\dots,d_q]$$
shall be concatenated at $x$, we get the list
\[L' = [a_1,\dots,a_n,d_1,\dots,d_q,c_1,\dots,c_p,b_1,\dots,b_m]\]
In addition we eliminate two identical elements if they are consecutive in a cyclic list.
If everything is concatenated to a list $P$, we choose in it a pair $(\alpha,\alpha)$ of entries that separates another pair $(\beta,\beta)$:
\[P=[\alpha,\dots,\beta,\dots,\alpha,\dots,\beta,\dots]\]
In the polygon corresponding to $P$, the pair $(\alpha,\alpha)$ corresponds to the closed curve $K$ that intersects the edge $\alpha$. Since  $(\alpha,\alpha)$ separates the pair $(\beta,\beta)$, $K$ is a separating curve. It is clear that such a pair $(\alpha,\alpha)$ separating another pair $(\beta,\beta)$ exists. We can see this as follows: Assume that such a pair does not exist. Let $(\alpha,\alpha)$ be any pair in $P$. Such a pair exists because the surface $X'$ has genus $g\geq1$ and so $P$ is not empty. This pair cannot be a pair of consecutive entries since such pairs have already been eliminated. Therefore, there is an element $\beta$ in between. By assumption, the other $\beta$ lies on the same side of the pair $(\alpha,\alpha)$. Between the pair $(\beta,\beta)$ $-$ on the other side as the pair $(\alpha,\alpha)$ $-$ lies again a pair $(\gamma,\gamma)$. So we get a finite sequence of pairs that come closer and closer together, until we have a pair of consecutive entries, a contradiction.\\
To reconstruct $K$, as it lies in the   origami, we look at the pair $(\alpha,\alpha)$ and follow it back through the steps, where we concatenated the lists from $\mathcal L$. If we go back a step of concatenating two lists along $x$, where the pair is separated again, we replace $(\alpha,\alpha)$ by the pairs $(\alpha,x)$ and $(x,\alpha)$, and we proceed like this if these pairs are separated and so on. Meanwhile, we always remember the lists that contain the pairs. At the end we have a chain of pairs
\begin{equation}
\label{kette1}
(\alpha_0,\alpha_1), (\alpha_1,\alpha_2), \dots, (\alpha_{r-2},\alpha_{r-1}), (\alpha_{r-1},\alpha_r)
\end{equation}
that correspond in this order to the pieces that build up $K$.\\
More precisely: A pair $(\alpha_s,\alpha_{s+1})$ is contained in a list from $L'\in\mathcal L$. Two cases are now possible:
\begin{enumerate}
\item[1)] $L'=L(Z,o)\in\mathcal L$ or $L'=L(Z,u)\in\mathcal L$ for a $Z\in\overline{\mathcal Z}$.
\item[2)] $L'=L(Z)=[a_Z,L(Z,u),a_Z,L(Z,o)]$ for a $Z\in\mathcal Z'$.
\end{enumerate}
In case 1 we define $L:=L'$. In case 2 $(\alpha_s,\alpha_{s+1})$ is either completely contained in the part $L(Z,u)$ or completely contained in the part $L(Z,o)$, as will follow from Lemma \ref{schonhorizontal}. We define $L:=L(Z,u)$ or $L:=L(Z,o)$ depending on where $(\alpha_s,\alpha_{s+1})$ lies.\\
We now want to represent the new curve $K$ by a starting point $q$ and a word $w$ in $x$ and $y$ that describes the horizontal and vertical steps of the curve. We choose $q:=\alpha_0$ if $\alpha_0\in L\in\mathcal L_u$ and $q:=p_2^{-1}(\alpha_0)$ if $\alpha_0\in L\in\mathcal L_o$.
\begin{itemize}
\item If $L\in \mathcal L_u$, choose a $t_s\in\Z$ with minimal absolute value such that $p_1^{t_s}(\alpha_s)=\alpha_{s+1}$.\\
We then define $w_s:=x^{t_s}y^{-1}$.
\item If $L\in \mathcal L_o$, choose a $t_s\in\Z$ with minimal absolute value such that $p_1^{t_s}p_2^{-1}(\alpha_s)=p_2^{-1}(\alpha_{s+1})$.\\ We then define $w_s:=x^{t_s}y$.
\end{itemize}
The new curve in the Schottky cut system is:
$$K:=(q,w), \text{ with } w:= w_0\cdots w_{r-1}$$
$K$ is horizontal, see Lemma \ref{schonhorizontal} below. If we still have less than $g$ curves, proceed with step 3.
\paragraph{Step 3: Add a new horizontal curve to a system of non-separating horizontal curves}~\\
In this step, it is possible that we cut out pieces from the upper or lower halves of the cylinders when we cut the origami along the curves already given. The polygons we get do not come all either from half or whole cylinders, but also from pieces that have been cut out from the latter, as for example in figure \ref{Origamiaufgeschnittenzus3}a, the polygon with the edges $1''$ and $2'$ has been cut out from $(Z_1,u)$ or the polygon with the edges $1'$ and $2''$ has been cut out from $Z_2$. If then the new curve that we'll construct in this step enters a polygon of this kind, it can proceed in the origami only in one direction to the point where it exits the horizontal cylinder because otherwise that curve would intersect another curve. Therefore, we cannot write these pieces as \textit{cyclic} lists. The other polygons, that come from the parts left after cutting out these pieces, can be written as cyclic lists as before (for example, in figure \ref{Origamiaufgeschnittenzus3}a all the other polygons).\\
In other words: If we cut an origami along the horizontal sides of the squares and around the middles of (some) cylinders, we get twice connected bounded surfaces that correspond to the horizontal cylinders or half cylinders. If we then cut out of them pieces along arcs, we get simply connected pieces \linebreak($\ra$ non-cyclic lists) and twice connected pieces ($\ra$ cyclic lists). That's the reason why we have to choose another procedure as in step 2, where we constructed the first horizontal curve that is not a curve around a cylinder. Here we'll have to proceed as follows:\\
Take the lists from the set $\mathcal L$ in the last step and change them using the chain of pairs
\[(\alpha_0,\alpha_1), (\alpha_1,\alpha_2), \dots, (\alpha_{r-2},\alpha_{r-1}), (\alpha_{r-1},\alpha_r).\]
Every list that contains a pair $(\alpha_s,\alpha_{s+1})$ shall be changed according to the following rules: Let $(\alpha_s,\alpha_{s+1})\in L$:
\begin{itemize}
\item If $L\in\mathcal L_{Z,u}$ and $t(s)>0$, replace in $\mathcal L_{Z,u}$ and in $\mathcal L$ the list
$$L=[a_1,\dots,a_n,\alpha_s,b_1,\dots,b_m,\alpha_{s+1},c_1,\dots,c_p]$$
by the lists
$$L_0:=[a_1,\dots,a_n,\alpha'_s,\alpha''_{s+1},c_1,\dots,c_p]\;\text{ and }\;L_1:=[\alpha''_s,b_1,\dots,b_m,\alpha'_{s+1}],$$
where $L_1$ shall be non-cyclic and $L_0$ shall be cyclic if and only if $L$ is cyclic.\\
If $t(s)<0$, do the same with the roles of $\alpha_s$ and $\alpha_{s+1}$ interchanged.
\item If $L\in\mathcal L_{Z,o}$ and $t(s)<0$, replace in $\mathcal L_{Z,o}$ and in $\mathcal L$ the list
$$L=[a_1,\dots,a_n,\alpha_s,b_1,\dots,b_m,\alpha_{s+1},c_1,\dots,c_p]$$
by the lists
$$L_0:=[a_1,\dots,a_n,\alpha''_s,\alpha'_{s+1},c_1,\dots,c_p]\;\text{ and }\;L_1:=[\alpha'_s,b_1,\dots,b_m,\alpha''_{s+1}],$$
where $L_1$ shall be non-cyclic and $L_0$ shall be cyclic if and only if $L$ is cyclic.\\
If $t(s)>0$, do the same with the roles of $\alpha_s$ and $\alpha_{s+1}$ interchanged.
\end{itemize}
If $L=L(Z)=[a_Z,L(Z,u),a_Z,L(Z,o)]\in\mathcal L\setminus(\mathcal L_u\cup\mathcal L_o)$, the pair $(\alpha_s,\alpha_{s+1})$ lies in $L(Z,u)$ or in $L(Z,o)$ (see Lemma \ref{schonhorizontal}): Change the list $L(Z,u)$ respectively $L(Z,o)$ as above and redefine in $\mathcal L$ the list $L(Z)$ by $L(Z):=[a_Z,L_0,a_Z,L(Z,o)]$ or $L(Z):=[a_Z,L(Z,u),a_Z,L_0]$, depending on if the pair $(\alpha_s,\alpha_{s+1})$ was before in $L(Z,u)$ or in $L(Z,o)$, and add the list $L_1$ to the sets $\mathcal L$ and $\mathcal L_{Z,u}\subseteq\mathcal L_u$ respectively $\mathcal L_{Z,o}\subseteq\mathcal L_o$.\\
As in step 2, we concatenate the lists from $\mathcal L$ to one list and get a chain of pairs
\begin{equation}
\label{kette2}
(\alpha_0,\alpha_1), (\alpha_1,\alpha_2), \dots, (\alpha_{r-2},\alpha_{r-1}), (\alpha_{r-1},\alpha_r).
\end{equation}
As starting point $q$ for the new curve $K$ we choose $\alpha_0$ if $\alpha_0\in L\in\mathcal L_u$, and we choose $p_2^{-1}(\alpha_0)$ if $\alpha_0\in L\in\mathcal L_o$. If the pair $(\alpha_s,\alpha_{s+1})$ lies in a cyclic list, we get the pieces $w_s$ of the new curve as in step 2. If $(\alpha_s,\alpha_{s+1})$ lies not in a cyclic list, we proceed as follows:
\begin{itemize}
\item If $L\in \mathcal L_{Z,u}$ and $\alpha_s$ appears before $\alpha_{s+1}$, let $t(s)\in\N$ be minimal such that $p_1^{t(s)}(\alpha_s)=\alpha_{s+1}$, and let the corresponding piece of the curve be $w_s:=x^{t(s)}y^{-1}$.
\item If $L\in \mathcal L_{Z,u}$ and $\alpha_{s+1}$ appears before $\alpha_s$, let $t(s)\in\N$ be minimal such that $p_1^{t(s)}(\alpha_{s+1})=\alpha_s$, and let the corresponding piece of the curve be $w_s:=x^{-t(s)}y^{-1}$.
\item If $L\in \mathcal L_{Z,o}$ and $\alpha_s$ appears before $\alpha_{s+1}$, let $t(s)\in\N$ be minimal such that \linebreak$p_1^{t(s)}p_2^{-1}(\alpha_{s+1})=p_2^{-1}(\alpha_s)$, and let the corresponding piece of the curve be $w_s:=x^{-t(s)}y$.
\item If $L\in \mathcal L_{Z,o}$ and $\alpha_{s+1}$ appears before $\alpha_s$, let $t(s)\in\N$ be minimal such that \linebreak$p_1^{t(s)}p_2^{-1}(\alpha_s)=p_2^{-1}(\alpha_{s+1})$, and let the corresponding piece of the curve be $w_s:=x^{t(s)}y$.
\end{itemize}
The new curve in the Schottky cut system is
$$K:=(q,w), \text{ with } w:= w_0\cdots w_{r-1}.$$
It follows again from the following lemma that $K$ is horizontal.\\
We repeat step 3 until we have $g$ curves.

\begin{lemma}
\label{schonhorizontal}
The curves that come from the steps 2 and 3 are horizontal.
\end{lemma}

\begin{Bew}
Let $K$ be such a curve. Assume that $K$ is not horizontal. Then there is a cylinder $Z$ crossed vertically by $K$ (i.e.\ such that $K$ enters $Z$ from the upper side and exits from the lower side or the other way round). But $Z$ cannot be crossed this way only once by $K$ since we assumed that we already have cut the surface $X$ along a maximal non-separating system in the set of the horizontal curves around the cylinders and that the horizontal cut around $Z$ is not in this system. That means that we would separate the surface if, in addition to our maximal non-separating system, we cut $X$ also along this curve around $Z$. But if $K$ would cross vertically the cylinder $Z$  only once, there would be a path from one side of $Z$ to the other, and, therefore, the surface would be still path-connected, a contradiction. Therefore, $K$ crosses $Z$ vertically more than once.\\
But $K$ was constructed as follows: Let $P$ be the list that we get if we concatenate all lists in $\mathcal L$ to one. There we find, as we explained above, a pair $(\alpha,\alpha)$ that separates another pair. If we now follow back these steps of concatenating the lists in our algorithm, we get a chain of pairs as in (\ref{kette1}) and (\ref{kette2}), corresponding to the new curve. This curve cuts the polygons belonging to the lists from $\mathcal L$ at most once since the pair
$(\alpha,\alpha)$ corresponds to the arc between the two sides $\alpha$ in the polygon corresponding to $P$. This arc splits in two parts if and only if $P$ is divided in two polygons in a way that the two $\alpha$ get separated and we have then two lists, each with one such pair corresponding to an arc. The same applies to pairs of the form $(\alpha_s,\alpha_{s+\nu})$. Therefore, this algorithm gives no possibilities that a new piece of a curve appears in a polygon.
But there is only one list in $\mathcal L$ that concatenates the upper with the lower side of $Z$, namely $L(Z)$, and the curve $K$ crosses the cylinder $Z$ only once, a contradiction.
\end{Bew}

\begin{figure}[!h]
\begin{center}
\includegraphics[scale=0.5]{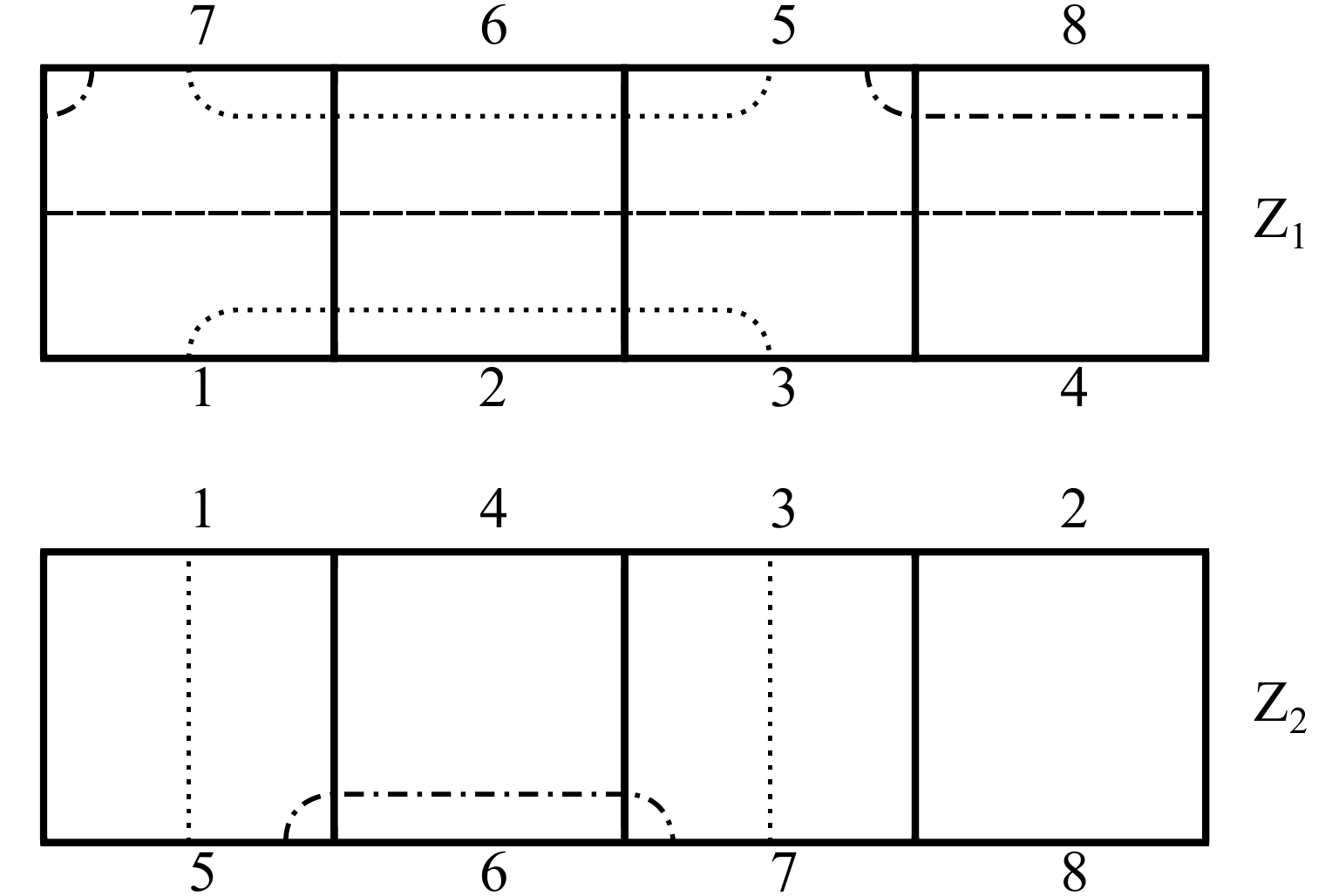}
\end{center}
\caption{A non-horizontal Schottky cut system}
\label{wollmilchsaunichthor}
\end{figure}

\begin{bem}
The last lemma does not imply that we cannot add a non-horizontal curve to a horizontal cut system to get a Schottky cut system, as figure \ref{wollmilchsaunichthor} shows.
\end{bem}

\section{Two examples for the algorithm}
\subsection*{The Wollmilchsau}
The Wollmilchsau (described in \cite{woll}) is the origami from figure \ref{Origami3}-\ref{neuekurve}. It is given by $Q=\{1,\dots,8\}$,
\[p_1=
\left(\begin{array}{cccccccc}
1&2&3&4&5&6&7&8\\
2&3&4&1&6&7&8&5\\
\end{array}\right)\quad\text{ and }\quad
p_2=\left(\begin{array}{cccccccc}
1&2&3&4&5&6&7&8\\
7&6&5&8&1&4&3&2\\
\end{array}\right)\]
\begin{itemize}
\item \textbf{cylinder:} $Z_1 = \{1,2,3,4\}; Z_2 = \{5,6,7,8\}$
\item \textbf{vertices:} orbits of $p_2p_1p_2^{-1}p_1^{-1}:\{1,3\},\{2,4\},\{5,7\},\{6,8\}$
\item \textbf{genus:} $g=\frac{d-\sharp\text{vertices}}{2}+1=\frac{8-4}{2}+1=3$
\end{itemize}
\paragraph{Step 1:}~\\
Build a graph:
\[
\begin{xy}
\xymatrix{
z_1^o\bullet\ar@{-}[r]&\bullet z_2^u\ar@<-5pt>@{.}[d]\\
z_1^u\bullet\ar@{-}[r]&\bullet z_2^o
}
\end{xy}
\]
We connect $z_2^o$ with $z_2^u$, to get a connected  graph.\\
\paragraph{Step 2:}~\\
Set of cylinders which have been cut horizontally: $\overline{\mathcal Z}=\{Z_1\}$\\
Set of cylinders which have not been cut horizontally: $\mathcal Z'=\{Z_2\}$\\
Let $K_1$ be the curve around $Z_1$.
\begin{center}
\begin{tabular}[c]{lll}
Lists in $\mathcal L_u$:&\;\;\;&Lists in $\mathcal L_o$:\\
$L(Z_1, u) = [1,2,3,4]$&&$L(Z_1, o) = [8,5,6,7]$\\
$L(Z_2, u) = [5,6,7,8]$&&$L(Z_2, o) = [2,3,4,1]$\\
&\;&
\end{tabular}\\
Lists in $\mathcal L\setminus(\mathcal L_u\cup\mathcal L_o)$:\\
$L(Z_2) = [a_2,5,6,7,8,a_2,2,3,4,1]$
\end{center}
We now have:
$$\mathcal L=\{L(Z_1, u),L(Z_1, o),L(Z_2)\}$$
and concatenate all to one list (always at the entries in bold) as in the left column of the table below. When all lists are concatenated to one, we take a pair that separates another (for example the pair (1,1) that separates the pair (3,3)) and follow it back through the steps of concatenating the lists (see right column of the table below).
\begin{center}
\begin{tabular}[c]{cc|l|l|lc}
\cline{3-4}
\multirow{7}{*}{\begin{sideways}concatenate lists\end{sideways}}&\multirow{7}{*}{\begin{sideways}$\longleftarrow$\end{sideways}}&\footnotesize$L(Z_1, u) = [1,\textbf{2},3,4]$&\footnotesize$(1,2)(2,1)$&\multirow{7}{*}{\begin{sideways}$\longrightarrow$\end{sideways}}&\multirow{7}{*}{\begin{sideways}new curve \end{sideways}}\\
&&\footnotesize$L(Z_1, o) = [8,5,6,7]$&&&\\
&&\footnotesize$L(Z_2) = [a_2,5,6,7,8,a_2,\textbf{2},3,4,1]$&&&\\
\cline{3-4}
&&\footnotesize$[1,3,4,1,a_2,5,6,7,\textbf{8},a_2 ,3,4]$&\footnotesize$(1,1)$&&\\
&&\footnotesize$[\textbf{8},5,6,7]$&&&\\
\cline{3-4}
&&\footnotesize$[1,3,4,1,a_2,5,6,7,5,6,7,a_2 ,3,4]$&\footnotesize$(1,1)$&&\\
\cline{3-4}
\end{tabular}
\end{center}
$(1,2)$ is in $L(Z_1,u)$. The number $t_0$ with minimal absolute value such that $p_1^{t_0}(1)=2$ is $t_0=1$. So we have $w_0=xy^{-1}$.\\
$(2,1)$ is in $L(Z_2)$ in the part $L(Z_2,o)$. The number $t_1$ with minimal absolute value such that $p_1^{t_1}p_2^{-1}(2)=p_2^{-1}(1)$ is $t_1=1$. So we have $w_1=xy$.\\
The new curve is, hence,
\[K_2:= (1,w_0w_1) = (1,xy^{-1}xy).\]
This corresponds until now to the situation as in figure \ref{Origami3}.\\
\paragraph{Step 3:}~\\
The new sets of lists are now (cyclic lists in square brackets, non-cyclic list in round brackets):
\begin{center}
\begin{tabular}[c]{lll}
Lists in $\mathcal L_u$&\;\;\;&Lists in $\mathcal L_o$\\
$L(Z_1, u) = [1',2'',3,4]$&&$L(Z_1, o) = [8,5,6,7]$\\
$L_1(Z_1, u) = (1'',2')$&&$L(Z_2, o) = [2',3,4,1'']$\\
$L(Z_2, u) = [5,6,7,8]$&&$L_1(Z_2, o) = (1',2'')$\\
\;&&
\end{tabular}\\
Lists in $\mathcal L\setminus(\mathcal L_u\cup\mathcal L_o)$:\\
$L(Z_2) = [a_2,5,6,7,8,a_2,2',3,4,1'']$\\
\end{center}
So we have:
$$\mathcal L=\{L(Z_1, u),L_1(Z_1, u),L(Z_1, o),L_1(Z_2, o),L(Z_2)\}$$
and concatenate everything to one list:
\begin{center}
\begin{tabular}[c]{|l|l|}
\hline
\footnotesize$L(Z_1,u) = [1',\textbf{2}'',3,4]$&\footnotesize$(5,8)(8,5)$\\
\footnotesize$L_1(Z_1,u) = (1'',2')$&\\
\footnotesize$L(Z_1,o) = [8,5,6,7]$&\\
\footnotesize$L_1(Z_2,o) = (1',\textbf{2}'')$&\\
\footnotesize$L(Z_2) = [a_2,5,6,7,8,a_2,2',3,4,1'']$&\\
\hline
\footnotesize$[\textbf{3},4]$&\footnotesize$(5,8)(8,5)$\\
\footnotesize$(1'',2')$&\\
\footnotesize$[8,5,6,7]$&\\
\footnotesize$[a_2,5,6,7,8,a_2,2',\textbf{3},4,1'']$&\\
\hline
\footnotesize$[\textbf{1}'',a_2,5,6,7,8,a_2,2']$&\footnotesize$(5,8)(8,5)$\\
\footnotesize$(\textbf{1}'',2')$&\\
\footnotesize$[8,5,6,7]$&\\
\hline
\footnotesize$[5,6,7,\textbf{8}]$&\footnotesize$(5,8)(8,5)$\\
\footnotesize$[\textbf{8},5,6,7]$&\\
\hline
\footnotesize$[5,6,7,5,6,7]$&\footnotesize$(5,5)$\\
\hline
\end{tabular}
\end{center}
$(5,8)$ is in $L(Z_2)$ in the part $L(Z_2,u)$. The number $t_0$ with minimal absolute value such that $p_1^{t_0}(5)=8$, is $t_0=-1$.
So we have $w_0=x^{-1}y^{-1}$.\\
$(8,5)$ is in $L(Z_1,o)$. The number $t_1$ with minimal absolute value such that $p_1^{t_1}p_2^{-1}(8)=p_2^{-1}(5)$, is $t_1=-1$. So we have $w_1=x^{-1}y$.\\
The new curve is
\[K_3:= (5,w_0w_1) = (5,x^{-1}y^{-1}x^{-1}y).\]
\subsection*{Another example}
Take the origami defined by $Q=\{1,\dots,14\}$,
\[p_1=
\left(\begin{array}{cccccccccccccc}
1&2&3&4&5&6&7&8&9&10&11&12&13&14\\
2&3&4&1&6&7&5&9&10&11&12&13&14&8
\end{array}\right)\]
and
\[p_2=\left(\begin{array}{cccccccccccccc}
1&2&3&4&5&6&7&8&9&10&11&12&13&14\\
13&14&11&9&8&10&12&5&3&4&1&2&6&7
\end{array}\right)\]
We'll denote it by $O_{14}$.
\begin{itemize}
\item \textbf{cylinder:} $Z_1 = \{1,2,3,4\}; Z_2 = \{5,6,7\}; Z_3 = \{8,9,10,11,12,13,14\}$
\item \textbf{vertices:} orbits of $p_2p_1p_2^{-1}p_1^{-1}\colon \{1\},\{2\},\{3,6\},\{4\},\{5\},\{7\},\{8,11,12,9,10,13\},\{14\}$
\item \textbf{genus:} $g=\frac{d-\sharp\text{vertices}}{2}+1=\frac{14-8}{2}+1=4$.
\end{itemize}
\paragraph{Step 1:}~\\
Build a graph:
\[
\begin{xy}
\xymatrix{
\underset{\bullet}{z_1^o}\ar@{-}[d]&\underset{\bullet}{z_2^o}\ar@{-}[dl]&\underset{\bullet}{z_3^o}\ar@{-}[dl]\ar@{-}[d]\ar@{.}[dll]\\
\stackrel{\bullet}{z_3^u}&\stackrel{\bullet}{z_1^u}&\stackrel{\bullet}{z_2^u}
}
\end{xy}
\]
We connect $z_3^o$ with $z_3^u$ to get a connected  graph.\\
\paragraph{Step 2:}~\\
Set of cylinders which have been cut horizontally: $\overline{\mathcal Z}=\{Z_1,Z_2\}$\\
Set of cylinders which have not been cut horizontally: $\mathcal Z'=\{Z_3\}$\\
Define the curves $K_1$ and $K_2$ as the horizontal cuts around $Z_1$ and $Z_2$.

\begin{center}
\begin{tabular}[c]{lll}
Lists in $\mathcal L_u$:&\;\;\;&Lists in $\mathcal L_o$:\\
$L(Z_1, u) = [1,2,3,4]$&&$L(Z_1, o) = [9,11,14,13]$\\
$L(Z_2, u) = [5,6,7]$&&$L(Z_2, o) = [12,10,8]$\\
$L(Z_3, u) = [8,9,10,11,12,13,14]$&&$L(Z_3, o) = [7,6,2,1,4,3,5]$\\
&\;&
\end{tabular}\\
Lists in $\mathcal L\setminus(\mathcal L_u\cup\mathcal L_o)$:\\
$L(Z_3) = [a_3,8,9,10,11,12,13,14 ,a_3,7,6,2,1,4,3,5]$
\end{center}
So we have
$$\mathcal L=\{L(Z_1, u),L(Z_2, u),L(Z_1, o),L(Z_2, o),L(Z_3)\}$$
and concatenate everything to one list (always on the entries in bold):
\begin{center}
\begin{tabular}[c]{|l|l|}
\hline
\footnotesize$L(Z_1, u) = [1,\textbf{2},3,4]$&\footnotesize$(8,12)(12,8)$\\
\footnotesize$L(Z_2, u) = [5,6,7]$&\\
\footnotesize$L(Z_1, o) = [9,11,14,13]$&\\
\footnotesize$L(Z_2, o) = [12,10,8]$&\\
\footnotesize$L(Z_3) = [a_3,8,9,10,11,12,13,14 ,a_3,7,6,\textbf{2},1,4,3,5]$&\\
\hline
\footnotesize$[\textbf{5},a_3,8,9,10,11,12,13,14 ,a_3,7,6]$&\footnotesize$(8,12)(12,8)$\\
\footnotesize$[\textbf{5},6,7]$&\\
\footnotesize$[9,11,14,13]$&\\
\footnotesize$[12,10,8]$&\\
\hline
\footnotesize$[8,\textbf{9},10,11,12,13,14]$&\footnotesize$(8,12)(12,8)$\\
\footnotesize$[\textbf{9},11,14,13]$&\\
\footnotesize$[12,10,8]$&\\
\hline
\footnotesize$[8,11,14,13,10,11,\textbf{12},13,14]$&\footnotesize$(8,12)(12,8)$\\
\footnotesize$[\textbf{12},10,8]$&\\
\hline
\footnotesize$[8,11,14,13,10,11,10,8,13,14]$&\footnotesize$(8,8)$\\
\hline
\end{tabular}
\end{center}
$(8,12)$ is in $L(Z_3)$ in the part $L(Z_3,u)$. The number $t_0$ with minimal absolute value such that $p_1^{t_0}(8)=12$ is $t_0=-3$. So we have $w_0=x^{-3}y^{-1}$.\\
$(12,8)$ is in $L(Z_2,o)$. The number $t_1$ with minimal absolute value such that $p_1^{t_1}p_2^{-1}(12)=p_2^{-1}(8)$ is $t_1=1$. So we have $w_1=xy$.\\
The new curve is
$$K_3:= (8, w_0w_1) = (8,x^{-3}y^{-1}xy).$$
\paragraph{Step 3:}~\\
The new sets of lists are now (cyclic lists in square brackets, non-cyclic in round brackets):
\begin{center}
\begin{tabular}[c]{lll}
Lists in $\mathcal L_u$&\;\;\;&Lists in $\mathcal L_o$\\
$L(Z_1, u) = [1,2,3,4]$&&$L(Z_1, o) = [9,11,14,13]$\\
$L(Z_2, u) = [5,6,7]$&&$L_0(Z_2, o) = [12',10,8'']$\\
$L_0(Z_3, u) = [8'',9,10,11,12']$&&$L_1(Z_2, o) = (8',12'')$\\
$L_1(Z_3, u) = (12'',13,14,8')$&&$L(Z_3, o) = [7,6,2,1,4,3,5]$\\
\;&&
\end{tabular}\\
Lists in $\mathcal L\setminus(\mathcal L_u\cup\mathcal L_o)$:\\
$L(Z_3) = [a_3,8'',9,10,11,12',a_3,7,6,2,1,4,3,5]$\\
\end{center}
So we have
$$\mathcal L=\{L(Z_1, u),L(Z_2, u),L_1(Z_3, u),L(Z_1, o),L_0(Z_2, o),L_1(Z_2,o),L(Z_3)\}$$
and concatenate all to one list:
\begin{center}
\begin{tabular}[c]{|l|l||l|l|}
\hline
\footnotesize$L(Z_1, u) = [1,\textbf{2},3,4]$&\footnotesize$(8'',10)(10,8'')$&\footnotesize$[8'',11,8',12'',\textbf{10},11,12']$&\footnotesize$(8'',10)(10,8'')$\\
\footnotesize$L(Z_2, u) = [5,6,7]$&&\footnotesize$[12',\textbf{10},8'']$&\\
\footnotesize$L_1(Z_3, u) = (12'',13,14,8')$&&\footnotesize$(8',12'')$&\\
\footnotesize$L(Z_1, o) = [9,11,14,13]$&&&\\
\footnotesize$L_0(Z_2, o) = [12',10,8'']$&&&\\
\footnotesize$L_1(Z_2, o) = (8',12'')$&&&\\
\footnotesize$L(Z_3) = [a_3,8'',9,10,11,12',a_3,7,6,\textbf{2},1,4,3,5]$&&&\\

\hline
\footnotesize$[\textbf{5},a_3,8'',9,10,11,12',a_3,7,6]$&\footnotesize$(8'',10)(10,8'')$&\footnotesize$[8'',11,\textbf{8}',12'',8'',12',11,12']$&\footnotesize$(8'',8'')$\\
\footnotesize$[\textbf{5},6,7]$&&\footnotesize$(\textbf{8}',12'')$&\\
\footnotesize$(12'',13,14,8')$&&&\\
\footnotesize$[9,11,14,13]$&&&\\
\footnotesize$[12',10,8'']$&&&\\
\footnotesize$(8',12'')$&&&\\
\hline
\end{tabular}
\end{center}

\begin{center}
\begin{tabular}[c]{|l|l||l|l|}
\hline
\footnotesize$[8'',\textbf{9},10,11,12']$&\footnotesize$(8'',10)(10,8'')$&\footnotesize$[8'',11,8'',12',11,12']$&\footnotesize$(8'',8'')$\\
\footnotesize$(12'',13,14,8')$&&&\\
\footnotesize$[\textbf{9},11,14,13]$&&&\\
\footnotesize$[12',10,8'']$&&&\\
\footnotesize$(8',12'')$&&&\\

\hline
\footnotesize$[8'',11,14,\textbf{13},10,11,12']$&\footnotesize$(8'',10)(10,8'')$&&\\
\footnotesize$(12'',\textbf{13},14,8')$&&&\\
\footnotesize$[12',10,8'']$&&&\\
\footnotesize$(8',12'')$&&&\\
\hline
\end{tabular}
\end{center}
$(8'',10)$ is in $L(Z_3)$ in the part $L(Z_3,u)$. The number $t_0$ with minimal absolute value such that $p_1^{t_0}(8)=10$ is $t_0=2$. So we have $w_0=x^2y^{-1}$.\\
$(10,8'')$ is in $L(Z_2,o)$. The number $t_1$ with minimal absolute value such that $p_1^{t_1}p_2^{-1}(10)=p_2^{-1}(8)$ is $t_1=-1$. So we have $w_1=x^{-1}y$.\\
The new curve is
$$K_4:= (8, x^2y^{-1}x^{-1}y).$$
\end{spacing}

\begin{figure}[h]
\begin{center}
\includegraphics[scale=0.6]{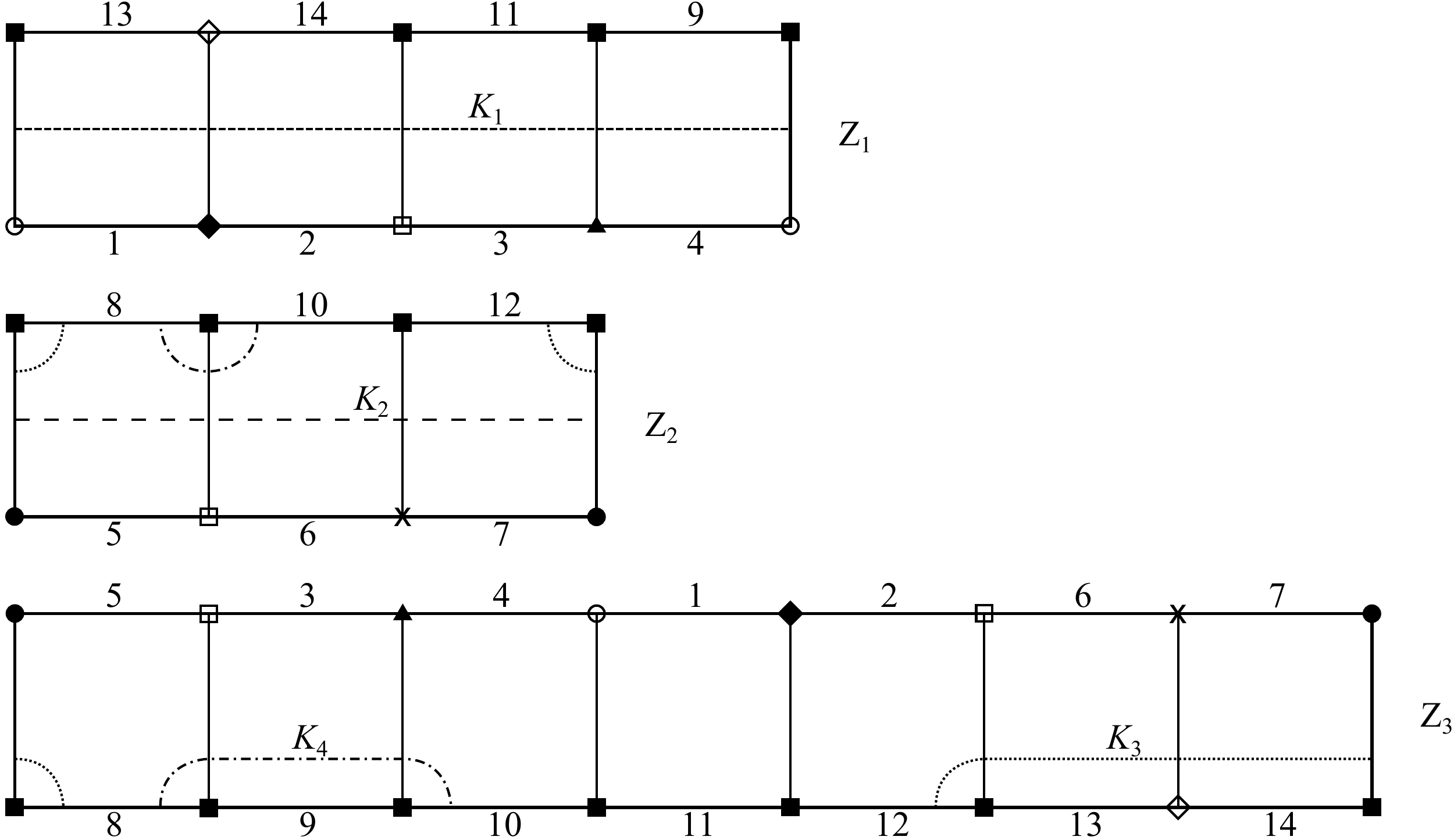}
\end{center}
\caption{The origami $O_{14}$ with the HSS $\{K_1,\dots,K_4\}$}
\end{figure}~\\

\chapter{Examples}
\label{kapbsp}
\section{L-origamis}
\label{lori}
\begin{definition}
For $m,n\geq2$, let \begriff{L(m,n)} be the origami that consists of $m+n-1$ squares and is defined by the following permutations:
\[p_1=(1,\cdots,m)\quad\text{ and }\quad p_2=(1,m+1,\dots,m+n-1)\]
We refer to the origamis $L(m,n)$ also as \begriff{L-origamis}.
\end{definition}
\begin{figure}[h]
\begin{center}
\includegraphics[scale=0.8]{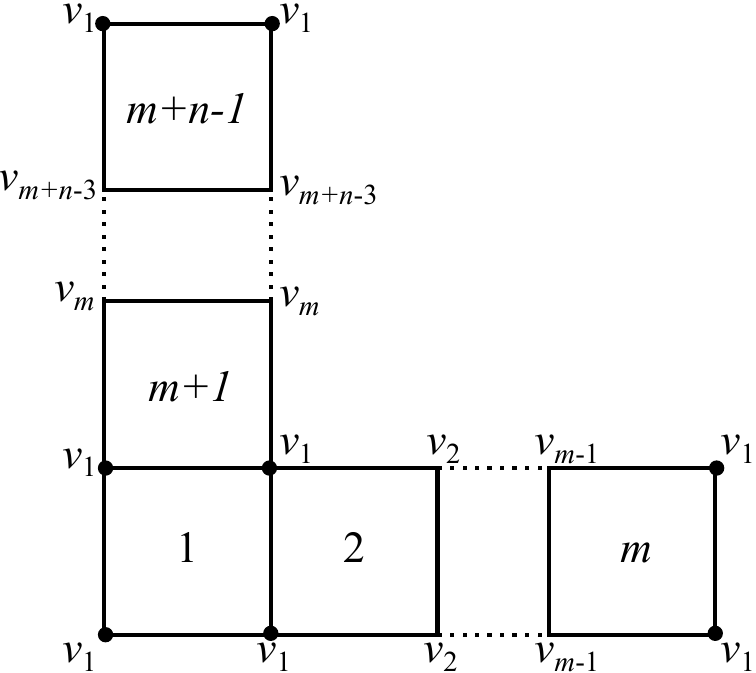}
\end{center}
\caption{Opposed sides are identified}
\label{Lmn}
\end{figure}~\\
An $L$-origami $L(m,n)$ has always genus $g=2$ and the ramification point $v_1$ of order 3. The punctures $v_2,\dots,v_{m+n-3}$ are not ramificated.\\
The subgroup of $F_2$ associated to $L(m,n)$ is $H:=\langle h_1,\dots,h_{m+n}\rangle$, where
\[h_1=x^m, \quad h_i=\begin{cases}
x^{m-(i-1)}yx^{-m+i-1},& 2\leq i\leq m\\
y^{i-m}xy^{m-i},& m+1\leq i\leq m+n-1
\end{cases}\text{ and }\quad h_{m+n}=y^n
\]
The paths around the punctures are $r_1=(x^{-1}y^{-1}xy)^3=h_1^{-1}h_2^{-1}h_1h_{m+1}^{-1}h_mh_{m+n}^{-1}h_{m+n-1}h_{m+n}$ and
\[r_i=\begin{cases}
h_ih_{i+1}^{-1},& 2\leq i\leq m-1\\
h_{i+1}h_{i+2}^{-1},& m\leq i\leq m+n-3.
\end{cases}
\]
Let $X$ be the surface corresponding to $L(m,n)$, and let $X^*$ be the surface corresponding to the punctured origami. We have
\[
\pi_1(X^*)= H \quad \text{ and }\quad\pi_1(X)=\langle H\mid R\rangle,\quad\text{ with }\quad R:=\langle\langle r_1,\dots,r_{m+n-3}\rangle\rangle.
\]
A HSS is given by the horizontal cut around the horizontal cylinder of length $m$ and a cut around a horizontal cylinder of length 1.
A set of symplectic generators $(a_1,a_2,b_1,b_2)$ such that $a_1$ and $a_2$ are conjugate horizontal is the following:
\[a_1 = \overline{x^{-m}} ~~~ a_2 = \overline{yxy^{-1}} ~~~ b_1 = \overline{xyx^{-1}} ~~~ b_2 = \overline{yx^{-1}y^{-1}xyx^{-1}y^{-n}}\]
or
\[a_1 = \overline{h_1^{-1}} ~~~ a_2 = \overline{h_{m+1}} ~~~ b_1 = \overline{h_m} ~~~ b_2 = \overline{h_{m+1}^{-1}h_mh_{m+n}^{-1}}.\]
\begin{figure}[h]
\begin{center}
\includegraphics[scale=0.8]{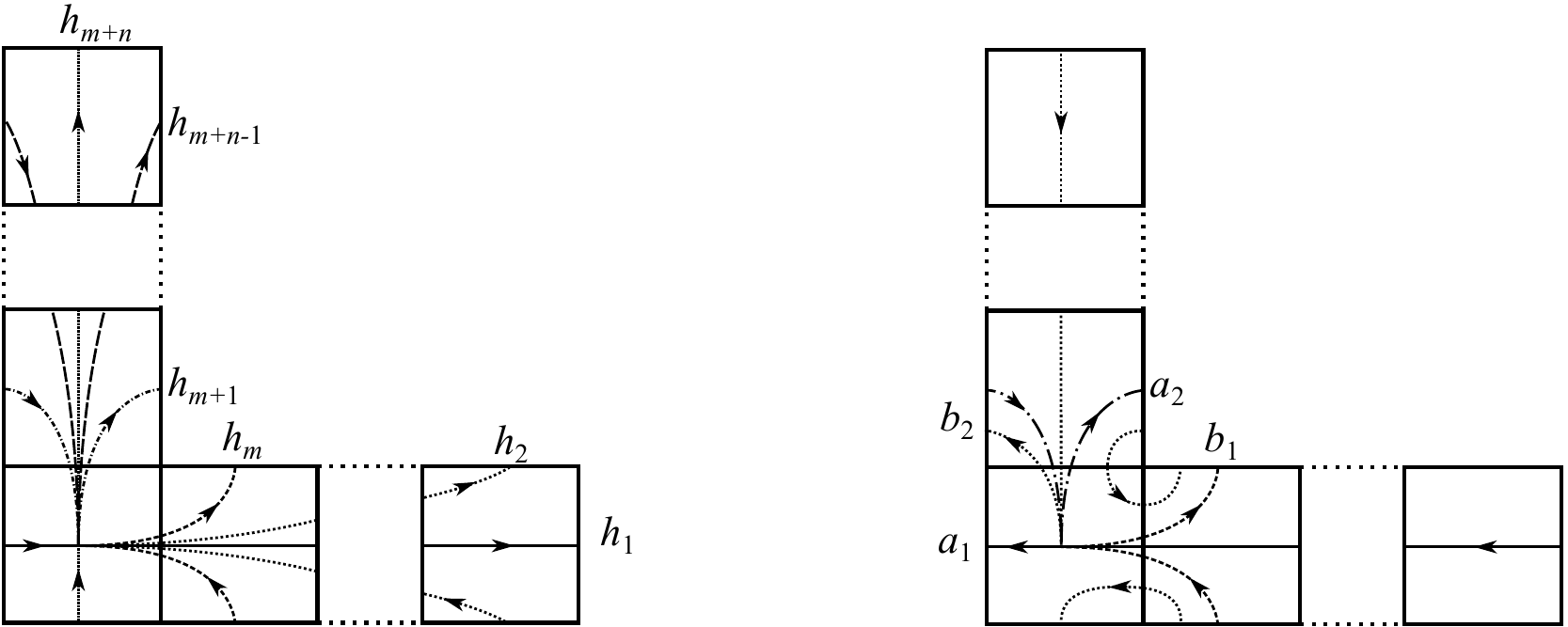}
\end{center}
\caption{The sets of generators for $\pi_1(X^*)=H$ (left) and of $\pi_1(X)$ (right)}
\label{Lmnerzsys}
\end{figure}\\
In $\Aff^+(X)$, we have the element
\[f\colon z\mapsto \begin{pmatrix} 1&m\\0&1\end{pmatrix}\cdot z,\]
that can be lifted to the element $\hat{f}_\star\in\Stab_{\Aut^+(F_2}(H)$, where
\begin{align*}
\hat{f}_\star:\;&F_2\ra F_2\\
&x\mapsto x,y\mapsto x^my.
\end{align*}
We get
\[\hat{f}_\star(h_1)=x^m, \quad \hat{f}_\star(h_{m+n})=(x^my)^n\]
and
\[\hat{f}_\star(h_i)=\begin{cases}
                      x^{m-(i-1)}x^myx^{-m+i-1},& 2\leq i\leq m\\
                     (x^my)^{i-m}x(x^my)^{m-i},& m+1\leq i\leq m+n-1.
\end{cases}\]
That means:
\[\textstyle\hat{f}_\star(h_1)=h_1,\quad \hat{f}_\star(h_{m+n})=h_1(\prod_{k=m+1}^{m+n-1}h_k^m)h_{m+n}\]
and
\[\hat{f}_\star(h_i)=\begin{cases}
                     h_1h_i,& 2\leq i\leq m\\
                     h_1(\prod_{k=m+1}^{i-1}h_k^m)h_i(\prod_{k=m+1}^{i-1}h_k^m)^{-1}h_1^{-1},& m+1\leq i\leq m+n-1.
\end{cases}\]
So we get
\begin{align*}
f_*(a_1)& = \overline{\hat{f}_\star(h_1^{-1})}=\overline{h_1^{-1}}=a_1\\
~\\
f_*(a_2)& = \overline{\hat{f}_\star(h_{m+1})}\\
&=\textstyle\overline{h_1(\prod_{k=m+1}^mh_k^m)h_{m+1}(\prod_{k=m+1}^mh_k^m)^{-1}h_1^{-1}}\\
&=\overline{h_1h_{m+1}h_1^{-1}}\\
&=a_1^{-1}a_2a_1\\
~\\
f_*(b_1)& = \overline{\hat{f}_\star(h_m)}=\overline{h_1h_m}=a_1^{-1}b_1\\
~\\
f_*(b_2) &= \overline{\hat{f}_\star(h_{m+1}^{-1}h_mh_{m+n}^{-1})}\\
&=\textstyle\overline{h_1h_{m+1}^{-1}h_1^{-1}}\cdot\overline{h_1h_m}\cdot\overline{(h_1(\prod_{k=m+1}^{m+n-1}h_k^m)h_{m+n})}^{-1}\\
&=\textstyle\overline{h_1h_{m+1}^{-1}h_mh_{m+n}^{-1}(\prod_{k=m+1}^{m+n-1}h_k^m)^{-1}h_1^{-1}}\\
&=\overline{h_1h_{m+1}^{-1}h_mh_{m+n}^{-1}h_{m+1}^{(1-n)m}h_1^{-1}}\quad\quad(\text{since } r_m,\dots,r_{m+n-3})\\
&=a_1^{-1}a_2^{-1}b_1b_1^{-1}a_2b_2a_2^{(1-n)m}a_1\quad\quad(\text{since } \overline{h_{m+n}} = b_2^{-1}a_2^{-1}b_1)\\
&=a_1^{-1}b_2a_2^{(1-n)m}a_1.
\end{align*}
So for the symplectic homomorphism $\alpha\colon \pi_1(X)\ra F_g$ with $\alpha(a_i)=1,\alpha(b_i)=\gamma_i$ for $i\in\{1,\dots,g\}$, we finally get
\begin{align*}
\alpha\circ f_*(a_1)&=\alpha(a_1)=1\\
\alpha\circ f_*(a_2)&=\alpha(a_1^{-1}a_2a_1)=1\\
\alpha\circ f_*(b_1)&=\alpha(a_1^{-1}b_1)=\gamma_1\\
\alpha\circ f_*(b_2)&=\alpha(a_1^{-1}b_2a_2^{(1-n)m}a_1)=\gamma_2.
\end{align*}

\section{X-origamis}
\begin{definition}
For $n\geq1$, let \begriff{O(n)} be the origami that consists of $2n$ squares and is defined by the following permutations:
\[p_1=(1,\cdots,n)\quad\text{ and }\quad p_2=(1,2)(3,4)\dots(2n-1,2n)\]
\end{definition}
\begin{figure}[h]
\begin{center}
\includegraphics[scale=0.8]{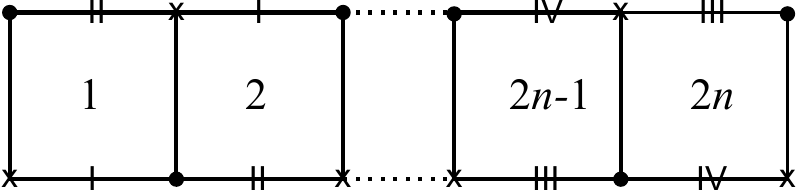}
\end{center}
\label{On}
\end{figure}~\\
We refer to the origamis $O(n)$ also as \begriff{X-origamis}. They are introduced in \cite[§4.5]{gabidiss} and their Veech group is given in \cite[§5.2]{gabidiss}.\\
The origami $O(n)$ has genus $g=n$ and 2 ramification points of order $n$. The subgroup of $F_2$ associated to $O(n)$ is $H:=\langle h_1,\dots,h_{2n+1}\rangle$, where
\[h_{2i-1}=x^{2i-2}yx^{1-2i}, \quad h_{2i}=x^{2i-1}yx^{2-2i},\text{ for }1\leq i\leq n \quad\text{ and } \quad h_{2n+1}=x^{2n}.\]
The paths around the punctures are \[r_1=(x^{-1}y^{-1}xy)^n=h_{2n+1}^{-1}\prod_{k=1}^nh_{2(n-k)+1}^{-1}h_{2(n-k)+2}\quad\text{ and }\quad
r_2=(yx^{-1}y^{-1}x)^n=h_1\prod_{k=1}^nh_{2k}^{-1}h_{2k+1}\]
We have
\[
\pi_1(X^*)= H \quad \text{ and }\quad\pi_1(X)=\langle H\mid R\rangle,\quad\text{ with }\quad R:=\langle\langle r_1, r_2\rangle\rangle.
\]
\begin{figure}[h]
\begin{center}
\includegraphics[scale=0.9]{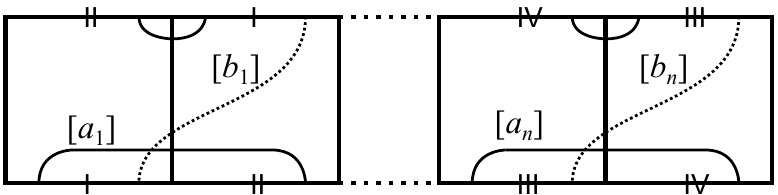}
\end{center}
\caption{Horizontal Schottky cut system $[a_i]$, extended to a symplectic system}
\label{Onsymplungerade2}
\end{figure}~\\
A HSS is given by the curves
\[K_i=(2i-1,xy^{-1}xy)\text{ for } i\in\{1,\dots,n\}.\]
With the curves
\[K'_i=(2i-1,xy),\]
we get a symplectic system, and we have then a symplectic set of generators of $\pi_1(X)$ with base point in square 1:
\[a_i = \overline{v_ixy^{-1}xyv_i^{-1}}\qquad b_i = \overline{v_ixyv_i^{-1}}\]
or
\[a_i = \overline{w_ih_{2i-1}^{-1}h_{2i}w_i^{-1}}\qquad b_i = \overline{w_ih_{2i}w_i^{-1}}.\]
for suitable $v_1,\dots,v_n\in F_2$, $w_1,\dots,w_n\in H$ and $i\in\{1,\dots,n\}$.\\
In $\Aff^+(X)$, we have the element
\[f\colon z\mapsto \begin{pmatrix} 1&2n\\0&1\end{pmatrix}\cdot z,\]
that can be lifted to the element $\hat{f}_\star\in\Stab_{\Aut^+(F_2}(H)$, where
\begin{align*}
\hat{f}_\star:\;&F_2\ra F_2\\
&x\mapsto x,y\mapsto x^{2n}y.
\end{align*}
We get
\[\hat{f}_\star(h_i)=h_{2n+1}h_i, \quad i\in\{1,\dots,2n\}\quad\text{ and } \quad\hat{f}_\star(h_{2n+1})=h_{2n+1}\]
So we get
\begin{align*}
f_*(a_i)& = \overline{\hat{f}_\star(w_ih_{2i-1}^{-1}h_{2i}w_i^{-1})}=\\
&=\overline{\hat{f}_\star(w_i)\cdot h_{2i-1}^{-1}h_{2i}\cdot\hat{f}_\star(w_i^{-1})}\\
&=\overline{\underbrace{\hat{f}_\star(w_i)w_i^{-1}}_{:=c_i\in\pi_1(X)}}\cdot\overline{w_ih_{2i-1}^{-1}h_{2i}w_i^{-1}}\cdot\overline{w_i\hat{f}_\star(w_i^{-1})}\\
&=c_ia_ic_i^{-1}\\
~\\
f_*(b_i) &= \overline{\hat{f}_\star(w_ih_{2i}w_i^{-1})}\\
&= \overline{\hat{f}_\star(w_i)h_{2n+1}h_{2i}\hat{f}_\star(w_i^{-1})}\\
&= \overline{\hat{f}_\star(w_i)h_{2n+1}w_i^{-1}w_ih_{2i}w_i^{-1}w_i\hat{f}_\star(w_i^{-1})}\\
&= \overline{\hat{f}_\star(w_i)h_{2n+1}w_i^{-1}}\cdot b_i\cdot\overline{w_i\hat{f}_\star(w_i^{-1})}
\end{align*}
and so we get for the symplectic homomorphism $\alpha\colon \pi_1(X)\ra F_g$ with $\alpha(a_i)=1,\alpha(b_i)=\gamma_i$ for $i\in\{1,\dots,g\}$:
\begin{align*}
\alpha\circ f_*(a_i)&=\alpha(c_ia_ic_1^{-1})=1\\
\alpha\circ f_*(b_i)&=\alpha(\overline{\hat{f}_\star(w_{i})h_{2n+1}w_{i}^{-1}}\cdot b_i\cdot\overline{w_i\hat{f}_\star(w_i^{-1})})&&(\overline{h_{2n+1}}\in\Kern(\alpha)\text{ by Thm. }\ref{zylhorimkern})\\
&=\alpha(\overline{\hat{f}_\star(w_i)w_i^{-1}})\cdot \gamma_i\cdot\alpha(\overline{w_i\hat{f}_\star(w_i^{-1})})\\
&=\gamma_i
\end{align*}
The last equality holds because of $\overline{w_i\hat{f}_\star(w_i^{-1})}\in\Kern(\alpha)$, by proof of Thm. \ref{schnittnichtleer}.
\section{$\Aff^+(X)\cap\Mod_g(\alpha)$ for a Schottky cut system in direction $(1,1)$}
\begin{figure}
\begin{center}
\includegraphics[scale=0.6]{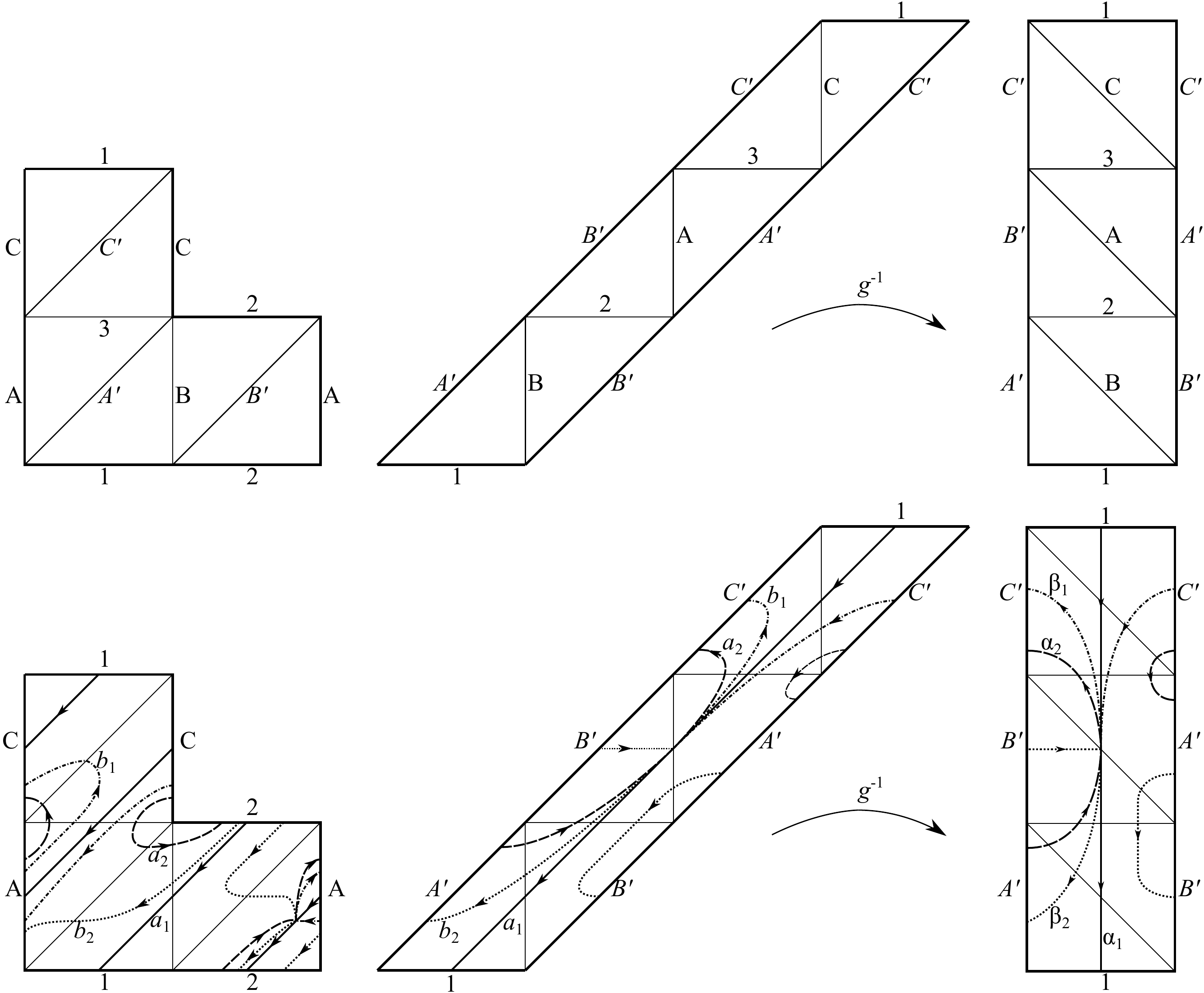}
\caption{A suitable set of symplectic generators of the fundamental group of $L(2,2)$}
\label{L22schief}
\end{center}
\end{figure}

Take the origami $L(2,2)$ and the vector $v=(1,1)\in\Z^2$. We look for a parabolic element $A\in\Gamma(L(2,2))$ with eigenvector $v$ and for a symplectic homomorphism $\alpha\colon \pi_1(X)\ra F_2=\langle\gamma_1,\gamma_2\rangle$ such that for a lift $f\in\Aff^+(X)$ of $A$ we have $\Aff^+(X)\cap\Mod_g(\alpha)\cong\Z$, where $X$ is the surface corresponding to this origami.\\
First, we construct a set of symplectic generators $(a_1,b_1,a_2,b_2)$ of the fundamental group with $\alpha(a_i)=1$ and $\alpha(b_i)=\gamma_i$ for $i\in\{1,2\}$.\\ As in the construction in the proof of Thm. \ref{schief}, we get from $L(2,2)$ by shearing with
\begin{align*}
g^{-1}:\; X&\longrightarrow X'\\
z&\longmapsto
\begin{pmatrix}
1&-1\\
0&1
\end{pmatrix}\cdot z
\end{align*}
a new origami $L(2,2)'$, see figure \ref{L22schief} at the top on the right. There we find a set of symplectic generators $(\alpha_1,\beta_1,\alpha_2,\beta_2)$ such that $\alpha_1$ and $\alpha_2$ are conjugate vertical (see figure \ref{L22schief} on the bottom right):
\begin{align*}
\alpha_1 &= \overline{y^{-3}}\\
\alpha_2 &= \overline{yx^{-1}y^{-1}xy}\\
\beta_1 &= \overline{yx^{-1}y^{-1}}\\
\beta_2 &= \overline{y^{-1}x^{-1}y^{-1}x} 
\end{align*}
Let $(a_1,b_1,a_2,b_2)$ be its preimage (figure \ref{L22schief} on the bottom left). We have:
\begin{align*}
a_1 &= \overline{(xy)^{-3}}\\
a_2 &= \overline{xyx^{-1}y^{-1}xy}\\
b_1 &= \overline{xyx^{-1}y^{-1}x^{-1}}\\
b_2 &= \overline{y^{-1}x^{-2}y^{-1}}
\end{align*}
With the notations from Thm. \ref{schief} we, furthermore, have $\der(f')=\begin{pmatrix}1&0\\3&1\end{pmatrix}$ and, therefore,
\begin{align*}
A&=\der(g)\der(f')\der(g^{-1})\\
&=\begin{pmatrix}1&1\\0&1\end{pmatrix}\begin{pmatrix}1&0\\3&1\end{pmatrix}\begin{pmatrix}1&-1\\0&1\end{pmatrix}\\
&=\begin{pmatrix}4&-3\\3&-2\end{pmatrix}
\end{align*}
A lift $\hat{f}_{\star}$ of $f$ to $\Aut^+(F_2)$ that maps to itself the subgroup $H\leq F_2$ corresponding to the punctured origami is given by the algorithm from \cite[§ 4.1]{gabidiss}, followed by conjugation with $x$:
\begin{align*}
\hat{f}_{\star}:\;F_2&\ra F_2\\
x&\mapsto (xy)^3x,\; y\mapsto (yx)^{-2}x^{-1}
\end{align*}
A free set of generators of $H$ is
\begin{align*}
h_1 &= x^2\\
h_2 &= y\\
h_3 &= x^{-1}yxy^{-1}x\\
h_4 &= x^{-1}y^2x.
\end{align*}
We get in the same way as in §\ref{lori}
\begin{align*}
f_*(a_1)&=a_1\\
f_*(b_1)&=b_1a_1\\
f_*(a_2)&=a_2\\
f_*(b_2)&=b_2,
\end{align*}
and get $\alpha\circ f_*=\alpha$.

\newpage

\section{Example for a flat surface}
\begin{figure}[h]
\begin{center}
\includegraphics[scale=0.55]{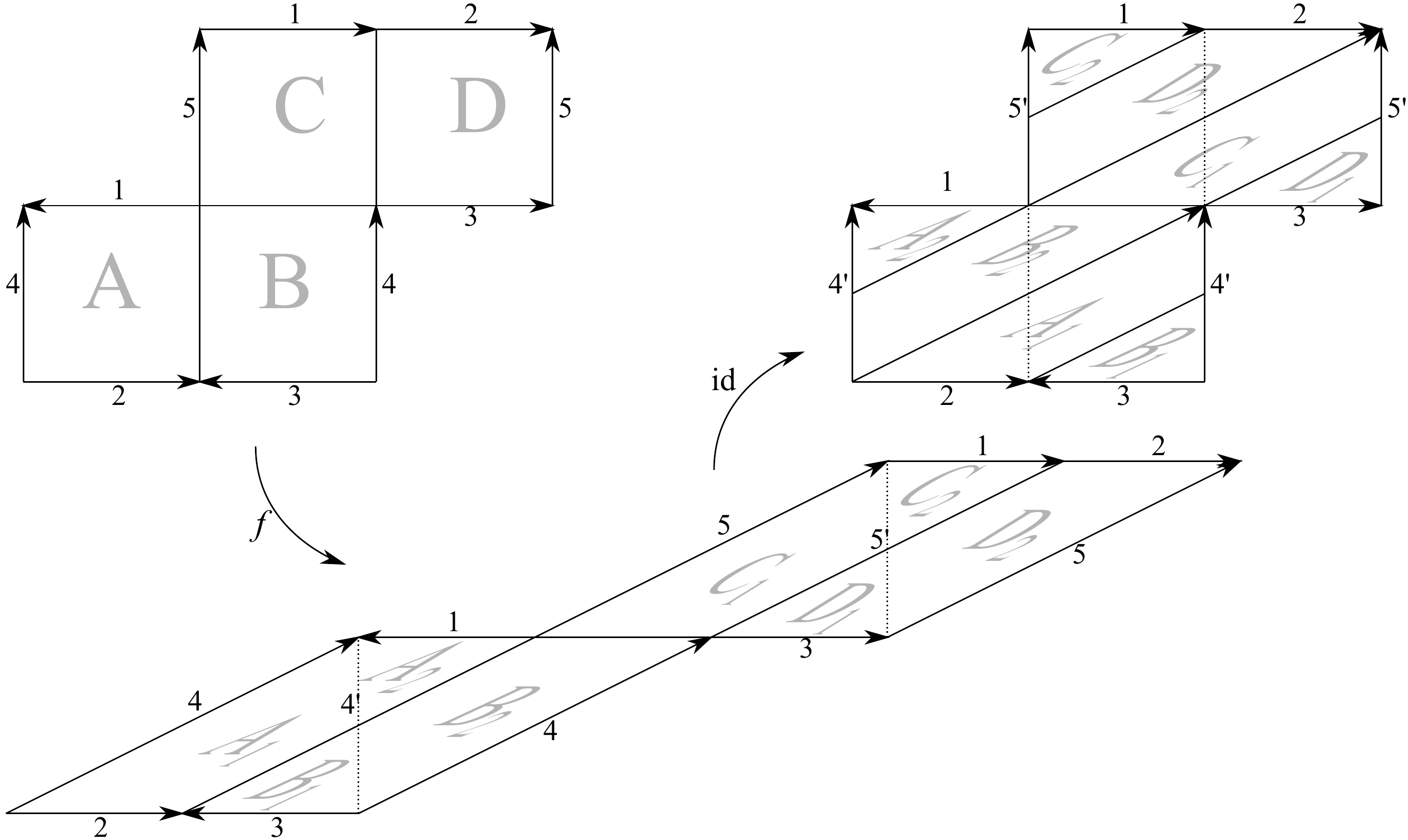}
\end{center}
\caption{The affine element $f$}
\label{Halbtransl1}
\end{figure}
\begin{figure}[h]
\begin{center}
\includegraphics[scale=0.55]{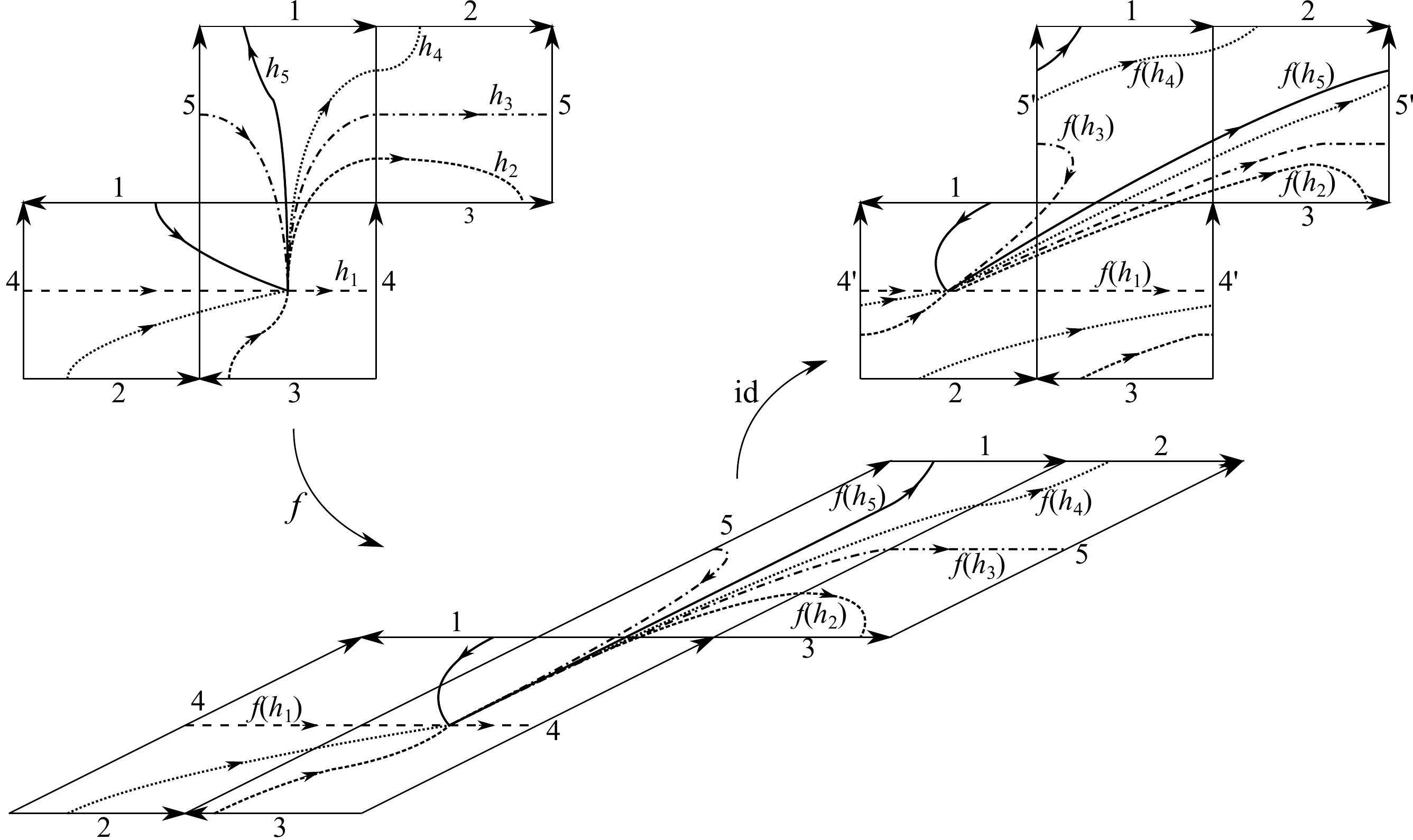}
\end{center}
\caption{The action of $f$ on $\pi_1(X^*)=\langle h_1,\dots,h_5\rangle$}
\label{Halbtransl2}
\end{figure}
\begin{figure}[h]
\begin{center}
\includegraphics[scale=0.55]{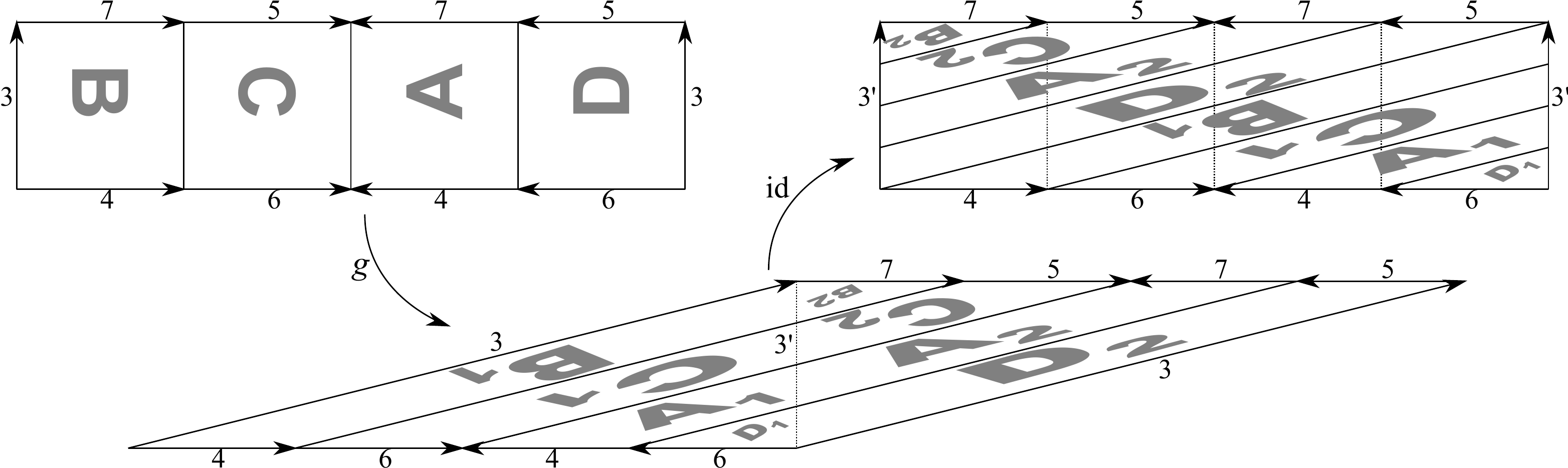}
\end{center}
\caption{The affine element $g$}
\label{Halbtransl3}
\end{figure}
\begin{figure}[h]
\begin{center}
\includegraphics[scale=0.55]{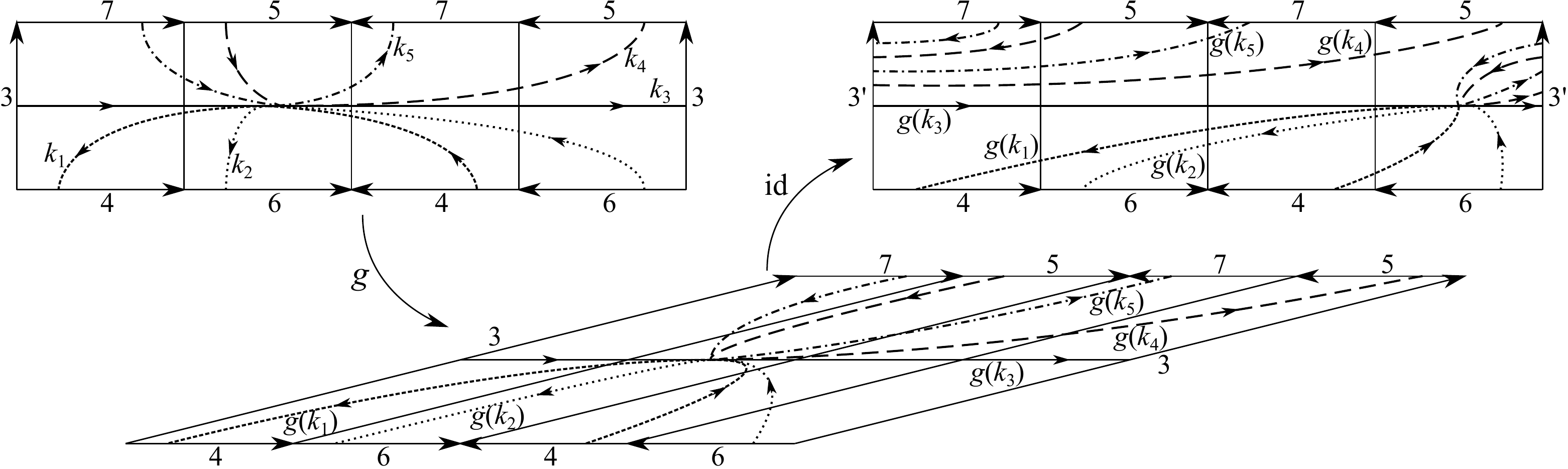}
\end{center}
\caption{The action of $g$ on $\pi_1(X^*)=\langle k_1,\dots,k_5\rangle$}
\label{Halbtransl4}
\end{figure}\begin{figure}[h]
\begin{center}
\includegraphics[scale=0.55]{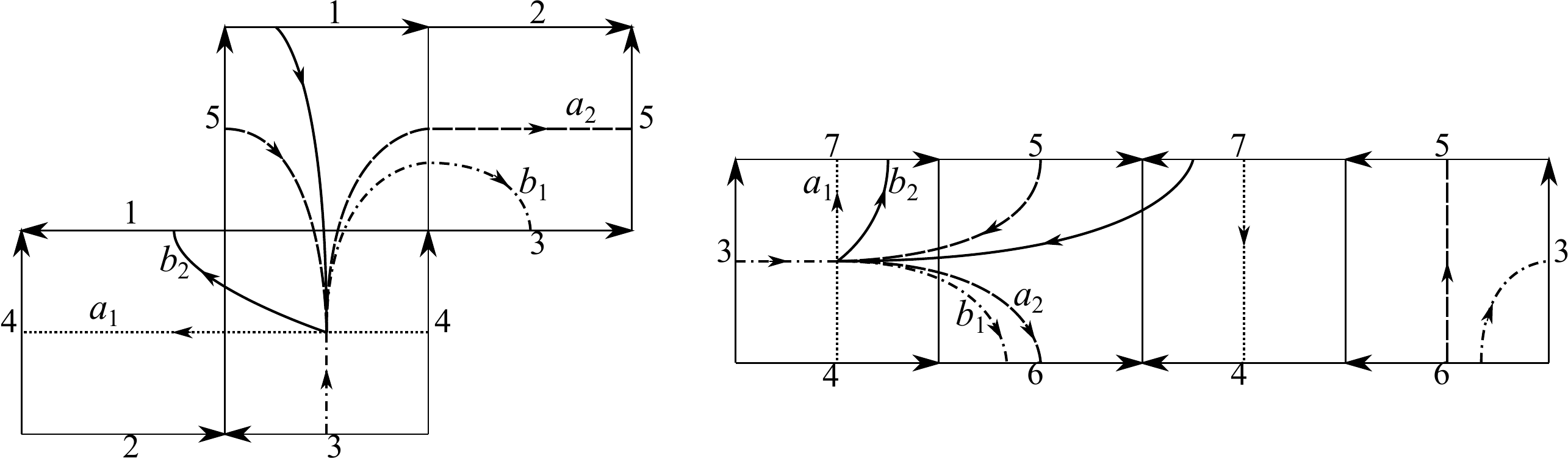}
\end{center}
\caption{A symplectic set of generators of $\pi_1(X)$}
\label{Halbtransl5}
\end{figure}~\\
Take the flat surface $X$ given on the top left in figure \ref{Halbtransl1} and in \ref{Halbtransl3} (after rotating by $90^\circ$). In figure \ref{Halbtransl1}, we see how the element $f\in\Aff^+(X)$ that is locally given by
\[f\colon z\mapsto \begin{pmatrix}1&2\\0&1\end{pmatrix}\cdot z\]
acts on $X$, and in figure \ref{Halbtransl3} we see how the element $g\in\Aff^+(X)$ that is locally given by
\[g\colon z\mapsto \begin{pmatrix}1&0\\4&1\end{pmatrix}\cdot z\]
acts on $X$. Figure \ref{Halbtransl5} shows a symplectic set of generators $(a_1, b_1, a_2, b_2)$.
We will show that the element $g$ acts trivially on the homology and give three symplectic homomorphisms $\alpha_1,\alpha_2$ and $\alpha_3$ such that $g\in\Aff^+(X)\cap\Mod_2(\alpha_i)$ for $i\in\{1,2,3\}$. Moreover, we'll give a further symplectic homomorphism $\alpha_4$ such that $f\in\Aff^+(X)\cap\Mod_2(\alpha_4)$, and we'll show that the eigenvalues of $\der(f\circ g)$ are not eigenvalues of $M_{f\circ g}=M_f$. (Since we have a flat surface, only $D(f), D(g), D(f\circ g)\in\PSL_2(\R)$ are well defined, but the eigenvalues of their preimages in $\SL_2(\R)$ do (up to sign) not depend from  the choice of the preimages. Hence, we can talk about eigenvalues (up to sign) of elements in $\PSL_2(\R)$. Therefore, we have here examples for elements contained in groups of the form $\Aff^+(X)\cap\Mod_g(\alpha)$ and examples which show that Thm. \ref{nichthypallgemein}, and Cor. \ref{nichtintorelli} can \textit{not} be generalised to the case of flat surfaces.\\
In figure \ref{Halbtransl2}, a set of generators $h_1,\dots,h_5$ of the fundamental group of the punctured surface $X^*$ is given and another one -- $k_1,\dots,k_5$ -- is given in figure \ref{Halbtransl4}. The paths around the two punctures are given by
\[r_1 = h_1h_4^{-1}h_5h_1^{-1}h_2,\quad\quad\text{ and }\quad\quad r_2= h_2^{-1}h_3h_5^{-1}h_3^{-1}h_4\]
or, with the generators $k_1,\dots,k_5$, given by
\[r_1' = k_1k_2^{-1}k_1^{-1}k_2k_3,\quad\quad\text{ and }\quad\quad r_2' = k_5^{-1}k_4k_5k_4^{-1}k_3.\]
We get $\pi_1(X)$ if we factorize out of $\pi_1(X^*)$ the normal subgroup generated by the paths around the punctures. The elements $a_1, a_2, b_1, b_2$ in figure \ref{Halbtransl5} are given by
\begin{center}
$\begin{array}[c]{cll}
a_1&=\overline{k_5^{-1}k_1^{-1}}&=\overline{h_1^{-1}}\\
a_2&=\overline{k_2k_4}&=\overline{h_3}\\
b_1&=\overline{k_2k_3}&=\overline{h_2}\\
b_2&=\overline{k_5^{-1}}&=\overline{h_5^{-1}}
\end{array}$
\end{center}
Conversely, we can also express $\overline{k_1},\dots,\overline{k_5},\overline{h_1},\dots,\overline{h_5}$ with $a_1,a_2,b_1,b_2$:
\begin{center}
$\begin{array}[c]{clccll}
\overline{h_1}& = a_1^{-1}&&\overline{k_1}&=a_1^{-1}\overline{k_5^{-1}}= a_1^{-1}b_2\\
\overline{h_2}& = b_1&&\overline{k_2}&\stackrel{r_1'}{=}\overline{k_1^{-1}k_2k_3k_1}=\overline{k_1^{-1}}b_1\overline{k_1}=b_2^{-1}a_1b_1a_1^{-1}b_2\\
\overline{h_3}& = a_2&\qquad&\overline{k_3}&= k_2^{-1}b_1=b_2^{-1}a_1b_1^{-1}a_1^{-1}b_2b_1\\
\overline{h_4}& \stackrel{r_1}{=} h_5h_1^{-1}h_2h_1=b_2^{-1}a_1b_1a_1^{-1}&&\overline{k_4}&= k_2^{-1}a_2=b_2^{-1}a_1b_1^{-1}a_1^{-1}b_2a_2\\
\overline{h_5}& = b_2^{-1}&&\overline{k_5}&= b_2^{-1}
\end{array}$
\end{center}
In figure \ref{Halbtransl2}, we see how $f$ acts on $h_1,\dots,h_5$ and in figure \ref{Halbtransl4} we see how $g$ acts on $k_1,\dots,k_5$:
\begin{center}
$\begin{array}[c]{clccll}
f_*(h_1) & = h_1&&g_*(k_1) & = k_1\\
f_*(h_2) & = h_2h_1&&g_*(k_2) & = k_2\\
f_*(h_3) & = h_3&\quad\quad\qquad\quad\quad&g_*(k_3) & = k_3\\
f_*(h_4) & = h_3h_4h_1&&g_*(k_4) & = k_3k_4k_3^{-1}\\
f_*(h_5) & = h_3h_5&&g_*(k_5) & = k_3k_5k_3^{-1}
\end{array}$
\end{center}
We now see
\begin{align*}
g_*(a_1)&=g_*(\overline{k_5^{-1}k_1^{-1}})\\
&= \overline{k_3k_5^{-1}k_3^{-1}}\cdot\overline{k_1^{-1}}\\
&= (b_2^{-1}a_1b_1^{-1}a_1^{-1}b_2b_1)b_2(b_1^{-1}b_2^{-1}a_1b_1a_1^{-1}b_2)(b_2^{-1}a_1)\\
&= b_2^{-1}a_1b_1^{-1}a_1^{-1}b_2b_1b_2b_1^{-1}b_2^{-1}a_1b_1\\
g_*(a_2)&= g_*(\overline{k_2k_4}) \\
&= \overline{k_2}\cdot\overline{k_3k_4k_3^{-1}}\\
& = (b_2^{-1}a_1b_1a_1^{-1}b_2)(b_2^{-1}a_1b_1^{-1}a_1^{-1}b_2b_1)(b_2^{-1}a_1b_1^{-1}a_1^{-1}b_2a_2)(b_1^{-1}b_2^{-1}a_1b_1a_1^{-1}b_2)\\
& = b_1b_2^{-1}a_1b_1^{-1}a_1^{-1}b_2a_2b_1^{-1}b_2^{-1}a_1b_1a_1^{-1}b_2\\
g_*(b_1)&=g_*(\overline{k_2k_3})=\overline{k_2k_3}=b_1\\
g_*(b_2)&= g_*(\overline{k_5^{-1}}) = \overline{k_3k_5^{-1}k_3^{-1}}=(b_2^{-1}a_1b_1^{-1}a_1^{-1}b_2b_1)b_2(b_1^{-1}b_2^{-1}a_1b_1a_1^{-1}b_2)
\end{align*}
and
\begin{align*}
f_*(a_1)&=f_*(\overline{h_1^{-1}})= \overline{h_1^{-1}}=a_1\\
f_*(a_2)&= f_*(\overline{h_3}) = \overline{h_3} = a_2\\
f_*(b_1)&=f_*(\overline{h_2})=\overline{h_2h_1}=b_1a_1^{-1}\\
f_*(b_2)&= f_*(\overline{h_5^{-1}}) = \overline{h_5^{-1}h_3^{-1}}=b_2a_2^{-1}.
\end{align*}
We have that $M_g$ is the identity, i.e.\ $g$ acts trivially on the homology and is, therefore, in the Torelli group. So this example shows that the assertion for translation surfaces in Cor. \ref{nichtintorelli} cannot be true for flat surfaces. Moreover, we have
\[M_f=
\begin{pmatrix}
1&0&-1&0\\
0&1&0&-1\\
0&0&1&0\\
0&0&0&1
\end{pmatrix} = M_{f\circ g},\]
and $M_{f\circ g}$ has the eigenvalue 1 with multiplicity 4. But
\[\der(f\circ g)=\der(f)\cdot\der(g)=
\begin{pmatrix}
1&2\\
0&1
\end{pmatrix}\cdot
\begin{pmatrix}
1&0\\
4&1
\end{pmatrix}=
\begin{pmatrix}
9&2\\
4&1
\end{pmatrix}\]
has the eigenvalues $5+2\sqrt{6}$ and $5-2\sqrt{6}$. So the assertion for translation surfaces in Thm. \ref{nichthypallgemein} cannot be true for flat surfaces.\\
Now, let  the symplectic homomorphisms $\alpha_1,\dots,\alpha_4\colon \pi_1(X)\ra F_2=\langle\gamma_1,\gamma_2\rangle$ be given as follows:
\[\begin{array}[c]{clclclcl}
\alpha_1:
&\alpha_1(a_1)=\gamma_1&\alpha_2:&\alpha_2(a_1)=\gamma_1&\alpha_3:&\alpha_3(a_1)=1&\alpha_4:&\alpha_4(a_1)=1\\
&\alpha_1(a_2)=\gamma_2&&\alpha_2(a_2)=1&&\alpha_3(a_2)=\gamma_1&&\alpha_4(a_2)=1\\
&\alpha_1(b_1)=1&&\alpha_2(b_1)=1&&\alpha_3(b_1)=\gamma_2&&\alpha_4(b_1)=\gamma_1\\
&\alpha_1(b_2)=1&&\alpha_2(b_2)=\gamma_2&&\alpha_3(b_2)=1&&\alpha_4(b_2)=\gamma_2\\
\end{array}\]
We have
\begin{align*}
\alpha_1\circ g_*(a_1)&=\alpha_1(b_2^{-1}a_1b_1^{-1}a_1^{-1}b_2b_1b_2b_1^{-1}b_2^{-1}a_1b_1)= \gamma_1\\
\alpha_1\circ g_*(a_2)&=\alpha_1(b_1b_2^{-1}a_1b_1^{-1}a_1^{-1}b_2a_2b_1^{-1}b_2^{-1}a_1b_1a_1^{-1}b_2)= \gamma_2\\
\alpha_1\circ g_*(b_1)&=\alpha_1(b_1)= 1\\
\alpha_1\circ g_*(b_2)&=\alpha_1(b_2^{-1}a_1b_1^{-1}a_1^{-1}b_2b_1b_2b_1^{-1}b_2^{-1}a_1b_1a_1^{-1}b_2)= 1
\end{align*}
\begin{align*}
\alpha_2\circ g_*(a_1)&=\alpha_2(b_2^{-1}a_1b_1^{-1}a_1^{-1}b_2b_1b_2b_1^{-1}b_2^{-1}a_1b_1) = \gamma_1\\
\alpha_2\circ g_*(a_2)&=\alpha_2(b_1b_2^{-1}a_1b_1^{-1}a_1^{-1}b_2a_2b_1^{-1}b_2^{-1}a_1b_1a_1^{-1}b_2) = 1\\
\alpha_2\circ g_*(b_1)&=\alpha_2(b_1)=1\\
\alpha_2\circ g_*(b_2)&=\alpha_2(b_2^{-1}a_1b_1^{-1}a_1^{-1}b_2b_1b_2b_1^{-1}b_2^{-1}a_1b_1a_1^{-1}b_2)= \gamma_2
\end{align*}
\begin{align*}
\alpha_3\circ g_*(a_1)&=\alpha_3(b_2^{-1}a_1b_1^{-1}a_1^{-1}b_2b_1b_2b_1^{-1}b_2^{-1}a_1b_1) = 1\\
\alpha_3\circ g_*(a_2)&=\alpha_3(b_1b_2^{-1}a_1b_1^{-1}a_1^{-1}b_2a_2b_1^{-1}b_2^{-1}a_1b_1a_1^{-1}b_2) = \gamma_1\\
\alpha_3\circ g_*(b_1)&=\alpha_3(b_1) = \gamma_2\\
\alpha_3\circ g_*(b_2)&=\alpha_3(b_2^{-1}a_1b_1^{-1}a_1^{-1}b_2b_1b_2b_1^{-1}b_2^{-1}a_1b_1a_1^{-1}b_2) = 1
\end{align*}
\begin{align*}
\alpha_4\circ f_*(a_1)&= \alpha_4(a_1)=1\\
\alpha_4\circ f_*(a_2)&= \alpha_4(a_2)=1\\
\alpha_4\circ f_*(b_1)&= \alpha_4(b_1a_1^{-1})=\gamma_1\\
\alpha_4\circ f_*(b_2)&= \alpha_4(b_2a_2^{-1})=\gamma_2.
\end{align*}
So we have that $g$ lies in $\Aff^+(X)\cap\Mod_2(\alpha_1)$,
$\Aff^+(X)\cap\Mod_2(\alpha_2)$, $\Aff^+(X)\cap\Mod_2(\alpha_3)$ and
that $f$ lies in $\Aff^+(X)\cap\Mod_2(\alpha_4)$.\\
Unfortunately, we have $g\notin\Aff^+(X)\cap\Mod_2(\alpha_4)$ and
$f\notin\Aff^+(X)\cap\Mod_2(\alpha_i)$, $i\in\{1,2,3\}$. But for flat
surfaces we cannot exclude the possibility that two parabolic elements
with different eigendirections lie in such a group, which would then be
not cyclic.\\ For a flat surface $X$ of genus $g$ with a symplectic set
of generators $(a_1,b_1,\dots,a_g,b_g)$ for $\pi_1(X)$ and a symplectic
homomorphism $\alpha$ defined by $\alpha(a_i)=1$ and
$\alpha(b_i)=\gamma_i$, we can only find an element
$\varphi\in\Out^+(\pi_1(X))$ in $\Mod_g(\alpha)$ if we have standard
unit vectors in the first $g$ columns of the matrix $M_\varphi$, the
columns that correspond to elements from $\Kern(\alpha)$, see Lemma
\ref{EW1} d.\\ Elements that act trivially on the homology, as above
$g$, are, therefore, good candidates for elements in a group of the
form$\Aff^+(X)\cap\Mod_g(\alpha)$ since they satisfy this necessary
condition for all symplectic homomorphisms $\alpha$.
\printindex

\bibliographystyle{amsalpha}
\bibliography{Teichmueller_Referat_verbessert.bib}
\end{document}